\documentclass[11pt]{article}

\usepackage[utf8]{inputenc}
\usepackage[T1]{fontenc}
\usepackage{geometry}
\usepackage[german]{babel}
\usepackage{enumerate}
\usepackage[ampersand]{easylist}

\usepackage{amsmath}
\usepackage[lite]{amsrefs}
\usepackage{amssymb}
\usepackage{amsthm}

\makeatletter\newcommand{\MSC}[1]{\let\@dummy\@title\gdef\@title{\@dummy\footnotetext{#1}}} 

\theoremstyle{plain}
\newtheorem{theorem}{Satz}[section]
\newtheorem{lemma}[theorem]{Lemma}
\newtheorem{corollary}[theorem]{Korollar}
\newtheorem{definition}[theorem]{Definition}

\title{
Imagin\"arquadratische Einbettung von Ordnungen rationaler Quaternionenalgebren,
und die nichtzyklischen endlichen Untergruppen der Bianchi-Gruppen
}

\author{Norbert Kr\"amer }

\date{\today}

\MSC{
\emph{Mathematics subject classification (2010):}
\\
\hspace*{1 cm}
11S45 Algebras and orders, and their zeta functions
\\
\hspace*{1 cm}
11R52 Quaternion and other division algebras: arithmetic, zeta functions
\\
\hspace*{1 cm}
11F06 Structure of modular groups and generalizations; arithmetic groups
\\
\hspace*{1 cm}
20G07 Structure theory (Linear algebraic groups and related topics)
\\
\hspace*{1 cm}
20G30 Linear algebraic groups over global fields and their integers
 }

\begin{document}

\maketitle

\begin{abstract}

\noindent
Seien $k$ ein imagin\"arquadratischer Zahlk\"orper mit Hauptordnung $\mathfrak{o}$,
$F$ eine $\mathbb{Q}$-Qua\-ternionenalgebra und $M$ eine Erweiterung von $F$ zur $k$-Algebra.
Wir klassifizieren die $F$-Ordnungen $\mathfrak{G}$, die sich optimal in eine M-Maximalordnung $\mathfrak{N}$ einbetten lassen,
d. h. f\"ur die $\mathfrak{G} = F \cap \mathfrak{N}$ gilt (Definition \ref{Item05-03}, Satz \ref{Item05-06}).
Der Isomorphietyp von $\mathfrak{N}$ h\"angt nur von der Diskriminante von $\mathfrak{G}$ ab.
Ist $F' \subset M$ eine zweite $\mathbb{Q}$-Quaternionenalgebra und $\mathfrak{N}'$ eine zweite $M$-Maximalordnung,
so ist unser Hauptergebnis der Nachweis einer Beziehung zwischen den Einbettungen
$\mathfrak{G} = F \cap \mathfrak{N}$ und $\mathfrak{G}' = F' \cap \mathfrak{N}'$ (Satz \ref{Item05-09}).
\\[2pt]
Auf dieser Grundlage k\"onnen wir speziell ermitteln, ob die Bianchi-Gruppe $PSL_2(\mathfrak{o})$
$3$-Dieder-, Tetraeder- oder maximalendliche $2$-Diedergruppen enth\"alt (Satz \ref{Item06-08}) 
\\[2pt]
Wir ermitteln die Anzahl optimaler Einbettungen $\mathfrak{G} \hookrightarrow \mathfrak{N}$
bis auf Konjugation mit Elementen aus $\mathfrak{N}^\times$ (Satz \ref{Item07-03}).
Damit berechnen wir speziell die Konjugationsklassenzahlen
nichtzyklischer endlicher Untergruppen der Bianchi-Gruppe neu (Satz \ref{Item07-05}).
\\[2pt]
Weiter zeigen wir f\"ur $m = 2,3$, wie sich $m$-Diederuntergruppen in Gruppen der Ordnung $m$ schneiden
(S\"atze \ref{Item08-03} - \ref{Item08-05}).
Schlie{\ss}lich dehnen wir unsere Ergebnisse auf Eichler-Ordnungen und speziell Kongruenzuntergruppen aus
(S\"atze \ref{Item09-03} - \ref{Item09-09}).

\end{abstract}


\tableofcontents


\section{Einleitung} \label{Chapter01}

\noindent
Sei $d \in \mathbb{N}$ quadratfrei, und sei $k = \mathbb{Q}(i\sqrt{d})$ mit Hauptordnung $\mathfrak{o}$.
Sei $F$ eine $\mathbb{Q}$-Quater\-nionenalgebra, und sei $M \cong k \otimes_\mathbb{Q} F$ eine Erweiterung von $F$ zur $k$-Quaternionenalgebra.
\\[4pt]
Ist $\mathfrak{G}$ eine $F$-Ordnung mit Diskriminante $n^2 \mathbb{Z}$ und $n \in \mathbb{N}$,
dann sei $\Delta(\mathfrak{G}) := +n$ oder $-n$, je nachdem, ob $F$ an der Stelle $\infty$ zerf\"allt oder verzweigt ist.
Ist $\mathfrak{F}$ eine $F$-Maximalordnung mit $\mathfrak{G} \subset \mathfrak{F}$,
dann sei $\Lambda(\mathfrak{G}): = [\mathfrak{F} : \mathfrak{G}]$.
Wir bezeichnen $\mathfrak{G}$ als $\mathfrak{o}$-kompatibel,
wenn $\mathfrak{G}_p$ an allen Stellen $p$ von $\mathbb{Q}$ mit $p \mid \Lambda(\mathfrak{G})$
eine zu $\mathfrak{o}_p$ isomorphe Ordnung enth\"alt. 
\\[4pt]
Wir untersuchen die Existenz und die Struktur der $\mathfrak{o}$-kompatiblen $F$-Ordnungen
(Satz \ref{Item05-04} mit den S\"atzen \ref{Item04-02}, \ref{Item04-03} und \ref{Item04-05}),
und wir zeigen, dass eine $F$-Ordnung $\mathfrak{G}$ genau dann $\mathfrak{o}$-kompatibel ist,
wenn es eine $M$-Maximalordnung $\mathfrak{M}$ mit $\mathfrak{G} = F \cap \mathfrak{M}$ gibt,
und dass der Isomorphietyp von $\mathfrak{M}$ dann nur von $\Delta(\mathfrak{G})$ abh\"angt (Satz \ref{Item05-06}).
\\[4pt]
Unser Hauptergebnis ist der nachstehende Satz \ref{Item05-09},
der die Durchschnitte $F \cap \mathfrak{M}$ und $F' \cap \mathfrak{M}'$ zweier $\mathbb{Q}$-Quaternionenalgebren $F, F'$
mit zwei $M$-Maximalordnungen $\mathfrak{M}, \mathfrak{M}'$ zueinander in Beziehung setzt.
Wir verwenden darin die folgenden Notationen:
F\"ur eine Stelle $p$ von $\mathbb{Q}$ und $a, b \in \mathbb{Q}_p^\times$ sei $\left(\dfrac{a,b}{p}\right)$ das Hilbert-Symbol.
Sei $\left(\dfrac{F}{p}\right) := +1$ oder $-1$, je nachdem, ob $F$ an der Stelle p zerf\"allt oder verzweigt ist.
Sei $\Sigma_k(F)$ das Produkt der Verzweigungsstellen von $F$, die in $k$ zerlegt sind.
Und sei $\Lambda(\mathfrak{M} \cap \mathfrak{M}'): = [\mathfrak{M} : (\mathfrak{M} \cap \mathfrak{M}')]$.
\\[4pt]
$Bemerkung$: Der Isomorphietyp von $\mathfrak{M}$ relativ zu $\mathfrak{M}'$
h\"angt nur von der Gesamtheit der Hilbertsymbole $\left(\dfrac{\Lambda(\mathfrak{M} \cap \mathfrak{M}'),-d}{p}\right)$ ab.
$\mathfrak{M}$, $\mathfrak{M}'$ sind genau dann isomorph, wenn es ein $f \in \mathbb{N}$ gibt mit
$f \mid \Sigma_k(F)$ und $\left( \dfrac{f \Lambda(\mathfrak{M} \cap \mathfrak{M}'),-d} {p} \right) = 1$
f\"ur alle Stellen $p$ (Lemma \ref{Item05-05}).

\newtheorem*{theorem_Item05-09}{Satz \ref{Item05-09}} 
\begin{theorem_Item05-09} 

Sei $d \in \mathbb{N}$ quadratfrei, und sei $k = \mathbb{Q}(i\sqrt{d})$.
\\
Seien $F, F'$ zwei $\mathbb{Q}$-Quaternionenalgebren,
sei $M$ eine gemeinsame Erweiterung von $F$ und $F'$ zur $k$-Quaternionenalgebra,
und seien $\mathfrak{M}, \mathfrak{M}'$ zwei $M$-Maximalordnungen.
\\
Dann gibt es ein $f \in \mathbb{N}$ mit $f \mid \Sigma_k(F)$,
so dass f\"ur alle Stellen $p$ von $\mathbb{Q}$ gilt:
\\
$\left(\dfrac{\Delta(F \cap \mathfrak{M}) f \Lambda(\mathfrak{M} \cap \mathfrak{M}') \Delta(F' \cap \mathfrak{M}'),-d}{p}\right)
= \left(\dfrac{F}{p}\right) \left(\dfrac{F'}{p}\right)$.

\end{theorem_Item05-09} 

\noindent
Wir zeigen die oben genannten Aussagen \"uber die Einbettungen $\mathfrak{G} = F \cap \mathfrak{N}$
einschlie{\ss}lich Satz \ref{Item05-09}, indem wir erstens
je zwei $\mathbb{Q}$-Quaternionenalgebren $F, F' \subset M$ zueinander in Beziehung setzen (Satz \ref{Item03-03}),
zweitens  in Kapitel \ref{Chapter04} die Einbettungen
$\mathfrak{G}_p = F_p \cap \mathfrak{N}_p$ an den endlichen Stellen $p$ von $\mathbb{Q}$ untersuchen, an denen $M_p$ zerf\"allt,
und drittens in Kapitel \ref{Chapter05} diese lokalen Ergebnisse zusammenf\"ugen und mit Satz \ref{Item03-03} kombinieren.
Falls $k$ Zerf\"allungsk\"orper von $F$ ist,
folgen aus den S\"atzen \ref{Item03-03} und \ref{Item05-09} die Korollare \ref{Item03-04} und \ref{Item05-10}.
\\[4pt]
Als Anwendung von Korollar \ref{Item05-10} k\"onnen wir dann unter anderem kl\"aren, welche nichtzyklischen endlichen Untergruppen die Bianchi-Gruppe $PSL_2(\mathfrak{o}) = SL_2(\mathfrak{o})/\{ \pm 1 \}$ enth\"alt:

\newtheorem*{theorem_Item06-08}{Satz \ref{Item06-08}} 
\begin{theorem_Item06-08} 

Sei $d \in \mathbb{N}$ quadratfrei, und sei $k = \mathbb{Q}(i\sqrt{d})$ mit Hauptordnung $\mathfrak{o}$. Dann gilt:
\begin{enumerate}[(i)]

\item [(\ref{Item06-08.1})]
$PSL_2(\mathfrak{o})$ enth\"alt genau dann eine $3$-Diedergruppe $\mathcal{D}_3$,
\\
wenn $p \equiv 1 \pmod{3}$ f\"ur alle Primteiler $p \not= 3$ von $d$.

\item [(\ref{Item06-08.2})]
$PSL_2(\mathfrak{o})$ enth\"alt genau dann eine Tetraedergruppe $\mathcal{T}$,
\\
wenn $p \equiv 1$ oder $p \equiv 3 \pmod{8}$ f\"ur alle Primteiler $p \not= 2$ von $d$.

\item [(\ref{Item06-08.3})]
$PSL_2(\mathfrak{o})$ enth\"alt genau dann eine maximalendliche $2$-Diedergruppe $\mathcal{D}_2$,
\\
wenn $p \equiv 1 \pmod{4}$ f\"ur alle Primteiler $p \not= 2$ von $d$.

\end{enumerate}

\end{theorem_Item06-08}

\noindent
Wir f\"uhren Satz \ref{Item06-08} wie folgt auf Satz \ref{Item05-09} zur\"uck:
Enth\"alt $PSL_2(k)$ eine Untergruppe vom Isomorphietyp $\mathcal{D}_3$, $\mathcal{T}$ oder $\mathcal{D}_2$,
dann erzeugt deren Urbild in $SL_2(k)$
eine $\mathbb{Q}$-Qua\-ternionenalgebra $F \subset M_2(k)$ und \"uber $\mathbb{Z}$ eine $F$-Ordnung $\mathfrak{F}$.
Wir wenden Satz \ref{Item05-09} auf die Einbettungen $M_2(\mathbb{Z}) = M_2(\mathbb{Q}) \cap M_2(\mathfrak{o})$
und $\mathfrak{F} = F \cap \mathfrak{N}$ mit einer $M_2(k)$-Maximalordnung $\mathfrak{N}$ an.
Wir erhalten die S\"atze \ref{Item06-04} und \ref{Item06-06}, sowie Satz \ref{Item06-07} mit dem Spezialfall Satz \ref{Item06-08}.
\\[4pt]
F\"ur eine $\mathfrak{o}$-kompatible $F$-Ordnung $\mathfrak{G}$ ermitteln wir zus\"atzlich
die Anzahl der $M$-Maxi\-malordnungen $\mathfrak{N}$ mit $\mathfrak{G} = F \cap \mathfrak{N}$
und die Anzahl der $\mathfrak{N}^\times$-Konjugationsklassen von Einbettungen
$j: \mathfrak{G} \hookrightarrow \mathfrak{N}$ mit $j(\mathfrak{G}) = j(F) \cap \mathfrak{N}$
(Satz \ref{Item07-03}, auf Basis des lokalen Satzes \ref{Item04-09}).
\\[4pt]
Als Anwendung von Satz \ref{Item07-03} berechnen wir dann unter anderem die Konjugationsklassenanzahlen
der nichtzyklischen maximalendlichen Untergruppen von $PSL_2(\mathfrak{o})$ (Satz \ref{Item07-05}).
\\[4pt]
Unter Verwendung der Ergebnisse von \cite{Kr} lassen sich dann auch unter anderem die Konjugationsklassenanzahlen der zyklischen Untergruppen einer Bianchi-Gruppe angeben.
\\
Rahm hat in \cite{Ra} die homologische Torsion und die Farrell-Tate-Kohomologie der Bianchi-Gruppen ermittelt und als Funktion dieser Konjugationsklassenanzahlen ausgedr\"uckt.
\\[4pt]
In Kapitel \ref{Chapter08} untersuchen wir die Durchschnitte nichtzyklischer endlicher Untergruppen (wird fortgesetzt).
F\"urs Erste zeigen wir f\"ur $m = 2,3$, ob und wie sich $m$-Diedergruppen in Gruppen der Ordnung $m$ schneiden
(S\"atze \ref{Item08-03} - \ref{Item08-05} und Korollar \ref{Item08-06}).
\\
Berkove und Rahm haben gezeigt, dass daraus
eine Klassifikation der Zusammenhangskomponenten reduzierter 2-Torsionsunterkomplexe folgt,
was es erlaubt, Formeln f\"ur die 2-Torsionskohomologie aufzustellen, siehe \cite[Corollary 7]{Be}.
\\[4pt]
Im Schlusskapitel \ref{Chapter09} dehnen wir die Ergebnisse der Kapitel \ref{Chapter06} bis \ref{Chapter08}
auf Eichler-Ordnungen und speziell auf die Bianchi-Kongruenzuntergruppen aus
(S\"atze \ref{Item09-03} - \ref{Item09-09}).
Nur relativ wenige Kongruenzuntergruppen enthalten maximalendliche nichtzyklische Untergruppen,
und wenn, dann jeweils nur Untergruppen eines einzigen Isomorphietyps.
\\[4pt]
Ich danke Alexander D. Rahm f\"ur seine Ratschl\"age und die technische Unterst\"utzung,
und J\"urgen Rohlfs f\"ur die kritische Durchsicht des Manuskripts und hilfreiche Vorschl\"age.


\section{Bezeichnungen und Grundlagen} \label{Chapter02}

F\"ur die in diesem Artikel benutzen zahlentheoretischen Grundlagen verweisen wir summarisch auf
\cite{Bo} f\"ur imagin\"arquadratische Zahlk\"orper und \cite{Ma} f\"ur Quaternionenalgebren.
F\"ur die speziellen Bezeichnungen und Grundlagen zu Bianchi-Gruppen siehe Kapitel \ref{Chapter06}.
\\[4pt]
In einer abelschen Gruppe $G$ sei $G^{(2)}$ die Untergruppe der Quadrate.
\\
In einem Ring $R$ mit Eins sei $R^\times$ die Gruppe der multiplikativ invertierbaren Elemente.
\\
F\"ur kommutatives $R$ sei $M_2(R)$ die $R$-Algebra der $2$$\times$$2$-Matrizen mit Koeffizienten in $R$.
\\
F\"ur $z \in \mathbb{C}$ bezeichne $z \mapsto \overline{z}$ die komplexe Konjugation und $|z| = \sqrt{z\overline{z}}$ den Absolutbetrag.
\\[4pt]
Sei im Folgenden stets $d \in \mathbb{N}$ quadratfrei
und $k = \mathbb{Q}(i\sqrt{d}) \subset \mathbb{C}$ ein imagin\"arquadratischer Zahlk\"orper mit Hauptordnung $\mathfrak{o}$ und Diskriminante $D$.
Die Einschr\"ankung der komplexen Konjugation auf $k$ stimmt mit der nichttrivialen Galois-Involution von $k/\mathbb{Q}$ \"uberein.
\\
Bezeichne $I$ die Gruppe der $\mathfrak{o}$-Ideale und $H \subset I$ die Untergruppe der Hauptideale.
\\[4pt]
Sei $p$ stets eine Stelle von $\mathbb{Q}$, d.h. eine Primzahl $p \in \mathbb{N}$ oder die unendliche Stelle $p = \infty$.
Sei $\mathfrak{p}$ stets eine Stelle von $k$, d.h.
ein Primideal $\mathfrak{p} \subset \mathfrak{o}$ oder die unendliche Stelle $\mathfrak{p} = \infty$.
\\[4pt]
Ist $p$ eine endliche Stelle von $\mathbb{Q}$, dann ist $k$ auf nat\"urliche Weise eingebettet in die halbeinfache quadratische Erweiterung $k_p \cong \mathbb{Q}_p \otimes_\mathbb{Q}k$ von $\mathbb{Q}_p$, mit der Hauptordnung $\mathfrak{o}_p = \mathbb{Z}_p\mathfrak{o}$.
Die Galois-Involution $a \mapsto \overline{a}$ von $k/\mathbb{Q}$ setzt sich auf nat\"urliche Weise auf $k_p/\mathbb{Q}_p$ fort.
\\
Ist $p$ verzweigt oder tr\"age in $k$ mit Primteiler $\mathfrak{p}$,
dann ist auf nat\"urliche Weise $k_p = k_\mathfrak{p}$ und $\mathfrak{o}_p = \mathfrak{o}_\mathfrak{p}$.
F\"ur die jeweils einzige unendliche Stelle von $\mathbb{Q}$ und $k$ gilt: $\mathbb{Q}_\infty = \mathbb{R}$, $k_\infty = \mathbb{C}$.
\\
Ist $p$ zerlegt in $k$ mit den Primteilern $\mathfrak{p}$ und $\overline{\mathfrak{p}}$,
dann ist $k_p \cong k_\mathfrak{p} \times k_{\overline{\mathfrak{p}}}$
und $\mathfrak{o}_p \cong \mathfrak{o}_\mathfrak{p} \times \mathfrak{o}_{\overline{\mathfrak{p}}}$.
\\[4pt]
F\"ur eine Stelle $p$ von $\mathbb{Q}$ und $a, b \in \mathbb{Q}_p^\times$ sei $\left(\dfrac{a,b}{p}\right)$ das Hilbert-Symbol.
\\
F\"ur Rechenregeln zum Hilbert-Symbol siehe \cite[Kapitel I, § 6, Formeln (9)-(13), Satz 7]{Bo}.
\\[4pt]
Eine $\mathbb{Q}$-Quaternionenalgebra $F$ ist auf nat\"urliche Weise eingebettet in ihre $p$-Komponente $F_p \cong \mathbb{Q}_p \otimes_\mathbb{Q} F$.
Wir identifizieren $\mathbb{Q}$ und $\mathbb{Q}_p$ mit den Zentren von $F$ und $F_p$.
F\"ur $p \in \mathbb{N}$ und einen $\mathbb{Z}$-Modul $\mathfrak{F} \subset F$
ist $\mathfrak{F}_p = \mathbb{Z}_p\mathfrak{F}$ die $p$-Komponente von $\mathfrak{F}$.
\\
F\"ur eine $k$-Quaternionenalgebra $M$ ist analog $k \subset M \subset M_\mathfrak{p} \cong k_\mathfrak{p} \otimes_k M$
und $k_\mathfrak{p} \subset M_\mathfrak{p}$.
Ist $\mathfrak{p}$ endlich und $\mathfrak{M} \subset M$ ein $\mathfrak{o}$-Modul,
so ist $\mathfrak{M}_\mathfrak{p} = \mathfrak{o}_\mathfrak{p}\mathfrak{M}$.
Wir identifizieren $M_\infty$ mit $M_2(\mathbb{C})$.
\\[4pt]
Sind $p \in \mathbb{N}$ eine Primzahl, $F$ eine $\mathbb{Q}$-Quaternionenalgebra
und $M \cong k \otimes_\mathbb{Q} F$ eine Erweiterung von $F$ zur $k$-Quaternionenalgebra,
dann bezeichne $M_p \cong k_p \otimes_k M \cong k_p \otimes_\mathbb{Q} F$ die $p$-Komponente von $M$.
Auf nat\"urliche Weise ist $M \subset M_p$ und $k_p$ das Zentrum von $M_p$.
\\[4pt]
Ist $p$ verzweigt oder tr\"age in $k$ mit Primteiler $\mathfrak{p}$, so ist auf nat\"urliche Weise $M_p = M_\mathfrak{p}$,
und f\"ur einen $\mathfrak{o}$-Modul $\mathfrak{M} \subset M$ ist $\mathfrak{M}_p = \mathfrak{M}_\mathfrak{p}$.
Ist $p$ zerlegt in $k$ mit den Primteilern $\mathfrak{p}$ und $\overline{\mathfrak{p}}$,
so ist $M_p \cong M_\mathfrak{p} \times M_{\overline{\mathfrak{p}}} \cong F_p \times F_p$.
Wir nennen dann $M_p$ eine $k_p$-Quaternionenalgebra, obwohl $k_p \cong k_\mathfrak{p} \times k_{\overline{\mathfrak{p}}}$ kein K\"orper ist.
Ist $\mathfrak{M}$ eine $M$-Ordnung, so nennen wir $\mathfrak{M}_p = \mathfrak{o}_p \mathfrak{M} \cong \mathfrak{M}_\mathfrak{p} \times \mathfrak{M}_{\overline{\mathfrak{p}}}$ eine $M_p$-Ordnung.
Diese Benennungen erm\"oglichen uns einheitliche und unkomplizierte Formulierungen von S\"atzen und Lemmata ab Kapitel \ref{Chapter04}.
\\[4pt]
Eine $F$-Ordnung $\mathfrak{G}$ hei{\ss}t optimal eingebettet in eine $M$-Ordnung $\mathfrak{M}$,
wenn $\mathfrak{G} = F \cap \mathfrak{M}$ gilt.
Eine $F$-Maximalordnung $\mathfrak{F}\subset \mathfrak{M}$ ist stets optimal in $\mathfrak{M}$ eingebettet.
\\[4pt]
F\"ur ein Element $A$ einer Quaternionenalgebra $Q$ bezeichne $A^*$ das konjugierte Element.
\\
$N(A) = AA^*$ und $S(A) = A + A^*$ sind dann die (reduzierte) Norm und die Spur von $A$.
\\[4pt]
Ist $K$ ein K\"orper und $Q \cong M_2(K)$ eine K-Quaternionenalgebra,
dann hei{\ss}en, in dieser Reihenfolge, $U_{11}, U_{12}, U_{21}, U_{22} \in Q$ Matrixeinheiten in $Q$,
wenn f\"ur alle $r, s, \in \{1, 2 \}$ gilt: $U_{r1}U_{1s} = U_{r2}U_{2s} = U_{rs}  \ne 0$ und $U_{r1}U_{2s} = U_{r2}U_{1s} = 0$.
Dann ist $\{ U_{11}, U_{12}, U_{21}, U_{22} \}$ eine Basis von $Q$ \"uber $K$,
und f\"ur ein $A = \alpha U_{11} + \beta U_{12} + \gamma U_{21} + \delta U_{22} \in Q$ mit $\alpha, \beta, \gamma, \delta \in K$
gilt bez\"uglich dieser Matrixeinheiten die Gleichung (oder Darstellung) $A = \begin{pmatrix} \alpha & \beta \\ \gamma & \delta \end{pmatrix}$.
F\"ur eine Darstellung bez\"uglich Matrixeinheiten gelten die bekannten Matrix-Rechenregeln.
\\[4pt]
Ist $K$ ein algebraischer oder $\mathfrak{p}$-adischer Zahlk\"orper mit Hauptordnung $\mathfrak{O}$, und ist $\mathfrak{A}$ ein $\mathfrak{O}$-Ideal,
dann bezeichne $N_{abs}(\mathfrak{A}) \in \mathbb{Q}^+$ die Absolutnorm von $\mathfrak{A}$,
und f\"ur $A \in K^\times$ sei $N_{abs}(A) := N_{abs}(A\mathfrak{O})$.
Falls $\mathfrak{A} \subset \mathfrak{O}$, ist $N_{abs}(\mathfrak{A}) = [\mathfrak{O} : \mathfrak{A}]$.

\begin{lemma} \label{Item02-01}

Sei $K$ algebraischer oder $\mathfrak{p}$-adischer Zahlk\"orper, und sei $Q$ eine K-Quater\-nionenalgebra.
Seien $A, B \in Q$ mit $N(A) = N(B)$, $S(A) = S(B)$ und $S(A)^2 \not= 4N(A)$.
\\
Dann gibt es ein $J \in Q^\times$ mit $JAJ^{-1} = B$.

\end{lemma}

\begin{proof}

Falls $K[A]$ K\"orper ist, siehe \cite[Theorem 2.9.8]{Ma}.
Sonst k\"onnen wir $Q = M_2(K)$ annehmen.
Das charakteristische Polynom von $A$ hat in $K$ zwei Nullstellen $a \not= b$.
Nun folgt, dass es $X, Y \in GL_2(K)$ gibt mit
$XAX^{-1} = YBY^{-1} = \begin{pmatrix} a & 0 \\ 0 & b \end{pmatrix}$.
Sei $J := Y^{-1}X$.

\end{proof}

\begin{lemma} \label{Item02-02}

Sei $K$ ein $\mathfrak{p}$-adischer Zahlk\"orper mit Hauptordnung $\mathfrak{O}$ und Primelement $\Pi$.
Seien $Q \cong M_2(K)$ eine K-Quaternionenalgebra und $\mathfrak{M}, \mathfrak{M}'$ zwei $Q$-Maximalordnungen.
\\
Dann gibt es genau ein $r \in \mathbb{N}_0$ mit den folgenden Eigenschaften.
\begin{enumerate}[(i)]

\item \label{Item02-02.1}
Es gibt Matrixeinheiten in $Q$, bez\"uglich derer gilt:
\\
$\mathfrak{M} = \begin{pmatrix} \mathfrak{O} & \mathfrak{O} \\ \mathfrak{O} & \mathfrak{O} \end{pmatrix}$ und
$\mathfrak{M}' = \begin{pmatrix} \mathfrak{O} & \Pi^{-r}\mathfrak{O} \\ \Pi^r\mathfrak{O} & \mathfrak{O} \end{pmatrix}
= J \mathfrak{M} J^{-1}$ mit $J = \begin{pmatrix} 1 & 0 \\ 0 & \Pi^r \end{pmatrix}$.

\item \label{Item02-02.2}
$N(\mathfrak{M}'\mathfrak{M})^{-1} = N(\mathfrak{M}\mathfrak{M}')^{-1} = \Pi^r\mathfrak{O}$.

\item \label{Item02-02.3}
Sei $J' \in \mathfrak{M} \smallsetminus \Pi\mathfrak{M}$ mit $\mathfrak{M}' = J'\mathfrak{M}J'^{-1}$. Dann ist $N(J') \in \Pi^r\mathfrak{O}^\times$.

\item \label{Item02-02.4}
$\mathfrak{M} \cap \mathfrak{M}'$ hat die Diskriminante $\Pi^{2r}\mathfrak{O}$, und es gilt
$[\mathfrak{M} : (\mathfrak{M} \cap \mathfrak{M}')] = N_{abs}( \Pi\mathfrak{O})^r$.

\end{enumerate}

\end{lemma}

\begin{proof}

Wir zeigen die Existenz von $r$. Die Eindeutigkeit folgt dann sofort, etwa mit (\ref{Item02-02.2}).
\begin{enumerate}[(i)]

\item[(\ref{Item02-02.1})]
F\"ur einen kurzen Beweis siehe \cite[Beweis von Satz 7, Fall b)]{Ei}.

\item[(\ref{Item02-02.2})]
Man pr\"uft leicht nach, dass $\mathfrak{M}\mathfrak{M}' = \mathfrak{M}J^{-1}$.
Also gilt $N(\mathfrak{M}\mathfrak{M}')^{-1} = N(J)\mathfrak{O} = \Pi^r\mathfrak{O}$.

\item[(\ref{Item02-02.3})]
Wegen $\mathfrak{M}' = J'\mathfrak{M}J'^{-1} = J\mathfrak{M}J^{-1}$
ist $J^{-1}J'\mathfrak{M} = \mathfrak{M}J^{-1}J'$ ein zweiseitiges $\mathfrak{M}$-Ideal.
Mit \cite[Theorem 6.5.3.2]{Ma} folgt $J' = \Pi^eJY$ f\"ur ein $e \in \mathbb{Z}$ und ein $Y \in \mathfrak{M}^\times$.
\\
Wegen $J, J' \in \mathfrak{M}$ und $J, J' \not\in \Pi\mathfrak{M}$ ist $e = 0$.
Also gilt $N(J')\mathfrak{O} = N(J)\mathfrak{O} = \Pi^r\mathfrak{O}$.

\item[(\ref{Item02-02.4})]
Nach (\ref{Item02-02.1}) hat $\mathfrak{M} \cap \mathfrak{M}'$ die Diskriminante $\Pi^{2r}\mathfrak{O}$,
und es ist $[\mathfrak{M} : (\mathfrak{M} \cap \mathfrak{M}')] = [\mathfrak{O} : \Pi^r\mathfrak{O}]$.

\end{enumerate}

\end{proof}

\begin{lemma} \label{Item02-03}

Sei $K$ ein $\mathfrak{p}$-adischer Zahlk\"orper, und sei $Q$ eine K-Quaternionenalgebra.
Sei $L \subset Q$ eine halbeinfache quadratische Erweiterung von $K$,
und seien $\mathfrak{M}, \mathfrak{M}'$ zwei $Q$-Maximalordnungen mit $L \cap \mathfrak{M} = L \cap \mathfrak{M}'$.
Dann gibt es ein $J \in L^\times$ mit $\mathfrak{M}' = J\mathfrak{M}J^{-1}$.

\end{lemma}

\begin{proof}

Spezialfall einer lokalen Version von \cite[Satz 7]{Ei}, Beweis siehe dort, Fall b).

\end{proof}

\begin{lemma} \label{Item02-04}

Sei $K$ ein $\mathfrak{p}$-adischer Zahlk\"orper,
und sei $Q \cong M_2(K)$ eine K-Quaternionen\-algebra.
Seien $\mathfrak{M} \ne \mathfrak{M}'$ sowie $\mathfrak{N} \ne \mathfrak{N}'$ jeweils $Q$-Maximalordnungen
mit $\mathfrak{M} \cap \mathfrak{M}' = \mathfrak{N} \cap \mathfrak{N}'$.
Dann gilt entweder $\mathfrak{M} = \mathfrak{N}$ und $\mathfrak{M}' = \mathfrak{N}'$
oder $\mathfrak{M} = \mathfrak{N}'$ und $\mathfrak{M}' = \mathfrak{N}$.

\end{lemma}

\begin{proof}

Sei $\mathfrak{O}$ die Hauptordnung von $K$, und sei $\Pi \in \mathfrak{O}$ ein Primelement.
Nach Lemma \ref{Item02-02}.(\ref{Item02-02.1}) gibt es ein $r \in \mathbb{N}$
und Matrixeinheiten in $Q$, bez\"uglich derer gilt:
$\mathfrak{M} = \begin{pmatrix} \mathfrak{O} & \mathfrak{O} \\ \mathfrak{O} & \mathfrak{O} \end{pmatrix}$,
$\mathfrak{M}' = \begin{pmatrix} \mathfrak{O} & \Pi^{-r}\mathfrak{O} \\ \Pi^r\mathfrak{O} & \mathfrak{O} \end{pmatrix}$
und $\mathfrak{M} \cap \mathfrak{M}' = \begin{pmatrix} \mathfrak{O} & \mathfrak{O} \\ \Pi^r\mathfrak{O} & \mathfrak{O} \end{pmatrix}$.
Dann ist $L := \begin{pmatrix} K & 0 \\ 0 & K \end{pmatrix}$ halbeinfache quadratische Erweiterung von $K$
mit Hauptordnung $\mathfrak{L} := \begin{pmatrix} \mathfrak{O} & 0 \\ 0 & \mathfrak{O} \end{pmatrix}
\subset \mathfrak{M} \cap \mathfrak{M}' = \mathfrak{N} \cap \mathfrak{N}'$.
Speziell gilt $\mathfrak{L} \subset L \cap \mathfrak{N}$.
Da $\mathfrak{L}$ die Hauptordnung von $L$ ist, folgt $L \cap \mathfrak{N} = \mathfrak{L} = L \cap \mathfrak{M}$.
Nach Lemma \ref{Item02-03} gibt es ein $J \in L^\times$ mit $\mathfrak{N} = J\mathfrak{M}J^{-1}$.
Wir k\"onnen o.B.d.A. annehmen, dass $J = \begin{pmatrix} 1 & 0 \\ 0 & \Pi^s \end{pmatrix}$ mit $s \in \mathbb{Z}$,
also dass $\mathfrak{N} = \begin{pmatrix} \mathfrak{O} & \Pi^{-s}\mathfrak{O} \\ \Pi^s\mathfrak{O} & \mathfrak{O} \end{pmatrix}$ gilt.
Aus $\mathfrak{M} \cap \mathfrak{M}' \subset \mathfrak{N}$ folgt $0 \le s \le r$.
Genauso gibt es ein $s' \in \mathbb{Z}$ mit $0 \le s' \le r$ und
$\mathfrak{N}' = \begin{pmatrix} \mathfrak{O} & \Pi^{-s'}\mathfrak{O} \\ \Pi^{s'}\mathfrak{O} & \mathfrak{O} \end{pmatrix}$.
Dann ist $\mathfrak{N} \cap \mathfrak{N}' =
\begin{pmatrix} \mathfrak{O} & \Pi^{-t'}\mathfrak{O} \\ \Pi^t\mathfrak{O} & \mathfrak{O} \end{pmatrix}$
mit $t = \max(s,s') = r$ und $t' = \min(s,s') = 0$,
also entweder $s=r$ und $s'=0$ oder $s=0$ und $s'=r$, was zu beweisen war.

\end{proof}


\section{Einbettung rationaler Quaternionenalgebren} \label{Chapter03}

\begin{definition} \label{Item03-01}

Sei $d \in \mathbb{N}$ quadratfrei, und sei $k = \mathbb{Q}(i\sqrt{d})$.
Sei $F$ eine $\mathbb{Q}$-Quaternionen\-algebra, und sei $M$ eine Erweiterung von $F$ zur $k$-Quaternionenalgebra.
\begin{enumerate}[(i)]

\item \label{Item03-01.1}
F\"ur eine Stelle $p$ von $\mathbb{Q}$ sei $\left(\dfrac{F}{p}\right) := +1$ oder $-1$,
\\
je nachdem, ob $F$ an der Stelle p zerf\"allt oder verzweigt ist.

\item \label{Item03-01.2}
Sei $\Sigma(F)$ das Produkt der endlichen Verzweigungsstellen von $F$,
\\
multipliziert mit $-1$, falls $F$ an der Stelle $\infty$ verzweigt ist.

\item \label{Item03-01.3}
Sei $\Sigma_k(F)$ das Produkt der Verzweigungsstellen von $F$, die in $k$ zerlegt sind.

\item \label{Item03-01.4}
Sei $\Phi_F: M \rightarrow M$ die Abbildung mit $\Phi_F(A + i\sqrt{d}B): = A -  i\sqrt{d}B$ f\"ur $A, B \in F$.

\item \label{Item03-01.5}
F\"ur $T \in M^\times$ mit $\Phi_F(T) = T^*$ sei $F(F,T): = \{ A \in M \mid T \Phi_F(A)T^{-1} = A \}$.

\end{enumerate}

\end{definition}

\begin{lemma} \label{Item03-02}

Sei $d \in \mathbb{N}$ quadratfrei, und sei $k = \mathbb{Q}(i\sqrt{d})$.
Sei $F$ eine $\mathbb{Q}$-Quaternionen\-algebra, und sei $M$ eine Erweiterung von $F$ zur $k$-Quaternionenalgebra.
Dann gilt:
\begin{enumerate}[(i)]

\item \label{Item03-02.1}
$M$ ist ist an einer Stelle $\mathfrak{p}$ von $k$ genau dann verzweigt,
wenn $\Sigma_k(F) \in \mathfrak{p}$ gilt.

\item \label{Item03-02.2}
Eine $\mathbb{Q}$-Quaternionenalgebra $E$ kann genau dann in $M$ eingebettet werden,
\\
wenn $\Sigma_k(E) = \Sigma_k(F)$ gilt.

\item \label{Item03-02.3}
Es gibt bis auf Isomorphie genau eine $\mathbb{Q}$-Quaternionenalgebra $E$
\\
mit $|\Sigma(E)| = \Sigma_k(E) = \Sigma_k(F)$.

\end{enumerate}

\end{lemma}

\begin{proof} ${}$
\begin{enumerate}[(i)]

\item[(\ref{Item03-02.1})]
Ist $p \ne \infty$ in $k$ nicht zerlegt, dann zerf\"allt $M_p$,
siehe \cite[Theorem 2.6.5]{Ma}, und es gilt $M_\infty = M_2(\mathbb{C})$.
Ist $p \ne \infty$ in $k$ zerlegt, $p \mathfrak{o} = \mathfrak{p} \overline{\mathfrak{p}}$,
dann gilt $M_p \cong M_\mathfrak{p} \times M_{\overline{\mathfrak{p}}} \cong F_p \times F_p$.

\item[(\ref{Item03-02.2})]
folgt aus (\ref{Item03-02.1}), siehe \cite[Theorem 2.7.5]{Ma}.

\item[(\ref{Item03-02.3})]
Sei $S$ die Menge der Stellen $p$ von $\mathbb{Q}$ mit $p \mid \Sigma_k(F)$.
Dann hat $S$ oder $S \cup \{ \infty \}$ eine gerade Anzahl von Elementen.
Die Behauptung folgt mit \cite[Theorem 7.3.6]{Ma}.

\end{enumerate}

\end{proof}

\begin{theorem} \label{Item03-03}

Sei $d \in \mathbb{N}$ quadratfrei, und sei $k = \mathbb{Q}(i\sqrt{d})$.
Sei $E$ eine $\mathbb{Q}$-Quaternionenalgebra, und sei $M$ eine Erweiterung von $E$ zur $k$-Quaternionenalgebra.
Dann gilt:

\begin{enumerate}[(i)]

\item \label{Item03-03.1}
$\Phi_E$ ist $\mathbb{Q}$-Algebrenautomorphismus von $M$
mit $\Phi_E(\lambda A) = \overline{\lambda} \Phi_E(A)$ und $\Phi_E^2(A) = A$ und $\Phi_E(A^*) = \Phi_E(A)^*$ f\"ur alle $\lambda \in k$, $A \in M$.
Es gilt $E = \{A \in M \mid \Phi_E(A) = A\}$.

\item \label{Item03-03.2}
Sei $T \in M^\times$ mit $\Phi_E(T) = T^*$. Dann gilt:
\\
$F(E,T)$ ist $\mathbb{Q}$-Quaternionenalgebra mit $\Phi_{F(E,T)}(A) = T \Phi_E(A) T^{-1}$ f\"ur $A \in M$,
und $N(T) \in \mathbb{Q}^\times$
sowie $\left(\dfrac{N(T),-d}{p}\right) = \left(\dfrac{E}{p}\right) \left(\dfrac{F(E,T)}{p}\right)$ f\"ur alle Stellen $p$ von $\mathbb{Q}$.

\item \label{Item03-03.3}
Sei $F \subset M$ eine $\mathbb{Q}$-Quaternionenalgebra.
Dann gibt es ein, bis auf einen Faktor aus $\mathbb{Q}^\times$ eindeutiges, $T \in M^\times$ mit $\Phi_E(T) = T^*$ und $F = F(E,T)$.
Falls $S(T)^2 \ne 4N(T)$, ist $E \cap F = \mathbb{Q}[i\sqrt{d}(T-T^*)]$ eine halbeinfache quadratische Erweiterung von $\mathbb{Q}$.

\end{enumerate}

\end{theorem}

\noindent
$Bemerkung$:
F\"ur jede Stelle $p$ von $\mathbb{Q}$ ist $\Phi_E$ auf nat\"urliche Weise zu einem
$\mathbb{Q}_p$-Algebren\-automorphismus von $M_p$ mit Fixalgebra $E_p$ fortsetzbar,
den wir auch mit $\Phi_E$ bezeichnen.

\begin{proof}

Wir k\"urzen $\Phi: = \Phi_E$ und $F(T) := F(E,T)$ ab.
\begin{enumerate}[(i)]

\item[(\ref{Item03-03.1})]
rechnet man leicht nach.

\item[(\ref{Item03-03.2})]
$F(T)$ ist eine $\mathbb{Q}$-Algebra und via Zuordnungsvorschrift $A \mapsto i\sqrt{d}A$
als $\mathbb{Q}$-Vektorraum isomorph zu $F(T,-) := \{ A \in M \mid T \Phi(A)T^{-1} = -A \}$.
Es gilt $T\Phi(T) \in \mathbb{Q}^\times$,
und mithilfe der Zerlegung $A = (A + T \Phi(A) T^{-1})/2 + (A - T \Phi(A) T^{-1})/2$ folgt, dass $M = F(T) \oplus F(T,-)$ gilt.
$F(T)$ hat \"uber $\mathbb{Q}$ die Dimension $4$.
Wegen $i\sqrt{d} \notin F(T)$ ist $F(T)$ eine zentrale $\mathbb{Q}$-Algebra,
und wegen $M \cong k \otimes_\mathbb{Q} F(T)$ auch einfache $\mathbb{Q}$-Algebra.
Sei $A \in M$. Dann gibt es $A', A'' \in F(T)$ mit $A = A' + i\sqrt{d}A''$, und damit gilt:
$\Phi_{F(T)}(A) = A' - i\sqrt{d}A'' = T\Phi(A')T^{-1} - i\sqrt{d}T\Phi(A'')T^{-1} = T\Phi(A' + i\sqrt{d}A'')T^{-1}$.

Sei $p$ eine Stelle von $\mathbb{Q}$.
Ist $p$ in $k$ zerlegt, so ist $\left(\dfrac{N(T),-d}{p}\right) = 1 = \left(\dfrac{E}{p}\right) \left(\dfrac{F(T)}{p}\right)$, siehe Lemma \ref{Item03-02}.(\ref{Item03-02.2}).
Wir k\"onnen nun also $p$ als nicht zerlegt in $k$ annehmen.
\begin{itemize}

\item
Im Fall $\left(\dfrac{E}{p}\right) \left(\dfrac{F(T)}{p}\right) = 1$ ist $E_p$ isomorph zu $F(T)_p$.
\\
Einen Isomorphismus $E_p \rightarrow F(T)_p$ k\"onnen wir zu einem $k_p$-Algebrenautomor\-phismus von $M_p$ fortsetzen.
Daher gibt es ein $J \in M_p^\times$ mit $F(T)_p = JE_pJ^{-1}$.
Dann gilt $JAJ^{-1} = \Phi_{F(T)}(JAJ^{-1}) = T\Phi(JAJ^{-1})T^{-1} = T\Phi(J)A\Phi(J)^{-1}T^{-1}$ f\"ur $A \in E_p$ 
Also gibt es $\lambda \in k_p^\times$ mit $T = \lambda J \Phi(J)^{-1}$ und $N(T) = T\Phi(T) = \lambda \overline{\lambda}$.

\end{itemize}

Sei $\epsilon \in \mathbb{Z}$ teilerfremd zu $D\Sigma(E)$ mit $\tau := \epsilon\Sigma(E) < 0$ und $\left(\dfrac{\tau,-d}{p}\right) = -1$.
\\
An allen Verzweigungsstellen $q$ von $E$ ist $\mathbb{Q}(\sqrt{\tau})_q$ quadratischer Erweiterungsk\"orper von $\mathbb{Q}_q$,
also ist $\mathbb{Q}(\sqrt{\tau})$ Zerf\"allungsk\"orper von $E$, siehe \cite[Theorems 2.6.5 and 2.7.3]{Ma}.
\\
Es gibt ein $U \in E^\times$ mit $U^2 = \tau$. Sei $T' := i\sqrt{d}U$. Dann gilt $\Phi(T') = -T' = T'^*$.
Aus $N(T') = d\tau$ folgt $\left(\dfrac{N(T'),-d}{p}\right) = -1$ und notwendig $\left(\dfrac{E}{p}\right) \left(\dfrac{F(T')}{p}\right) = -1$.
\begin{itemize}

\item
Im Fall $\left(\dfrac{E}{p}\right) \left(\dfrac{F(T)}{p}\right) = -1$ ist $F(T')_p$ isomorph zu $F(T)_p$.
\\
Daher gibt es ein $J \in M_p^\times$ mit $JAJ^{-1} = \Phi_{F(T)}(JAJ^{-1})$ f\"ur $A \in F(T')_p$,
\\
also $JAJ^{-1} = T\Phi(J)\Phi(A)\Phi(J)^{-1}T^{-1} = T\Phi(J) T'^{-1}\Phi_{F(T')}(A)T'\Phi(J)^{-1}T^{-1}$.
\\
Also gibt es $\lambda \in k_p^\times$ mit $T = \lambda JT'\Phi(J)^{-1}$ und $N(T) = T\Phi(T) = \lambda \overline{\lambda}N(T')$.

\end{itemize}

\item[(\ref{Item03-03.3})]
Wegen $\Phi_F(\Phi(\lambda A)) = \lambda \Phi_F(\Phi(A))$ f\"ur $\lambda \in k$, $A \in M$
ist $\Phi_F \circ \Phi$ ein $k$-Algebren\-automorphismus von $M$.
Nach \cite[Corollary 2.9.9]{Ma} gibt es also ein $T \in M^\times$ mit $\Phi_F(A) = T \Phi(A) T^{-1}$ f\"ur alle $A \in M$,
und aus $\Phi_F^2(A) = \Phi^2(A) = A$ folgt dann $T\Phi(T) \in k^\times$.
Sei $\Phi(T) = \lambda T^*$ mit $\lambda \in k^\times$.
Dann ist $T = \Phi(\lambda T^*) = \overline{\lambda} \lambda T$, und daher gilt $\overline{\lambda} \lambda = 1$.
Also gibt es ein $\mu \in k^\times$ mit $\lambda = \mu \overline{\mu}^{-1}$.
Wir ersetzen $T$ durch $\mu T$.
\\
Dann ist $\Phi(T) = T^*$, und offensichtlich gilt $\Phi_F = \Phi_{F(T)}$, also $F = F(T)$.
\\
F\"ur $S(T)^2 \ne 4N(T)$ sei $X \in E \cap F$. Dann gilt $X = \Phi_F(X) = T\Phi(X)T^{-1} = TXT^{-1}$, also $X \in k[T]$.
Mit $i\sqrt{d}(T-T^*) \in E \cap F$ folgt die Behauptung $X \in \mathbb{Q}[i\sqrt{d}(T-T^*)]$.

\end{enumerate}

\end{proof}

\begin{corollary} \label{Item03-04}

Sei $d \in \mathbb{N}$ quadratfrei, und sei $k = \mathbb{Q}(i\sqrt{d})$ mit Diskriminante $D$.
\\
Sei $F$ eine $\mathbb{Q}$-Quaternionenalgebra mit Zerf\"allungsk\"orper $k$.
\\
Dann gibt es ein $\tau \in \mathbb{Q}^\times$, so dass $F$ isomorph ist zu
$F(\tau) := \left\{ \left. \begin{pmatrix} a & b \\ \tau\overline{ b} & \overline{a} \end{pmatrix} \right| a, b \in k \right\}$.
\\
Damit gilt $\left(\dfrac{F}{p}\right) = \left(\dfrac{\tau,-d}{p}\right)$ f\"ur alle Stellen $p$ von $\mathbb{Q}$.
Man kann o.B.d.A. annehmen, dass $\tau \in \mathbb{Z}$ quadratfrei und teilerfremd zu $D$ ist,
und dass $\tau \equiv 1 \pmod{4}$ gilt, falls $2 \mid d$.
\\
Dann ist $\mathfrak{F}(\tau) := F(\tau) \cap M_2(\mathfrak{o})$ eine $F(\tau)$-Ordnung mit Diskriminante $\tau^2D^2\mathbb{Z}$.
\\
Falls $\tau \ne 1$, gilt $F(\tau) \cap M_2(\mathbb{Q}) \cong \mathbb{Q}(\sqrt{\tau})$
und $\mathfrak{F}(\tau) \cap M_2(\mathbb{Z}) = F(\tau) \cap M_2(\mathbb{Z}) \cong \mathbb{Z}[\sqrt{\tau}]$.

\end{corollary}

\begin{proof}

Seien $E = M_2(\mathbb{Q})$ und $M = M_2(k)$.
Dann sind $E$ eine $\mathbb{Q}$-Quaternionenalgebra
mit $\left(\dfrac{E}{p}\right) = 1$ f\"ur alle Stellen $p$ von $\mathbb{Q}$
und $M$ eine Erweiterung von $E$ zur $k$-Quaternionen\-algebra.
F\"ur $A = \begin{pmatrix} a' & b' \\ c' & d' \end{pmatrix} \in M$ gilt
$\Phi_E(A) = \begin{pmatrix} \overline{a'} & \overline{b'} \\ \overline{c'} & \overline{d'} \end{pmatrix}$.
Nach Voraussetzung kann $F$ in $M$ eingebettet werden,
und nach Satz \ref{Item03-03}.(\ref{Item03-03.3}) gibt es ein $T \in M^\times$ mit $\Phi_E(T) = T^*$, so dass $F$ isomorph ist zu $F(E,T)$.
Nach Satz \ref{Item03-03}.(\ref{Item03-03.2}) h\"angt die Isomorphieklasse von $F$ nur davon ab,
welche Werte $\left(\dfrac{N(T),-d}{p}\right)$ f\"ur die Stellen $p$ von $\mathbb{Q}$ annimmt.
\\
Sei $\tau = N(T)/d$. Wir k\"onnen dann zun\"achst $T$ durch $i\sqrt{d}\begin{pmatrix} 0 & 1 \\ \tau & 0 \end{pmatrix}$ ersetzen.
\\
Damit folgt leicht $F(E,T) = F(\tau)$ und $\left(\dfrac{F}{p}\right) = \left(\dfrac{\tau,-d}{p}\right)$ f\"ur alle Stellen $p$ von $\mathbb{Q}$.
\\
Es gibt genau ein $r \in \mathbb{Q}^+$, so dass $\tau r^2 \in \mathbb{Z}$ quadratfrei ist.
Wir ersetzen $\tau$ durch $\tau r^2$.
\\
Sei $g$ der gr\"o{\ss}te gemeinsame Teiler von $\tau$ und $d$, und sei $\tau = g\tau'$ und $d = gd'$.
Wegen $\left(\dfrac{d + g^2,-d}{p}\right) = 1$ gilt $\left(\dfrac{\tau,-d}{p}\right) = \left(\dfrac{g^2\tau'(d' + g),-d}{p}\right)$.
Falls $g > 1$ ist, ersetzen wir $\tau$ durch den quadratfreien Anteil von $\tau'(d'+g)$.
Dann ist $\tau$ teilerfremd zu $d$.
Falls $\tau$ und $D$ einen gemeinsamen Teiler $g' > 1$ haben, ist $g' = 2$ und $d\equiv 1 \pmod{4}$,
und wir ersetzen $\tau$ durch den quadratfreien Anteil von $\tau(d + 1) / 4$.
Dann ist $\tau$ teilerfremd zu $D$.
\\
Falls $2 \mid d$ und $\tau \not\equiv 1 \pmod{4}$, ist $d \equiv 2 \pmod{4}$ und $\tau \equiv 3 \pmod{4}$.
In diesem Fall  ersetzen wir schlie{\ss}lich $\tau$ durch den quadratfreien Anteil von $\tau(d + 1)$.
\\
Die Diskriminante von $\mathfrak{F}(\tau)$ kann man leicht elementar berechnen.
Die letzte Behauptung folgt aus $F(\tau) \cap M_2(\mathbb{Q})
= \left\{ \left. \begin{pmatrix} a & b \\ \tau  b & a \end{pmatrix} \right| a, b \in \mathbb{Q} \right\}
= \mathbb{Q} \left[ \begin{pmatrix} 0 & 1 \\ \tau & 0 \end{pmatrix} \right]$
sowie \smash{$\begin{pmatrix} 0 & 1 \\ \tau & 0 \end{pmatrix}^2 = \tau$}.

\end{proof} 


\section{Einbettung $\boldsymbol{p}$-adischer Quaternionenordnungen} \label{Chapter04}

\begin{definition} \label{Item04-01}

Seien $F$ eine $\mathbb{Q}$-Quaternionenalgebra und $p$ eine endliche Stelle von $\mathbb{Q}$.
Sei $r \in \mathbb{N}_0$, und sei $\mathfrak{G}_p$ eine $F_p$-Ordnung mit Diskriminante $p^{2r} \mathbb{Z}_p$.
\\
Dann sei $\Delta_p(\mathfrak{G}_p) := \pm p^r$, je nachdem ob $F$ an der Stelle $p$ zerf\"allt oder verzweigt ist.

\end{definition}

\begin{theorem} \label{Item04-02}

Sei $p$ eine endliche Stelle von $\mathbb{Q}$.
\\
Seien $F_p$ eine $\mathbb{Q}_p$-Quaternionenalgebra, $\mathfrak{F}_p$ eine $F_p$-Maximalordnung
und $K_p, K'_p \subset F_p$ zwei quadratische K\"orpererweiterungen von $\mathbb{Q}_p$
mit den Hauptordnungen $\mathfrak{O}_p, \mathfrak{O}'_p \subset \mathfrak{F}_p$.
\\
Sei $p$ verzweigt in $K_p$ und tr\"age in $K'_p$. Sei $\Pi \in \mathfrak{O}_p$ Primelement, und sei $\mathfrak{O}'_p = \mathbb{Z}_p[\Omega]$.
\begin{enumerate}[(i)]

\item \label{Item04-02.1}
Sei $\mathfrak{G}_p$ eine $F_p$-Ordnung mit $\mathfrak{O}_p \subset \mathfrak{G}_p \subset \mathfrak{F}_p$.
\\
Dann gibt es ein $r \in \mathbb{N}_0$
mit $\mathfrak{G}_p = \mathfrak{O}_p + \mathfrak{O}_p \Pi^r \Omega$
und $[\mathfrak{F}_p : \mathfrak{G}_p] = [\mathfrak{O}_p : \mathfrak{O}_p \Pi^r] = p^r$.

\item \label{Item04-02.2}
Sei $\mathfrak{G}_p$ eine $F_p$-Ordnung mit $\mathfrak{O}'_p \subset \mathfrak{G}_p \subset \mathfrak{F}_p$.
\\
Dann gibt es ein $r \in \mathbb{N}_0$
mit $\mathfrak{G}_p = \mathfrak{O}'_p + p^r \mathfrak{O}'_p \Pi$
und $[\mathfrak{F}_p : \mathfrak{G}_p] = [\mathfrak{O}'_p : p^r \mathfrak{O}'_p] = p^{2r}$.

\item  \label{Item04-02.3}
Sei $r \in \mathbb{N}_0$.
Dann sind $\mathfrak{O}_p + \mathfrak{O}_p \Pi^r \Omega$ und $\mathfrak{O}'_p + p^r \mathfrak{O}'_p \Pi$ jeweils $F_p$-Ordnungen.

\item  \label{Item04-02.4}
Sei $F_p \cong M_2(\mathbb{Q}_p)$. Dann gilt:
\\
$\Pi\mathfrak{F}_p\Pi^{-1}$ ist eine $F_p$-Maximalordnung
mit $\mathfrak{O}_p \subset \Pi\mathfrak{F}_p\Pi^{-1} \ne \mathfrak{F}_p$.
\\
Ist $\mathfrak{G}_p$ eine $F_p$-Ordnung mit $\mathfrak{O}_p \subset \mathfrak{G}_p \not\subset \mathfrak{F}_p$,
so ist $\mathfrak{G}_p = \Pi\mathfrak{F}_p\Pi^{-1}$.
\\
Ist $\mathfrak{F}'_p$ eine $F_p$-Maximalordnung mit $\mathfrak{O}'_p \subset \mathfrak{F}'_p$,
so ist $\mathfrak{F}'_p = \mathfrak{F}_p$.

\end{enumerate}

\end{theorem}

\noindent
$Bemerkung$: Man kann jede quadratische K\"orpererweiterung von $\mathbb{Q}_p$ in $F_p$
und ihre Hauptordnung in eine $F_p$-Maximalordnung einbetten.
Alle $F_p$-Maximalordnungen sind isomorph.
Daher gibt es stets die im Satz genannten $K_p, K'_p \subset F_p$ und $\mathfrak{O}_p, \mathfrak{O}'_p \subset \mathfrak{F}_p$.

\begin{proof} ${}$
\begin{enumerate}[(i)]

\item[(\ref{Item04-02.1})]
Wir zeigen zuerst durch Widerspruch, dass $\{ 1, \Pi, \Omega, \Pi \Omega \}$ eine $\mathbb{Z}_p$-Basis von $\mathfrak{F}_p$ ist.
Wir nehmen also an, dass es $\alpha, \beta, \gamma, \delta \in \mathbb{Z}_p$ gibt, nicht alle durch $p$ teilbar,
mit $A: = \alpha + \beta \Pi + \gamma \Omega + \delta \Pi \Omega \in p \mathfrak{F}_p$.
Wir multiplizieren $A$ von links mit $\Pi^*$.
\\
Wegen $p \mid N(\Pi) = \Pi^*  \Pi$ gilt dann auch $B: = \alpha \Pi^* + \gamma \Pi^* \Omega \in p \mathfrak{F}_p$, also $p^2 \mid N(B)$.
Wegen $p^2 \nmid N(\Pi)$ folgt $p \mid N(\alpha + \gamma \Omega)$.
Daraus folgt $\alpha + \gamma \Omega \in p\mathfrak{O}'_p$, also $\alpha, \gamma, \in p \mathbb{Z}_p$.
Dann gilt auch $C: = \beta \Pi + \delta \Pi \Omega \in p \mathfrak{F}_p$.
Wie f\"ur $B$ folgt $\beta, \delta \in p \mathbb{Z}_p$, Widerspruch.
\\
$\{ 1, \Pi, \Omega, \Pi \Omega \}$ ist eine $\mathbb{Z}_p$-Basis und $\{ 1, \Omega \}$ eine $\mathfrak{O}_p$-Basis des Linksmoduls $\mathfrak{F}_p$.

Sei $\mathfrak{B}_p = \{ B \in K_p \mid B \Omega \in \mathfrak{G}_p \}$.
\\
Wegen $\mathfrak{O}_p \subset \mathfrak{G}_p \subset \mathfrak{F}_p$ ist $\mathfrak{B}_p$ ein ganzes $\mathfrak{O}_p$-Ideal,
also gilt $\mathfrak{B}_p = \mathfrak{O}_p \Pi^r$ mit $r \in \mathbb{N}_0$.
Nach Konstruktion gilt $\mathfrak{B}_p \Omega \subset \mathfrak{G}_p$,
wegen $\mathfrak{O}_p \subset \mathfrak{G}_p$ also $\mathfrak{O}_p + \mathfrak{B}_p \Omega \subset \mathfrak{G}_p$.
\\
Seien umgekehrt $A, B \in K_p$ mit $A + B \Omega \in \mathfrak{G}_p$.
\\
Aus $\mathfrak{G}_p \subset \mathfrak{F}_p$ folgt dann $A, B \in \mathfrak{O}_p$,
und wegen $\mathfrak{O}_p \subset \mathfrak{G}_p$ ist $B \Omega \in \mathfrak{G}_p$, also $B \in \mathfrak{B}_p$.
\\
Da $p$ in $K_p$ verzweigt ist, gilt $[\mathfrak{O}_p : \mathfrak{O}_p \Pi] = p$
und $[\mathfrak{F}_p : \mathfrak{G}_p] = [\mathfrak{O}_p : \mathfrak{B}_p]= p^r$.

\item[(\ref{Item04-02.2})]
Wie in (\ref{Item04-02.1}) folgt, dass $\{ 1, \Omega, \Pi, \Omega \Pi \}$ eine $\mathbb{Z}_p$-Basis von $\mathfrak{F}_p$ ist.
(Man multipliziere mit $\Pi^*$ von rechts statt von links.)
Analog folgt dann, dass $\{ 1, \Pi \}$ eine $\mathfrak{O}'_p$-Basis des Linksmoduls $\mathfrak{F}_p$ ist,
und dass $\mathfrak{G}_p = \mathfrak{O}'_p + \mathfrak{B}'_p \Pi$
f\"ur ein ganzes $\mathfrak{O}'_p$-Ideal $\mathfrak{B}'_p = p^r \mathfrak{O}'_p$ ist, mit $r \in \mathbb{N}_0$.
Wegen $[\mathfrak{O}'_p : p \mathfrak{O}'_p] = p^2$
ist dann $[\mathfrak{F}_p : \mathfrak{G}_p] = [\mathfrak{O}'_p : \mathfrak{B}'_p]= p^{2r}$.

\item[(\ref{Item04-02.3})]
Zum Beweis, dass $\mathfrak{G}_p: = \mathfrak{O}_p + \mathfrak{O}_p \Pi^r \Omega$ eine $F_p$-Ordnung ist,
m\"ussen wir nur zeigen, dass $\mathfrak{G}_p$ multiplikativ abgeschlossen ist.
Seien also $A_1, A_2 \in \mathfrak{O}_p$ und $B_1, B_2 \in \mathfrak{O}_p \Pi^r$.
Da $\mathfrak{F}_p$ multiplikativ abgeschlossen ist, ist $\Omega  (A_2 + B_2 \Omega) = A' + B' \Omega$ mit $A', B' \in \mathfrak{O}_p$.
Dann ist $(A_1 + B_1 \Omega) (A_2 + B_2 \Omega) = A'' + B'' \Omega$
mit $A'' = A_1A_2 + B_1 A' \in \mathfrak{O}_p$ und $B'' = A_1 B_2 + B_1 B'\in \mathfrak{O}_p \Pi^r$.
Genauso folgt, dass $\mathfrak{O}'_p + p^r \mathfrak{O}'_p \Pi$ eine $F_p$-Ordnung ist.

\item[(\ref{Item04-02.4})]
Sei $\mathfrak{F}'_p$ eine $F_p$-Maximalordnung mit $\mathfrak{O}_p \subset \mathfrak{F}'_p$.
Aus $\mathfrak{O}_p = K_p \cap \mathfrak{F}_p = K_p \cap \mathfrak{F}'_p$
folgt mit Lemma \ref{Item02-03}, dass $\mathfrak{F}'_p = X\mathfrak{F}_pX^{-1}$ f\"ur ein $X \in K^\times_p$.
Wegen $\Pi^2p^{-1} \in \mathfrak{O}_p^\times$ gibt es $r \in \mathbb{Z}$ und $Y \in \mathfrak{O}_p^\times$
mit $X = p^r Y$ oder $X = \Pi p^r Y$.
Also ist $\mathfrak{F}'_p = \mathfrak{F}_p$ oder $\mathfrak{F}'_p = \Pi\mathfrak{F}_p\Pi^{-1}$,
und wegen $\Pi \not\in \mathbb{Q}_p^\times \mathfrak{F}_p^\times$
gilt $\Pi\mathfrak{F}_p\Pi^{-1} \ne \mathfrak{F}_p$, siehe \cite[Theorem 6.5.3.2]{Ma}.
Weiter folgt $\mathfrak{G}_p \subset \mathfrak{F}_p$ oder $\mathfrak{G}_p \subset \Pi\mathfrak{F}_p\Pi^{-1}$.
Falls $\mathfrak{G}_p \subsetneqq \Pi\mathfrak{F}_p\Pi^{-1}$,
gibt es nach (\ref{Item04-02.1}) ein $r \in \mathbb{N}$ mit
$\mathfrak{G}_p = \mathfrak{O}_p + \mathfrak{O}_p \Pi^r\Pi\Omega\Pi^{-1} =
\mathfrak{O}_p + \mathfrak{O}_p N(\Pi)^{-1}\Pi^{r+1}\Omega\Pi^* =
\mathfrak{O}_p + \mathfrak{O}_p \Pi^{r-1}\Omega\Pi^* \subset \mathfrak{F}_p$.

Die Behauptung $\mathfrak{F}'_p = \mathfrak{F}_p$ f\"ur $\mathfrak{O}'_p \subset \mathfrak{F}'_p$ folgt analog, da $p \in \mathfrak{O}'_p$ ein Primelement ist.

\end{enumerate}

\end{proof}

\begin{theorem} \label{Item04-03}

Sei $p$ eine endliche Stelle von $\mathbb{Q}$.
\\
Seien $F_p \cong M_2(\mathbb{Q}_p)$ eine $\mathbb{Q}_p$-Quaternionenalgebra,
$\mathfrak{F}_p$ eine $F_p$-Maximalordnung
und \mbox{$K_p \subset F_p$} eine halbeinfache quadratische Erweiterung von $\mathbb{Q}_p$
mit Hauptordnung $\mathfrak{O}_p \subset \mathfrak{F}_p$.
Sei $p$ zerlegt in $K_p$, sei also $K_p$ kein K\"orper.
\begin{enumerate} [(i)]

\item \label{Item04-03.1}
Sei $\mathfrak{G}_p$ eine $F_p$-Ordnung mit $\mathfrak{O}_p \subset  \mathfrak{G}_p \subset  \mathfrak{F}_p$.
\\
Dann gibt es Matrixeinheiten in $F_p$, bez\"uglich derer
$\mathfrak{O}_p = \begin{pmatrix} \mathbb{Z}_p & 0 \\ 0 & \mathbb{Z}_p \end{pmatrix}$.
Weiter gibt es dann $r, s \in \mathbb{Z}$ mit
$\mathfrak{G}_p = \begin{pmatrix} \mathbb{Z}_p & p^s \mathbb{Z}_p \\ p^r  \mathbb{Z}_p & \mathbb{Z}_p \end{pmatrix}$.
Damit ist $[\mathfrak{F}_p : \mathfrak{G}_p] = \Delta_p(\mathfrak{G}_p) = p^{r+s}$.

\item \label{Item04-03.2}
Seien $r,s \in \mathbb{Z}$ mit $r + s \ge 0$, und seien (beliebige) Matrixeinheiten in $F_p$ gegeben.
\\
Bez\"uglich dieser Matrixeinheiten ist dann 
$\begin{pmatrix} \mathbb{Z}_p & p^s \mathbb{Z}_p \\ p^r \mathbb{Z}_p & \mathbb{Z}_p \end{pmatrix}$ eine $F_p$-Ordnung.

\end{enumerate}

\end{theorem}

\begin{proof} ${}$
\begin{enumerate} [(i)]

\item [(\ref{Item04-03.1})]
Je zwei zerfallende halbeinfache quadratische Erweiterungen von $\mathbb{Q}_p$ sind isomorph zueinander.
Nach Lemma \ref{Item02-01} gibt es also Matrixeinheiten $U_{11}, U_{12}, U_{21}, U_{22} \in F_p$,
bez\"uglich derer $\mathfrak{O}_p = \begin{pmatrix} \mathbb{Z}_p & 0 \\ 0 & \mathbb{Z}_p \end{pmatrix}$ gilt.
F\"ur $A = \begin{pmatrix} \alpha & \beta \\ \gamma & \delta \end{pmatrix} \in \mathfrak{G}_p$
mit $\alpha, \beta, \gamma, \delta \in \mathbb{Q}_p$
gilt dann $\beta U_{12} = U_{11} A U_{22} \in \mathfrak{G}_p$
und $\gamma U_{21} = U_{22} A U_{11} \in \mathfrak{G}_p$
sowie $\alpha, \delta \in \mathbb{Z}_p$.
Also ist $\{ U_{11}, p^s U_{12}, p^r U_{21}, U_{22} \}$ mit geeigneten $r,s \in \mathbb{Z}$ eine $\mathbb{Z}_p$-Basis von $\mathfrak{G}_p$.

\item [(\ref{Item04-03.2})]
leicht nachzurechnen

\end{enumerate}

\end{proof}

\begin{lemma} \label{Item04-04}

Sei $p$ eine endliche Stelle von $\mathbb{Q}$, und sei $F_p$ eine $\mathbb{Q}_p$-Quaternionenalgebra.
\\
Seien $\mathfrak{G}_p, \mathfrak{G}'_p$ zwei $F_p$-Ordnungen mit $\Delta_p(\mathfrak{G}_p) =  \Delta_p(\mathfrak{G}'_p)$,
und seien $K_p, K'_p \subset F_p$ halbeinfache quadratische Erweiterungen von $\mathbb{Q}_p$
mit den Hauptordnungen $\mathfrak{O}_p \subset \mathfrak{G}_p$, $\mathfrak{O}'_p \subset \mathfrak{G}'_p$.
\\
Ein Isomorphismus $\mathfrak{O}_p \rightarrow \mathfrak{O}'_p$ ist dann zu einem Isomorphismus $\mathfrak{G}_p \rightarrow \mathfrak{G}'_p$ fortsetzbar.

\end{lemma}

\begin{proof}

$\mathfrak{F}_p, \mathfrak{F}'_p$ seien $F_p$-Maximalordnungen
mit $\mathfrak{G}_p \subset \mathfrak{F}_p$, $\mathfrak{G}'_p \subset \mathfrak{F}'_p$.
Nach Lemma \ref{Item02-01} ist ein Isomorphismus $\mathfrak{O}_p \rightarrow \mathfrak{O}'_p$
zu einem Automorphismus $j: F_p \rightarrow F_p$ fortsetzbar.
Nach Lemma \ref{Item02-03} gibt es wegen $K'_p \cap j(\mathfrak{F}_p) = \mathfrak{O}'_p = K'_p \cap \mathfrak{F}'_p$.
ein $X \in K'^\times_p$ mit $\mathfrak{F}'_p = Xj(\mathfrak{F}_p)X^{-1}$
und $\mathfrak{O}'_p \subset Xj(\mathfrak{G}_p)X^{-1} \subset \mathfrak{F}'_p$.
Mit Satz \ref{Item04-02} folgt $Xj(\mathfrak{G}_p)X^{-1} = \mathfrak{G}'_p$, falls $K_p$ K\"orper ist.
\\
Falls $K_p$ zerf\"allt, gibt es nach Satz \ref{Item04-03}.(\ref{Item04-03.1})
Matrixeinheiten in $F_p$ und $r, s, r', s' \in \mathbb{Z}$ mit $r + s = r' + s'$, bez\"uglich derer gilt:
$\mathfrak{O}'_p = \begin{pmatrix} \mathbb{Z}_p & 0 \\ 0 & \mathbb{Z}_p \end{pmatrix}$,
$Xj(\mathfrak{G}_p)X^{-1} = \begin{pmatrix} \mathbb{Z}_p & p^s \mathbb{Z}_p \\ p^r \mathbb{Z}_p & \mathbb{Z}_p \end{pmatrix}$
und $\mathfrak{G}'_p = \begin{pmatrix} \mathbb{Z}_p & p^{s'}\mathbb{Z}_p \\ p^{r'}\mathbb{Z}_p & \mathbb{Z}_p \end{pmatrix}$.
Mit $Y := \begin{pmatrix} p^r & 0 \\ 0 & p^{r'} \end{pmatrix} \in K'^\times_p$ folgt $YXj(\mathfrak{G}_p)X^{-1}Y^{-1} = \mathfrak{G}'_p$.

\end{proof}

\begin{theorem} \label{Item04-05}

Sei $p$ eine endliche Stelle von $\mathbb{Q}$, und sei $F_p$ eine $\mathbb{Q}_p$-Quaternionenalgebra.
\\
Sei $\mathfrak{F}_p$ eine $F_p$-Maximalordnung,
und sei $\mathfrak{G}_p \subset \mathfrak{F}_p$ eine $F_p$-Ordnung mit $[\mathfrak{F}_p : \mathfrak{G}_p] = p$.
\begin{enumerate}[(i)]

\item \label{Item04-05.1}
Falls $F_p$ Divisionsalgebra ist, gilt $\mathfrak{G}_p = \{ A \in \mathfrak{F}_p \mid S(A)^2  - 4N(A) \in p\mathbb{Z}_p \}$.

\item \label{Item04-05.2}
Falls $F_p$ zerf\"allt,
gibt es Matrixeinheiten in $F_p$, bez\"uglich derer
$\mathfrak{G}_p = \begin{pmatrix} \mathbb{Z}_p & \mathbb{Z}_p \\ p\mathbb{Z}_p & \mathbb{Z}_p \end{pmatrix}$.

\end{enumerate}

\end{theorem}

\begin{proof} ${}$
\begin{enumerate}[(i)]

\item[(\ref{Item04-05.1})]
Seien $K_p, K'_p \subset F_p$ quadratische K\"orpererweiterungen von $\mathbb{Q}_p$ mit den Hauptordnungen $\mathfrak{O}_p, \mathfrak{O}'_p$.
Sei $p$ verzweigt in $K_p$ und tr\"age in $K'_p$. Sei $\Pi \in \mathfrak{O}_p$ Primelement, und sei $\mathfrak{O}'_p = \mathbb{Z}_p[\Omega]$.
Da $\mathfrak{F}_p$ die einzige $F_p$-Maximalordnung ist, gilt $K_p \cap \mathfrak{F}_p = \mathfrak{O}_p$ und $K'_p \cap \mathfrak{F}_p = \mathfrak{O}'_p$.
Wir zeigen zun\"achst durch Widerspruch, dass $\Pi \in \mathfrak{G}_p$.

Es gibt $\alpha, \beta, \gamma, \beta', \gamma', \gamma'' \in \mathbb{Z}_p$, so dass
$\{ 1, \alpha \Omega + \beta \Omega \Pi + \gamma \Pi, \beta' \Omega \Pi + \gamma' \Pi, \gamma'' \Pi \}$
eine $\mathbb{Z}_p$-Basis von $\mathfrak{G}_p$ ist.
Falls $\Pi \notin \mathfrak{G}_p$, k\"onnen wir wegen $[\mathfrak{F}_p : \mathfrak{G}_p] = p$ annehmen,
dass $\alpha = \beta' = 1$, $\gamma'' = p$, $\beta = 0$.
Und wir k\"onnen $\Omega$ durch $\Omega + \gamma'$ ersetzen, also annehmen,
dass $\{ 1, \Omega + \gamma \Pi, \Omega \Pi, p\Pi \}$ eine $\mathbb{Z}_p$-Basis von $\mathfrak{G}_p$ ist.
Wegen $\Pi^{-1} \mathfrak{F}_p \Pi = \mathfrak{F}_p$
ist $N(\Omega) \Pi = (\Omega + \gamma \Pi)^* \Omega \Pi - \gamma \Pi^* \Omega \Pi =
(\Omega + \gamma \Pi)^* \Omega \Pi - \gamma N(\Pi) \Pi^{-1} \Omega \Pi \in \mathfrak{G}_p + p \mathfrak{F}_p$,
und deshalb $\Pi \in \mathfrak{G}_p + p \mathfrak{F}_p \subset \mathfrak{G}_p$, Widerspruch.

Nach Satz \ref{Item04-02}.(\ref{Item04-02.1}) gilt also $\mathfrak{G}_p = \mathfrak{O}_p + \mathfrak{O}_p \Pi \Omega$,
und $\{ 1, \Pi, p \Omega, \Pi \Omega \}$ ist eine $\mathbb{Z}_p$-Basis von $\mathfrak{G}_p$.
Damit folgt, dass $K_p \cap \mathfrak{G}_p = \mathfrak{O}_p = \{ A \in K_p \cap \mathfrak{F}_p \mid (S(A)^2  - 4N(A))\in p\mathbb{Z}_p \}$
und $K'_p \cap \mathfrak{G}_p = \mathbb{Z}_p + p\mathbb{Z}_p\Omega = \{ A \in K'_p \cap \mathfrak{F}_p \mid (S(A)^2  - 4N(A))\in p\mathbb{Z}_p \}$.

Sei nun $A \in \mathfrak{F}_p$ beliebig.
Dann ist $A$ in einer quadratischen K\"orpererweiterung $K''_p \subset F_p$ von $\mathbb{Q}_p$ enthalten.
$p$ ist entweder verzweigt oder tr\"age in $K''_p$.
Wir k\"onnen o.B.d.A. annehmen, dass $K''_p = K_p$ oder $K''_p = K'_p$, womit die Behauptung folgt.

\item[(\ref{Item04-05.2})]
siehe \cite[Satz 3]{Ei}

\end{enumerate}

\end{proof}

\begin{lemma} \label{Item04-06}

Sei $d \in \mathbb{N}$ quadratfrei, und sei $k = \mathbb{Q}(i\sqrt{d})$.
\\
Seien $F$ eine $\mathbb{Q}$-Quaternionenalgebra und $M$ eine Erweiterung von $F$ zur $k$-Quaternionen\-algebra.
Seien $\mathfrak{F}$ eine $F$-Maximalordnung und $p$ eine endliche Verzweigungsstelle von $F$.
\begin{enumerate}[(i)]

\item \label{Item04-06.1}
Falls $p$ in $k$ verzweigt ist, gibt es genau eine $M_p$-Maximalordnung $\mathfrak{M}_p$ mit $\mathfrak{F}_p \subset \mathfrak{M}_p$.

\item \label{Item04-06.2}
Falls $p$ in $k$ tr\"age ist, gibt es genau zwei $M_p$-Maximalordnungen $\mathfrak{M}_p \ne \mathfrak{M}'_p$ mit
$\mathfrak{F}_p \subset \mathfrak{M}_p$.
Es gilt $[\mathfrak{M}_p : (\mathfrak{M}_p \cap \mathfrak{M}'_p)] = p^2$.

\end{enumerate}

\end{lemma}

\begin{proof}

Sei $K_p \subset F_p$ quadratische K\"orpererweiterung von  $\mathbb{Q}_p$ mit Hauptordnung $\mathfrak{O}_p$.
Falls $p$ verzweigt in $k$ ist, sei $p$ tr\"age in $K_p$. Falls $p$ tr\"age in $k$ ist, sei $p$ verzweigt in $K_p$.
$L_p := k_p K_p$ ist quadratische K\"orpererweiterung von $k_p$ mit Hauptordnung $\mathfrak{L}_p = \mathfrak{o}_p \mathfrak{O}_p$.
Seien $\mathfrak{M}_p$ und $\mathfrak{M}'_p$ zwei $M_p$-Maximalordnungen
mit $\mathfrak{F}_p \subset \mathfrak{M}_p$ und $\mathfrak{F}_p \subset\mathfrak{M}'_p$.
Dann gilt $\mathfrak{L}_p = \mathfrak{o}_p \mathfrak{O}_p = \mathfrak{o}_p (K_p \cap \mathfrak{M}_p)
\subset L_p  \cap \mathfrak{M}_p \subset \mathfrak{L}_p$,
also gilt $\mathfrak{L}_p = L_p \cap \mathfrak{M}_p$ und genauso $\mathfrak{L}_p = L_p \cap \mathfrak{M}'_p$.
Daher gibt es nach Lemma \ref{Item02-03} ein $X \in L^\times_p$ mit $\mathfrak{M}'_p = X\mathfrak{M}_pX^{-1}$.
\begin{enumerate}[(i)]

\item[(\ref{Item04-06.1})]
Sei $\pi \in \mathfrak{o}_p$ Primelement. Dann ist $\pi$ auch Primelement in $\mathfrak{L}_p$.
Daher gibt es $r \in \mathbb{Z}$ und $Y \in \mathfrak{L}_p^\times \subset \mathfrak{M}_p^\times$ mit $X = \pi^r Y$.
Also gilt $\mathfrak{M}'_p = \pi^rY\mathfrak{M}_pY^{-1}\pi^{-r} = \mathfrak{M}_p$.

\item[(\ref{Item04-06.2})]
Sei $\Pi \in \mathfrak{O}_p$ Primelement, also auch Primelement in $\mathfrak{L}_p$. Es gilt $\Pi^2p^{-1} \in \mathfrak{O}_p^\times$.
Daher gibt es $r \in \mathbb{Z}$ und $Y \in \mathfrak{L}_p^\times \subset \mathfrak{M}_p^\times$ mit $X = p^r Y$ oder $X = \Pi p^r Y$.
Also gilt $\mathfrak{M}'_p = \mathfrak{M}_p$ oder $\mathfrak{M}'_p = \Pi \mathfrak{M}_p \Pi^{-1}$.
Sei nun $\mathfrak{M}'_p = \Pi \mathfrak{M}_p \Pi^{-1}$.
Die $k_p$-Algebra $M_p$ zerf\"allt, siehe Lemma \ref{Item03-02}.(\ref{Item03-02.1}).
Wegen $\Pi \not\in k_p^\times \mathfrak{M}_p^\times$ gilt daher $\mathfrak{M}'_p \ne \mathfrak{M}_p$, siehe \cite[Theorem 6.5.3.2]{Ma}.
\\
Aus $N(\Pi) \in p\mathfrak{O}_p^\times$ folgt mit Lemma \ref{Item02-02}, dass
$[\mathfrak{M}_p : (\mathfrak{M}_p \cap \mathfrak{M}'_p)] = N_{abs}( p\mathfrak{O}_p) = p^2$.

\end{enumerate}

\end{proof}

\begin{lemma} \label{Item04-07}

Sei $d \in \mathbb{N}$ quadratfrei, und sei $k = \mathbb{Q}(i\sqrt{d})$.
\\
Seien $F$ eine $\mathbb{Q}$-Quaternionenalgebra und $M$ eine Erweiterung von $F$ zur $k$-Quaternionen\-algebra.
Seien $\mathfrak{F}$ eine $F$-Maximalordnung und $p$ eine endliche Stelle von  $\mathbb{Q}$.
Falls $F_p$ zerf\"allt oder $p$ in $k$ zerlegt ist,
gibt es genau eine $M_p$-Maximalordnung $\mathfrak{M}_p$ mit $\mathfrak{F}_p \subset \mathfrak{M}_p$.

\end{lemma}

\begin{proof}

Falls $F_p$ zerf\"allt, hat $\mathfrak{F}_p$ die Diskriminante $\mathbb{Z}_p$. 
Dann hat $\mathfrak{M}_p := \mathfrak{o}_p \mathfrak{F}_p$ die Diskriminante $\mathfrak{o}_p$,
ist also die einzige $M_p$-Maximalordnung mit $\mathfrak{F}_p \subset \mathfrak{M}_p$.
\\
Falls $F_p$ Divisionsalgebra ist und $p$ in $k$ zerlegt ist, gilt $M_p \cong F_p \times F_p$,
und $\mathfrak{M}_p \cong \mathfrak{F}_p \times \mathfrak{F}_p$ ist die einzige $M_p$-Maximalordnung.
($\mathfrak{F}_p$ ist darin diagonal eingebettet.)

\end{proof}

\begin{lemma} \label{Item04-08}

Sei $d \in \mathbb{N}$ quadratfrei, und sei $k = \mathbb{Q}(i\sqrt{d})$ mit Hauptordnung $\mathfrak{o}$.
\\
Seien $F$ eine $\mathbb{Q}$-Quaternionenalgebra, $M$ eine Erweiterung von $F$ zur $k$-Quaternionen\-algebra
und $\mathfrak{N}$ eine $M$-Maximalordnung. Sei $p$ eine endliche Stelle von $\mathbb{Q}$ mit $p \nmid \Sigma_k(F)$.
\begin{enumerate}[(i)]

\item \label{Item04-08.1}
Dann enth\"alt $F_p \cap \mathfrak{N}_p$ eine zu $\mathfrak{o}_p$ isomorphe Ordnung.

\item \label{Item04-08.2}
Sei $\mathfrak{F}_p$ eine $F_p$-Maximalordnung mit $F_p \cap \mathfrak{N}_p \subset \mathfrak{F}_p$.
Dann gibt es genau eine \mbox{$M_p$-Maximalordnung} $\mathfrak{M}_p$ mit $\mathfrak{F}_p \subset \mathfrak{M}_p$
und $[\mathfrak{F}_p : (F_p \cap \mathfrak{N}_p)] = [\mathfrak{M}_p : (\mathfrak{M}_p \cap \mathfrak{N}_p)]$.
Die Inklusion $\mathfrak{F}_p \hookrightarrow \mathfrak{M}_p$ induziert einen Isomorphismus
$\mathfrak{F}_p / (F_p \cap \mathfrak{N}_p) \rightarrow \mathfrak{M}_p / (\mathfrak{M}_p \cap \mathfrak{N}_p)$.

\end{enumerate}

\end{lemma}

\begin{proof}

Sei $\mathfrak{G}_p = F_p \cap \mathfrak{N}_p$, und sei $\mathfrak{F}_p$ eine $F_p$-Maximalordnung mit $\mathfrak{G}_p \subset \mathfrak{F}_p$.
\begin{easylist}[checklist]
\ListProperties(Style*=\large$\bullet$,Style**=)
& ${ }$ Sei zun\"achst $p$ zerlegt in $k$.
\end{easylist}
\noindent
$F_p$ ist diagonal in $M_p  \cong F_p \times F_p$ eingebettet,
und wegen $p \nmid \Sigma_k(F)$ gilt $F_p \cong M_2(\mathbb{Q}_p)$.
Daher gilt $\mathfrak{N}_p \cong \mathfrak{F}_p \times \mathfrak{F}'_p$
mit zwei $F_p$-Maximalordnungen $\mathfrak{F}'_p$, $\mathfrak{F}''_p$.
Es gibt ein $r \in \mathbb{N}_0$ und Matrixeinheiten in $F_p$,
bez\"uglich derer $\mathfrak{F}'_p = \begin{pmatrix} \mathbb{Z}_p & \mathbb{Z}_p \\ \mathbb{Z}_p & \mathbb{Z}_p \end{pmatrix}$
und $\mathfrak{F}''_p = \begin{pmatrix} \mathbb{Z}_p & p^{-r} \mathbb{Z}_p \\ p^r \mathbb{Z}_p & \mathbb{Z}_p \end{pmatrix}$ gilt.
Damit gilt $\mathfrak{G}_p =
F_p \cap \mathfrak{N}_p = \begin{pmatrix} \mathbb{Z}_p & \mathbb{Z}_p \\ p^r \mathbb{Z}_p & \mathbb{Z}_p \end{pmatrix}$
und $\mathfrak{o}_p \cong \begin{pmatrix} \mathbb{Z}_p & 0 \\ 0 & \mathbb{Z}_p \end{pmatrix} \subset \mathfrak{G}_p$.
Nach Satz \ref{Item04-03}.(\ref{Item04-03.1}) gilt
$\mathfrak{F}_p = \begin{pmatrix} \mathbb{Z}_p & p^{-s} \mathbb{Z}_p \\ p^s \mathbb{Z}_p & \mathbb{Z}_p \end{pmatrix}$ mit $s \in \mathbb{Z}$.
Wegen $\mathfrak{G}_p \subset \mathfrak{F}_p$ gilt $0 \le s \le r$ und $[\mathfrak{F}_p : \mathfrak{G}_p)] = p^r$.
Nach Lemma \ref{Item04-07} ist $\mathfrak{M}_p \cong \mathfrak{F}_p \times \mathfrak{F}_p$
die einzige $M_p$-Maximalordnung mit $\mathfrak{F}_p \subset \mathfrak{M}_p$,
und aus $\mathfrak{M}_p \cap \mathfrak{N}_p \cong 
\begin{pmatrix} \mathbb{Z}_p & \mathbb{Z}_p \\ p^s \mathbb{Z}_p & \mathbb{Z}_p \end{pmatrix}
\times \begin{pmatrix} \mathbb{Z}_p & p^{-s} \mathbb{Z}_p \\ p^r \mathbb{Z}_p & \mathbb{Z}_p \end{pmatrix}$
folgt $[\mathfrak{M}_p : (\mathfrak{M}_p \cap \mathfrak{N}_p)] = p^s p^{r-s} = p^r$.
Damit induziert die Einbettung $\mathfrak{F}_p \hookrightarrow \mathfrak{M}_p$
einen Isomorphismus $\mathfrak{F}_p / \mathfrak{G}_p \rightarrow \mathfrak{M}_p / (\mathfrak{M}_p \cap \mathfrak{N}_p)$.
\begin{easylist}[checklist]
\ListProperties(Style*=\large$\bullet$,Style**=)
& ${ }$ Sei nun $p$ verzweigt oder tr\"age in $k$.
\end{easylist}
\noindent
$k_p$ ist quadratische K\"orpererweiterung von $\mathbb{Q}_p$,
und nach Lemma \ref{Item03-02}.(\ref{Item03-02.1}) gilt $M_p \cong M_2(k_p)$.
Sei $\mathfrak{o}_p = \mathbb{Z}_p[\omega]$,
und sei $\pi \in \mathfrak{o}_p$ ein Primelement.
Falls $p$ in $k$ verzweigt ist, sei $\omega = \pi$.
Falls $p$ in $k$ tr\"age ist, sei $\pi = p$.
Sei $\mathfrak{M}_p$ eine $M_p$-Maximalordnung mit $\mathfrak{F}_p\subset \mathfrak{M}_p$.
Falls $p$ in $k$ tr\"age ist, m\"ussen wir $\mathfrak{M}_p$ im Beweisverlauf eventuell ersetzen, siehe Lemma \ref{Item04-06}.
Die induzierte Abbildung $\mathfrak{F}_p / \mathfrak{G}_p \rightarrow \mathfrak{M}_p / (\mathfrak{M}_p \cap \mathfrak{N}_p)$ ist injektiv.
Es gibt ein $r \in \mathbb{N}_0$ und Matrixeinheiten $U_{11}, U_{12}, U_{21}, U_{22} \in M_p$, bez\"uglich derer
$\mathfrak{M}_p = \begin{pmatrix} \mathfrak{o}_p & \mathfrak{o}_p \\ \mathfrak{o}_p & \mathfrak{o}_p \end{pmatrix}$ und
$\mathfrak{N}_p = \begin{pmatrix} \mathfrak{o}_p & \pi^{-r}\mathfrak{o}_p \\ \pi^{r}\mathfrak{o}_p & \mathfrak{o}_p \end{pmatrix}$ gilt.
\\
Wir zeigen (\ref{Item04-08.1}), (\ref{Item04-08.2}) zuerst  im Fall $\pi^r \notin p\mathfrak{o}_p$:
Dann ist entweder $r = 0$, also $\mathfrak{M}_p = \mathfrak{N}_p$ und $\mathfrak{G}_p = \mathfrak{F}_p$, also nichts weiter zu zeigen.
Oder es ist $r = 1$, und $p$ ist verzweigt in $k$.
Dann ist $[\mathfrak{F}_p : \mathfrak{G}_p] \le [\mathfrak{M}_p : (\mathfrak{M}_p \cap \mathfrak{N}_p)] = p$.
Falls $\mathfrak{G}_p = \mathfrak{F}_p$, w\"are $\mathfrak{M}_p = \mathfrak{N}_p$, also doch $r = 0$.
Also ist $[\mathfrak{F}_p : \mathfrak{G}_p] = p$, und nach Satz \ref{Item04-05} ist $\mathfrak{G}_p$ bis auf Isomorphie eindeutig,
enth\"alt nach Satz \ref{Item04-02} also eine zu $\mathfrak{o}_p$ isomorphe Ordnung,
und $\mathfrak{F}_p / \mathfrak{G}_p \rightarrow \mathfrak{M}_p / (\mathfrak{M}_p \cap \mathfrak{N}_p)$ ist surjektiv.
\\[4pt]
Wir k\"onnen nun also annehmen, dass $\pi^r \in p\mathfrak{o}_p$ oder gleichwertig $[\mathfrak{M}_p : (\mathfrak{M}_p \cap \mathfrak{N}_p)] \ge p^2$.
\begin{enumerate}[(i)]

\item
Wir entwickeln (unabh\"angig von $r$) ein $U \in \mathfrak{G}_p$, so dass $\mathbb{Z}_p[U]$ isomorph zu $\mathfrak{o}_p$ ist.
$F_p$ und $\begin{pmatrix} k_p & k_p \\ 0 & k_p \end{pmatrix}$ haben als $\mathbb{Q}_p$-Algebren  die Dimensionen $4$ und $6$.
Ihr Durchschnitt hat also eine Dimension $\ge 2$.
Daher gibt es $a, b, c \in k_p$ mit $U: = \begin{pmatrix} a & b \\ 0 & c \end{pmatrix} \in F_p \smallsetminus \mathbb{Q}_p$.

Falls $c \ne \overline{a}$, gilt $c - \overline{a} = S(U) - (a + \overline{a}) \in \mathbb{Q}_p^\times$,
und wegen $a(c - \overline{a}) = N(U) - a\overline{a} \in \mathbb{Q}_p$ sind $a, c \in \mathbb{Q}_p$.
Wir ersetzen $U$ durch $U - c$, k\"onnen also $c = 0$ und $a \ne 0$ annehmen.
Damit folgt leicht $UM_pU^* \subset k_p U_{12}$ und $UU_{12}U^* = a^2U_{12}$.
Also gilt $UM_pU^* \ne  \{0\}$ und wegen $M_p = F_p + \omega F_p$ auch $UF_pU^* \ne  \{0\}$.
Sei $X \in F_p$ mit $UXU^* \ne 0$. Dann gilt $X \not\in \mathbb{Q}_p$.
Wir ersetzen $U$ durch  $UXU^* \in F_p \smallsetminus \mathbb{Q}_p$ und erhalten so $c = \overline{a} = 0$.

Wir k\"onnen also $U = \begin{pmatrix} a & b \\ 0 & \overline{a} \end{pmatrix}$ annehmen, mit $a, b \in \mathfrak{o}_p$, nicht beide durch $p$ teilbar.
Es gibt $\alpha \in \mathbb{Z}_p$, $\beta \in \mathbb{Z}_p^\times$ und $s \in \mathbb{N}_0$ mit $a = \alpha$ oder $a = \alpha + \beta p^s \omega$.
Sei $t \in \mathbb{N}_0$ maximal mit $p^t \mid (a - \alpha)$ und $p^t \mid b$.
Wir ersetzen $U$ durch $(U - \alpha)/(\beta p^t)$.
Dann gilt $a = 0$ oder $a = p^s \omega$ mit $s \in \mathbb{N}_0$.
Wir zeigen durch Widerspruch, dass $a \ne 0$, $s = 0$:
\\
Andernfalls gilt $p^{-1} U \in \begin{pmatrix} \mathfrak{o}_p & p^{-1}\mathfrak{o}_p \\ 0 & \mathfrak{o}_p \end{pmatrix}$,
wegen $\pi^r \in p\mathfrak{o}_p$ also $p^{-1} U \in \mathfrak{N}_p$.
Daher gilt $p^{-1}U \in F_p \cap \mathfrak{N}_p = \mathfrak{G}_p$,
andererseits gilt wegen $p \nmid b$ aber $p^{-1}U \notin \mathfrak{M}_p$, und deshalb $p^{-1}U \notin \mathfrak{F}_p$, Widerspruch.
Wir k\"onnen also $U = \begin{pmatrix} \omega & b \\ 0 & \overline{\omega} \end{pmatrix}$ annehmen, mit $b \in \mathfrak{o}_p$.

\item
Wir zeigen die Surjektivit\"at von $\mathfrak{F}_p / \mathfrak{G}_p \rightarrow \mathfrak{M}_p / (\mathfrak{M}_p \cap \mathfrak{N}_p)$.
Dazu gen\"ugt es, die Surjektivit\"at von $f_{21}: \mathfrak{F}_p \rightarrow \mathfrak{o}_p$
mit Zuordnungsvorschrift $\begin{pmatrix} \alpha & \beta \\ \gamma & \delta \end{pmatrix} \mapsto \gamma$ zu zeigen.

Falls $p$ in $k$ verzweigt ist,
sei $K_p \subset F_p$ quadratische K\"orpererweiterung von $\mathbb{Q}_p$ mit Hauptordnung $\mathfrak{O}_p \subset \mathfrak{F}_p$.
Sei $p$ tr\"age in $K_p$, und sei $\mathfrak{O}_p = \mathbb{Z}_p[\Omega]$.
Es gibt $\alpha, \beta, \gamma, \delta \in \mathfrak{o}_p$ mit
$\Omega = \begin{pmatrix} \alpha & \beta \\ \gamma & \delta \end{pmatrix}$
und $\Omega U = \begin{pmatrix} \alpha \omega & \alpha b + \beta \overline{\omega} \\ \gamma \omega & \gamma b + \delta \overline{\omega} \end{pmatrix}$.
Nach \ref{Item04-06}.(\ref{Item04-06.1}) und \ref{Item04-07} gilt
$\mathfrak{F}_p \not\subset \begin{pmatrix} \mathfrak{o}_p & \pi^{-1} \mathfrak{o}_p \\ \pi\mathfrak{o}_p & \mathfrak{o}_p \end{pmatrix}$,
also $\gamma \in \mathfrak{o}_p^\times$,
$\{ \gamma, \gamma \omega \}$ ist $\mathbb{Z}_p$-Basis von $\mathfrak{o}_p$,
$f_{21}$ ist surjektiv.

Falls $p$ in $k$ tr\"age ist,
sei $K_p \subset F_p$ quadratische K\"orpererweiterung von $\mathbb{Q}_p$ mit Hauptordnung $\mathfrak{O}_p \subset \mathfrak{F}_p$.
Sei $p$ verzweigt in $K_p$, und $\Pi \in \mathfrak{O}_p$ Primelement.
Es gibt $\alpha, \beta, \gamma, \delta \in \mathfrak{o}_p$ mit
$\Pi = \begin{pmatrix} \alpha & \beta \\ \gamma & \delta \end{pmatrix}$.
Falls $\gamma \in \mathfrak{o}_p^\times$, folgt wie eben, dass $f_{21}$ surjektiv ist.
\\
Wir d\"urfen also $\gamma = p^s \gamma'$ mit $\gamma' \in \mathfrak{o}_p^\times$, $s \in \mathbb{N}$ annehmen.
Mit $\Pi U = \begin{pmatrix} \alpha \omega & \alpha b + \beta \overline{\omega} \\
\gamma \omega & \gamma b + \delta \overline{\omega} \end{pmatrix}$
folgt $\mathfrak{F}_p \subset \begin{pmatrix} \mathfrak{o}_p & \mathfrak{o}_p \\ p^s \mathfrak{o}_p & \mathfrak{o}_p \end{pmatrix}$.
Also ist $F_p$ Divisionsalgebra, mit $\Delta_p(\mathfrak{F}_p) = -p$ folgt $s = 1$.
\\
Bez\"uglich der Matrixeinheiten $U_{11}, p^{-1}U_{12}, p U_{21}, U_{22}$ gilt
$\mathfrak{N}_p= \begin{pmatrix} \mathfrak{o}_p & p^{1-r}\mathfrak{o}_p \\ p^{r-1}\mathfrak{o}_p & \mathfrak{o}_p \end{pmatrix}$
mit $r \in \mathbb{N}$,
$\mathfrak{M}_p = \begin{pmatrix} \mathfrak{o}_p & p\mathfrak{o}_p \\ p^{-1}\mathfrak{o}_p & \mathfrak{o}_p \end{pmatrix}$ und
$\Pi = \begin{pmatrix} \alpha & p \beta \\ \gamma' & \delta \end{pmatrix}$,
$\Pi U = \begin{pmatrix} \alpha \omega & p(\alpha b + \beta \overline{\omega}) \\
\gamma' \omega & p \gamma' b + \delta \overline{\omega} \end{pmatrix}$
Wir ersetzen $\mathfrak{M}_p$ durch
$\begin{pmatrix} \mathfrak{o}_p & \mathfrak{o}_p \\ \mathfrak{o}_p & \mathfrak{o}_p \end{pmatrix}$.
Wegen $\gamma' \in \mathfrak{o}_p^\times$ ist $f_{21}$ dann surjektiv.

\item
Wir zeigen die Eindeutigkeit von $\mathfrak{M}_p$.
Sei $\mathfrak{M}'_p \ne \mathfrak{M}_p$ eine $M_p$-Maximalordnung mit $\mathfrak{F}_p \subset \mathfrak{M}'_p$,
also $p$ tr\"age in $k$ und $F_p$ Divisionsalgebra.
Seien $X, Y \in \mathfrak{M}_p \smallsetminus p\mathfrak{M}_p$ mit $N(X) = p, N(Y) = p^r$,
sowie $\mathfrak{M}'_p = X\mathfrak{M}_pX^{-1}$ und $\mathfrak{N}_p = Y\mathfrak{M}_pY^{-1}$.
Sei $t \in \mathbb{Z}$ mit $Z := p^tX^{-1}Y \in \mathfrak{M}'_p \smallsetminus p\mathfrak{M}'_p$.
Aus $Z\mathfrak{M}'_pZ^{-1} = \mathfrak{N}_p$ und $N(Z) = p^{2t-1+r} \ne p^r = N(Y)$
folgt nach Lemma \ref{Item02-02}.(\ref{Item02-02.4}) also
$[\mathfrak{M}'_p : (\mathfrak{M}'_p \cap \mathfrak{N}_p)] \ne [\mathfrak{M}_p : (\mathfrak{M}_p \cap \mathfrak{N}_p)]
= [\mathfrak{F}_p : \mathfrak{G}_p]$.

\end{enumerate}

\end{proof}

\begin{theorem} \label{Item04-09}

Sei $d \in \mathbb{N}$ quadratfrei, und sei $k = \mathbb{Q}(i\sqrt{d})$ mit Hauptordnung $\mathfrak{o}$.
\\
Seien $F$ eine $\mathbb{Q}$-Quaternionenalgebra und $M$ eine Erweiterung von $F$ zur $k$-Quaternionen\-algebra.
Sei $p$ eine endliche Stelle von $\mathbb{Q}$ mit $p \nmid \Sigma_k(F)$.
Seien $\mathfrak{O}_p \subset F_p$ eine zu $\mathfrak{o}_p$ isomorphe Ordnung
und $\mathfrak{G}_p$ eine $F_p$-Ordnung mit $\mathfrak{O}_p \subset \mathfrak{G}_p$.
Bezeichne $C_p(\mathfrak{G}_p)$ die Anzahl der $M_p$-Maximalordnungen $\mathfrak{N}_p$
mit $\mathfrak{G}_p = F_p \cap \mathfrak{N}_p$.
Dann gilt:
\begin{easylist}[checklist]
\ListProperties(Style*=\large$\bullet$,Style**=)
& ${ }$ Falls $\Delta_p(\mathfrak{G}_p) = 1$ ist, so ist $C_p(\mathfrak{G}_p) = 1$.
& ${ }$ Falls $p$ zerlegt oder tr\"age in $k$ ist und $\Delta_p(\mathfrak{G}_p) \ne 1$ ist, so ist $C_p(\mathfrak{G}_p) = 2$.
& ${ }$ Falls $p$ verzweigt in $k$ ist, so gilt:
\\
${ }$ ${ }$ Falls $\Delta_p(\mathfrak{G}_p) = \pm p$ ist, so ist $C_p(\mathfrak{G}_p) = 1$.
\\
${ }$ ${ }$ Falls $\Delta_p(\mathfrak{G}_p) = + p^2$ ist, so ist $C_p(\mathfrak{G}_p) = p-1$.
\\
${ }$ ${ }$ Falls $\Delta_p(\mathfrak{G}_p) = -p^2$ ist, so ist $C_p(\mathfrak{G}_p) = p+1$.
\\
${ }$ ${ }$ Falls $p \ne 2$ und $|\Delta_p(\mathfrak{G}_p)| \ge p^3$ ist, so ist $C_p(\mathfrak{G}_p) = 2p$.
\\
${ }$ ${ }$ Falls $p = 2$ und $|\Delta_p(\mathfrak{G}_p)| \ge p^3$ ist, so wird $C_p(\mathfrak{G}_p)$ durch folgende Tabelle dargestellt:
\end{easylist}

\begin{tabular}{|r|r|r|}
\hline $|\Delta_2(\mathfrak{G}_2)|$ &
\multicolumn{2}{|c|}{$C_2(\mathfrak{G}_2)$}
\\
\hline
\hline {} &
$d \equiv 1 \pmod 4$ & $d \equiv 2 \pmod 4$
\\
\hline
\hline $8$ &
$2$ & $2$
\\
\hline $16$ &
$4$ & $4$
\\
\hline $32$ &
$8$ & $4$
\\
\hline $64$ &
$8$ & $8$
\\
\hline $\ge 128$ &
$8$ & $16$
\\
\hline
\end{tabular}

\end{theorem}

\begin{proof}

Sei $\mathfrak{F}_p$ eine $F_p$-Maximalordnung mit $\mathfrak{G}_p \subset \mathfrak{F}_p$.
\begin{easylist}[checklist]
\ListProperties(Style*=\large$\bullet$,Style**=)
& ${ }$ Sei $\Delta_p(\mathfrak{G}_p) = 1$. Dann zerf\"allt $F_p$, es gilt $\mathfrak{G}_p = \mathfrak{F}_p$,
und mit Lemma \ref{Item04-07} folgt $C_p(\mathfrak{G}_p) = 1$.
& ${ }$ Sei $p$ zerlegt in $k$, und sei $\Delta_p(\mathfrak{G}_p) \ne 1$.
\end{easylist}
\noindent
Nach Satz \ref{Item04-03}.(\ref{Item04-03.1}) gibt es $r, s \in \mathbb{Z}$ und Matrixeinheiten in $F_p$, bez\"uglich derer gilt:
\\
$\mathfrak{O}_p = \begin{pmatrix} \mathbb{Z}_p & 0 \\ 0 & \mathbb{Z}_p \end{pmatrix}$ und
$\mathfrak{G}_p = \begin{pmatrix} \mathbb{Z}_p & p^s \mathbb{Z}_p \\ p^r \mathbb{Z}_p & \mathbb{Z}_p \end{pmatrix}$.
Da $\mathfrak{O}_p$ diagonal in $\mathfrak{N}_p$ enthalten ist,
gilt $\mathfrak{N}_p \cong \begin{pmatrix} \mathbb{Z}_p & p^{-x} \mathbb{Z}_p \\ p^x \mathbb{Z}_p & \mathbb{Z}_p \end{pmatrix}
\times \begin{pmatrix} \mathbb{Z}_p & p^{-y} \mathbb{Z}_p \\ p^y \mathbb{Z}_p & \mathbb{Z}_p \end{pmatrix}$
mit $x, y \in \mathbb{Z}$.
Man rechnet nun leicht nach, dass $\mathfrak{G}_p = F_p \cap \mathfrak{N}_p$ genau dann gilt,
wenn entweder $x = r$, $y = -s$ oder $x = -s$, $y = r$.
\begin{easylist}[checklist]
\ListProperties(Style*=\large$\bullet$,Style**=)
& ${ }$ Sei $p$ tr\"age in $k$, und sei $\Delta_p(\mathfrak{G}_p) \ne 1$.
\end{easylist}
\noindent
$F_p$ ist isomorph zur $p$-Komponente einer $\mathbb{Q}$-Algebra mit Zerf\"allungsk\"orper $k$.
Es gibt also $\tau \in \mathbb{Z}$ mit den Eigenschaften laut Korollar \ref{Item03-04} sowie Matrixeinheiten in $M_p$,
bez\"uglich derer $F_p = \left\{ \left. \begin{pmatrix} a & b \\ \tau\overline{ b} & \overline{a} \end{pmatrix} \right| a, b \in k_p \right\}$.
F\"ur $r \in \mathbb{N}_0$ ist $\mathfrak{G}(r)_p :=
\left\{ \left. \begin{pmatrix} a & p^rb \\ p^r\tau\overline{ b} & \overline{a} \end{pmatrix} \right| a, b \in \mathfrak{o}_p \right\}$
eine $F_p$-Ordnung und enth\"alt die zu $\mathfrak{o}_p$ isomorphe Ordnung
$\mathfrak{O}_p := \left\{ \left. \begin{pmatrix} a & 0 \\ 0 & \overline{a} \end{pmatrix} \right| a \in \mathfrak{o}_p \right\}$.
Sei $L_p := \begin{pmatrix} k_p & 0 \\ 0 & k_p \end{pmatrix}$
mit Hauptordnung $\mathfrak{L}_p = \begin{pmatrix} \mathfrak{o}_p & 0 \\ 0 & \mathfrak{o}_p \end{pmatrix}$.
Mit $\mathfrak{N}(r)_p := \begin{pmatrix} \mathfrak{o}_p & p^r \mathfrak{o}_p \\ p^{-r} \mathfrak{o}_p & \mathfrak{o}_p \end{pmatrix}$
und $\mathfrak{N}(- r)_p := \begin{pmatrix} \mathfrak{o}_p & (p^r \tau)^{-1} \mathfrak{o}_p \\ p^r \tau \mathfrak{o}_p & \mathfrak{o}_p \end{pmatrix}$ gilt $F_p \cap \mathfrak{N}(r)_p = F_p \cap \mathfrak{N}(-r)_p = \mathfrak{G}(r)_p$.
Sei umgekehrt $\mathfrak{N}_p$ eine $M_p$-Ordnung mit $F_p \cap \mathfrak{N}_p = \mathfrak{G}(r)_p$.
Da $p$ in $k$ tr\"age ist, gilt $\mathfrak{L}_p = \mathfrak{o}_p\mathfrak{O}_p \subset \mathfrak{N}_p$.
Nach Lemma \ref{Item02-03} gibt es ein $J \in L_p^\times$ mit $\mathfrak{N}_p = JM_2(\mathfrak{o}_p)J^{-1}$.
Dann ist $\mathfrak{N}_p = \begin{pmatrix} \mathfrak{o}_p & p^{-s}\mathfrak{o}_p \\ p^s\mathfrak{o}_p & \mathfrak{o}_p \end{pmatrix}$
mit $s \in \mathbb{Z}$,
also $\mathfrak{N}_p = \mathfrak{N}(r)_p$ oder $\mathfrak{N}_p = \mathfrak{N}(-r)_p$.
Nach Lemma \ref{Item04-04} ist $\mathfrak{G}_p$ isomorph zu $\mathfrak{G}(r)_p$ mit $r \in \mathbb{N}_0$,
und $\mathfrak{N}(r)_p = \mathfrak{N}(-r)_p$ gilt nur dann, wenn $r=0$ und $F_p$ zerlegt ist.
\begin{easylist}[checklist]
\ListProperties(Style*=\large$\bullet$,Style**=)
& ${ }$ Sei $p$ verzweigt in $k$, sei $\pi \in \mathfrak{o}_p$ ein Primelement, und sei $\Delta_p(\mathfrak{G}_p) \ne 1$.
\end{easylist}
\noindent
Sei $\mathfrak{F}_p$ eine $F_p$-Maximalordnung
und $\mathfrak{M}_p$ die $M_p$-Maximalordnung mit $\mathfrak{F}_p \subset \mathfrak{M}_p$.
Wie im Beweis von Lemma \ref{Item04-08} folgt,
dass es $b, \alpha, \beta, \delta \in \mathfrak{o}_p$, $\gamma \in \mathfrak{o}_p^\times$ und Matrixeinheiten in $M_p$ gibt,
bez\"uglich derer 
$\left\{ 1, \Pi:= \begin{pmatrix} \pi & b \\ 0 & \overline{\pi} \end{pmatrix},
\Omega := \begin{pmatrix} \alpha & \beta \\ \gamma & \delta \end{pmatrix}, \Pi\Omega \right\}$
eine $\mathbb{Z}_p$-Basis von $\mathfrak{F}_p$ ist,
\\
$p$ tr\"age in $K'_p = \mathbb{Q}_p(\Omega)$ mit Hauptordnung $ \mathfrak{O}'_p = \mathbb{Z}_p[\Omega]$ ist,
und $\mathfrak{M}_p = \begin{pmatrix} \mathfrak{o}_p & \mathfrak{o}_p \\ \mathfrak{o}_p & \mathfrak{o}_p \end{pmatrix}$ gilt.
Sei $\mathfrak{O}_p := \mathbb{Z}_p[\Pi]$
und $\mathfrak{G}(r)_p := \mathfrak{O}_p + \mathfrak{O}_p \Pi^r \Omega$ f\"ur $r \in \mathbb{N}$.
Es gen\"ugt, die  Anzahl der $M_p$-Maximalordnungen $\mathfrak{N}_p$ mit $\mathfrak{G}(r)_p = F_p \cap \mathfrak{N}_p$ zu berechnen,
siehe Lemma \ref{Item04-02}.(\ref{Item04-02.1}) und \ref{Item04-04}.
\\
Sei zun\"achst $\mathfrak{N}_p$ eine $M_p$-Maximalordnung mit $\mathfrak{G}(r)_p = F_p \cap \mathfrak{N}_p$.
Nach Lemma \ref{Item04-08}.(\ref{Item04-08.2}) und den Lemmata \ref{Item04-06}.(\ref{Item04-06.1}), \ref{Item04-07} ist dann
$[\mathfrak{F}_p : \mathfrak{G}(r)_p] = [\mathfrak{M}_p : (\mathfrak{M}_p \cap \mathfrak{N}_p)] = p^r$.
Nach Lemma \ref{Item02-02}.(\ref{Item02-02.1}) gibt es $J \in \mathfrak{M}_p \smallsetminus \pi\mathfrak{M}_p$
mit $J\mathfrak{M}_pJ^{-1} = \mathfrak{N}_p$ und $N(J) \in \pi^r\mathfrak{o}_p^\times$.
\\
F\"ur jedes $m \in \mathbb{N}_0$ w\"ahlen wir eine Repr\"asentantenmenge
$\mathfrak{r}_p(m) \subset \mathfrak{o}_p$ von $\mathfrak{o}_p/\pi^m\mathfrak{o}_p$.
Mit \cite[Theorem 6.5.3.3]{Ma} folgt dann leicht,
dass wir $J = J(m,n,c) := \begin{pmatrix} \pi^m & c \\ 0 & \overline{\pi}^n \end{pmatrix}$ annehmen k\"onnen,
mit $m, n \in \mathbb{N}_0$, $m + n = r$, $c \in \mathfrak{o}_p$, f\"ur die dar\"uberhinaus gilt:
\begin{enumerate}[(C1)]
\item
$c \in \mathfrak{r}_p(m)$, und $c \in \mathfrak{o}_p^\times$, falls $mn > 0$.
\end{enumerate}
\noindent
Die Abbildung mit Zuordnungsvorschrift $\mathfrak{N}_p \mapsto (m, n, c)$ ist wohldefiniert und injektiv,
und man \"uberpr\"uft leicht, dass $\Pi \in \mathfrak{N}_p$ oder gleichwertig $J^{-1} \Pi J \in \mathfrak{M}_p$ \"aquivalent ist zu:
\begin{enumerate}[(C2)]
\item
$(\pi - \overline{\pi})c + \overline{\pi}^nb \in \pi^m\mathfrak{o}_p$.
\end{enumerate}
\noindent
Wir berechnen f\"ur jedes $(m,n) \in \mathbb{N}_0 \times \mathbb{N}_0$ die Anzahl C(m,n) der $c \in \mathfrak{o}_p$,
die die Bedingungen (C1) und (C2) erf\"ullen. Dazu ben\"otigen wir drei Vorbemerkungen:
\begin{enumerate}[(P1)]

\item

Falls $p \ne 2$, gilt $(\pi - \overline{\pi}) \in \pi\mathfrak{o}_p^\times$.
\\
Falls $p = 2$ und $d \equiv 1 \pmod{4}$, gilt $(\pi - \overline{\pi}) \in \pi^2\mathfrak{o}_p^\times$.
\\
Falls $p = 2$ und $d \equiv 2 \pmod{4}$, gilt $(\pi - \overline{\pi}) \in \pi^3\mathfrak{o}_p^\times$.
\\
Zum Beweis ist nur zu beachten, dass wir $\pi - \overline{\pi} = 2i\sqrt{d}$ annehmen k\"onnen.

\item

Genau dann gilt $b \in \mathfrak{o}_p^\times$, wenn $F_p$ zerf\"allt.
\\
Zum Beweis sei $b \in \pi\mathfrak{o}_p$. Dann ist $\Pi/\pi \in \mathfrak{M}_p$,
also $\mathfrak{o}_p\mathfrak{F}_p \subsetneqq \mathfrak{M}_p$, also $\Delta_p(\mathfrak{F}_p) = -p$.
Sei umgekehrt $\Delta_p(\mathfrak{F}_p) = -p$, 
und sei $\mathfrak{M}'_p$ der $\mathfrak{o}_p$-Modul mit Basis $\left\{ 1, \Pi/\pi, \Omega, \Pi\Omega/\pi \right\}$.
Multiplikation von links mit $(\Pi/\pi)$ oder von rechts mit $\Omega$ \"uberf\"uhrt $\mathfrak{M}'_p$ in sich.
Mit $\Omega\Pi/\pi = (\Pi/\pi)(\Pi^{-1}\Omega\Pi) \in  (\Pi/\pi)\mathfrak{F}_p \subset \mathfrak{M}'_p$ folgt,
dass $\mathfrak{M}'_p$ multiplikativ abgeschlossen, also $M_p$-Maximalordnung ist.
Dann ist notwendig $\mathfrak{M}'_p = \mathfrak{M}_p$, also $b \in \pi\mathfrak{o}_p$.

\item

Falls $p = 2$, gilt $b \not\in \pi^2\mathfrak{o}_p$.
\\
Zum Beweis nehmen wir $b = \pi^2b'$ mit $b' \in \mathfrak{o}_p$ an.
Es gibt $x, u, v, w \in \mathfrak{o}_p$ mit $\begin{pmatrix} 0 & 1 \\ 0 & 0 \end{pmatrix} = x + u\Pi/\pi + v\Omega + w\Pi\Omega/\pi$.
Daraus folgt
$0 = v\gamma + w\gamma\overline{\pi}/\pi$, also $v+w \in \pi\mathfrak{o}_p$,
und $1 = (u + w\delta)\pi b' + (v + w)\beta \in \pi\mathfrak{o}_p$, Widerspruch.

\end{enumerate}

\noindent
Nun k\"onnen wir die $C(m,n)$ elementar berechnen.
Wir f\"uhren die Details hier nicht aus, sondern fassen nur die Ergebnisse in den beiden folgenden Tabellen zusammen.
\\[4pt]
\begin{tabular}{|c|c|c|}
\hline $m, n$ &
\multicolumn{2}{|c|}{$C(m,n)$, falls $p \ne 2$}
\\
\hline
\hline {} &
$\Delta_p(\mathfrak{F}_p) = 1$ & $\Delta_p(\mathfrak{F}_p) = -p$
\\
\hline
\hline $m = 0, n \ge 0$ &
$1$ & $1$
\\
\hline $m \ge 1, n = 0$ &
$0$ & $p$
\\
\hline $m = 1, n \ge 1$ &
$p-1$ & $p-1$
\\
\hline $m \ge 2, n = 1$ &
$p$ & $0$
\\
\hline $m \ge 2, n \ge 2$ &
$0$ & $0$
\\
\hline
\end{tabular}
\\[4pt]
\begin{tabular}{|c|r|r|r|r|}
\hline $m, n$ &
\multicolumn{4}{|c|}{$C(m,n)$, falls $p = 2$}
\\
\hline
\hline {} &
\multicolumn{2}{|c|}{$d \equiv 1 \pmod 4$} & \multicolumn{2}{|c|}{$d \equiv 2 \pmod 4$}
\\
\cline{2-5} {} &
$\Delta_p(\mathfrak{F}_p) = 1$ & $\Delta_p(\mathfrak{F}_p) = -p$ & $\Delta_p(\mathfrak{F}_p) = 1$ & $\Delta_p(\mathfrak{F}_p) = -p$
\\
\hline
\hline $m = 0, n \ge 0$ &
$1$ & $1$ & $1$ & $1$
\\
\hline $m = 1, n = 0$ &
$0$ & $p = 2$ & $0$ & $p = 2$
\\
\hline $m = 1, n \ge 1$ &
$p-1= 1$ & $p-1 = 1$ & $p-1 = 1$ & $p-1 = 1$
\\
\hline $m \ge 2, n = 0$ &
$0$ & $0$ & $0$ & $0$
\\
\hline $m = 2, n = 1$ &
$0$ & $p^2-p = 2$ & $0$ & $p^2-p = 2$
\\
\hline $m = 2, n \ge 2$ &
$p^2-p = 2$ & $p^2-p = 2$ & $p^2-p = 2$ & $p^2-p = 2$
\\
\hline $m \ge 3, n = 1$ &
$0$ & $p^2 = 4$ & $0$ & $0$
\\
\hline $m = 3, n = 2$ &
$p^2 = 4$ & $0$ & $0$ & $p^3-p^2 = 4$
\\
\hline $m = 3, n \ge 3$ &
$0$ & $0$ & $p^3-p^2 = 4$ & $p^3-p^2 = 4$
\\
\hline $m \ge 4, n = 2$ &
$p^2 = 4$ & $0$ & $0$ & $p^3 = 8$
\\
\hline $m \ge 4, n = 3$ &
$0$ & $0$ & $p^3 = 8$ & $0$
\\
\hline $m \ge 4, n \ge 4$ &
$0$ & $0$ & $0$ & $0$
\\
\hline
\end{tabular}
\\[4pt]
Seien nun umgekehrt $m, n \in \mathbb{N}_0$ mit $m + n = r$ und $c \in \mathfrak{o}_p$, f\"ur die (C1) und (C2) gelten,
und seien $\mathfrak{N}_p := J(m,n,c)\mathfrak{M}_pJ(m,n,c)^{-1}$ und $\mathfrak{H}_p := F_p \cap \mathfrak{N}_p$.
Dann gilt $\mathfrak{O}_p \subset \mathfrak{H}_p$.
\\
Falls $\mathfrak{H}_p \subset \mathfrak{F}_p$,
folgt mit Lemma \ref{Item04-08}.(\ref{Item04-08.2}) und \ref{Item02-02}, dass
$[\mathfrak{F}_p : \mathfrak{H}_p] = [\mathfrak{M}_p : (\mathfrak{M}_p \cap \mathfrak{N}_p)] = p^r$,
mit Lemma \ref{Item04-02}.(\ref{Item04-02.1}) also $\mathfrak{H}_p = \mathfrak{G}(r)_p$.
Falls $\mathfrak{H}_p \not\subset \mathfrak{F}_p$,
folgt mit Lemma \ref{Item04-02}.(\ref{Item04-02.4}) und \ref{Item04-06}.(\ref{Item04-06.1}),
dass $\mathfrak{H}_p = \Pi\mathfrak{F}_p\Pi^{-1}$ und $\mathfrak{N}_p = \Pi\mathfrak{M}_p\Pi^{-1}$.
O.B.d.A. sei $b \in \mathfrak{r}_p(1)$.
Wegen $\Pi = J(1,1,b)$ k\"onnen wir zusammenfassen:
\noindent
Sei $[\mathfrak{F}_p : \mathfrak{G}_p] = p^r$.
Dann ist $C_p(\mathfrak{G}_p)$ gleich der Summe der $C(m,n)$ mit $m, n \in \mathbb{N}_0$ und $m + n = r$,
vermindert um $1$, falls $F_p$ zerf\"allt und $r = 2$ ist.
\\
Zum Schluss addieren wir die passenden Tabelleneintr\"age unter $\Delta_p(\mathfrak{G}_p) = \Delta_p(\mathfrak{F}_p)p^r$.

\end{proof}


\section{Einbettung rationaler Quaternionenordnungen} \label{Chapter05}

\begin{definition} \label{Item05-01}

Ist $K$ ein algebraischer Zahlk\"orper, $Q$ eine $K$-Quaternionenalgebra,
\\
$\mathfrak{M}$ eine $Q$-Maximalordnung und $\mathfrak{G} \subset \mathfrak{M}$ eine $Q$-Ordnung,
dann sei $\Lambda(\mathfrak{G}): = [\mathfrak{M}:\mathfrak{G}]$.

\end{definition}

\begin{definition} \label{Item05-02}

Sei $F$ eine $\mathbb{Q}$-Quaternionenalgebra, und sei $\mathfrak{G}$ eine $F$-Ordnung.
\\
Ist $n^2 \mathbb{Z}$ die Diskriminante von $\mathfrak{G}$ mit $n \in \mathbb{N}$,
dann sei $\Delta(\mathfrak{G}) := +n$ oder $-n$,
\\
je nachdem, ob $F$ an der Stelle $\infty$ zerf\"allt oder verzweigt ist.

\end{definition}

\noindent
Wir bezeichnen $\Delta(\mathfrak{G})$ auch als reduzierte Diskriminante von $\mathfrak{G}$, siehe \cite[Definition 2.7.4]{Ma}.
\\
Es gilt $\Delta(\mathfrak{G}) = \Sigma(F) \Lambda(\mathfrak{G})$.
F\"ur eine $F$-Maximalordnung $\mathfrak{F}$ gilt speziell $\Delta(\mathfrak{F}) = \Sigma(F)$.
$\Delta(\mathfrak{G})$ ist das Produkt der $\Delta_p(\mathfrak{G}_p)$ \"uber alle Primzahlen $p \in \mathbb{N}$,
siehe \cite[Theorem 2.7.3]{Ma}.

\begin{definition} \label{Item05-03}

Sei $d \in \mathbb{N}$ quadratfrei, und sei $k = \mathbb{Q}(i\sqrt{d})$ mit Hauptordnung $\mathfrak{o}$.
\\
Sei $F$ eine $\mathbb{Q}$-Quaternionenalgebra.
Eine $F$-Ordnung $\mathfrak{G}$ hei{\ss}t $\mathfrak{o}$-kompatibel genau dann,
wenn f\"ur alle Stellen $p$ von $\mathbb{Q}$ mit $p \mid \Lambda(\mathfrak{G})$ gilt:
$\mathfrak{G}_p$ enth\"alt eine zu $\mathfrak{o}_p$ isomorphe Ordnung.

\end{definition}

\begin{theorem} \label{Item05-04}

Sei $d \in \mathbb{N}$ quadratfrei, und sei $k = \mathbb{Q}(i\sqrt{d})$ mit Hauptordnung $\mathfrak{o}$ und Idealgruppe $I$.
Sei $F$ eine $\mathbb{Q}$-Quaternionenalgebra, und sei $\mathfrak{F}$ eine $F$-Maximalordnung.
\begin{enumerate}[(i)]

\item \label{Item05-04.1}
Sei $\mathfrak{G} \subset \mathfrak{F}$ eine $\mathfrak{o}$-kompatible $F$-Ordnung.
\\
Dann ist $\Lambda(\mathfrak{G})$ teilerfremd zu $\Sigma_k(F)$, und es gilt $\Lambda(\mathfrak{G}) \in N_{abs}(I)$.

\item \label{Item05-04.2}
Sei $\lambda \in \mathbb{N}$ teilerfremd zu $\Sigma_k(F)$, und sei $\lambda \in N_{abs}(I)$.
\\
Dann gibt es eine $\mathfrak{o}$-kompatible $F$-Ordnung $\mathfrak{G} \subset \mathfrak{F}$ mit $\Lambda(\mathfrak{G}) = \lambda$.

\item \label{Item05-04.3}
Seien $\mathfrak{G}, \mathfrak{G}'$ zwei $\mathfrak{o}$-kompatible $F$-Ordnungen mit $\Lambda(\mathfrak{G}) = \Lambda(\mathfrak{G}')$.
\\
F\"ur alle endlichen Stellen $p$ von $\mathbb{Q}$ gilt dann: $\mathfrak{G}_p$ ist isomorph zu $\mathfrak{G}'_p$.

\end{enumerate}

\end{theorem}

\noindent
$Bemerkung$: Wir haben die Struktur der $\mathfrak{G}_p$ in den S\"atzen \ref{Item04-02}, \ref{Item04-03} und \ref{Item04-05} dargelegt.

\begin{proof}

jeweils durch Zusammenf\"ugen der Ergebnisse f\"ur alle endlichen Stellen $p$ von $\mathbb{Q}$.

\begin{enumerate}[(i)]

\item[(\ref{Item05-04.1})]
Falls $p \mid \Lambda(\mathfrak{G})$, l\"asst sich $\mathfrak{o}_p$ in $\mathfrak{G}_p$, also $k_p$ in $F_p$ einbetten.
Wenn dann auch $p \mid \Sigma_k(F)$ gelten w\"urde, so w\"are $F_p$ eine Divisionsalgebra, aber $k_p$ kein K\"orper, Widerspruch.
\\
Falls $p \nmid \Sigma_k(F)$,
gibt es nach den S\"atzen \ref{Item04-02} und \ref{Item04-03} ein $r \in \mathbb{N}_0$
mit $[\mathfrak{F}_p : \mathfrak{G}_p] = p^r$, falls $p$ in $k$ verzweigt oder zerlegt ist,
und mit $[\mathfrak{F}_p : \mathfrak{G}_p] = p^{2r}$, falls $p$ in $k$ tr\"age ist.

\item[(\ref{Item05-04.2})]
Falls $p \nmid \lambda$, sei $\mathfrak{G}_p = \mathfrak{F}_p$.
\\
Falls $p \mid \lambda$, sei $\lambda_p$ der $p$-Anteil von $\lambda$.
Dann gibt es ein $r \in \mathbb{N}$ mit $\lambda_p = p^r$ oder mit $\lambda_p = p^{2r}$,
je nachdem, ob $p$ in $k$ verzweigt oder zerlegt ist, oder ob $p$ in $k$ tr\"age ist.
Wegen $p \nmid \Sigma_k(F)$ enth\"alt $\mathfrak{F}_p$ eine zu $\mathfrak{o}_p$ isomorphe Ordnung $\mathfrak{O}_p$.
Nach den S\"atzen \ref{Item04-02} und \ref{Item04-03} gibt es eine $F_p$-Ordnung $\mathfrak{G}_p$
mit $\mathfrak{O}_p \subset \mathfrak{G}_p \subset \mathfrak{F}_p$ und $[\mathfrak{F}_p : \mathfrak{G}_p] = \lambda_p$.

\item[(\ref{Item05-04.3})]
Falls $p \mid \Sigma_k(F)$, ist $\mathfrak{G}_p = \mathfrak{G}'_p$ die einzige $F_p$-Maximalordnung.
\\
Falls $p \nmid \Sigma_k(F)$, gibt es nach Lemma \ref{Item04-04} einen Isomorphismus $\mathfrak{G}_p \rightarrow \mathfrak{G}'_p$.

\end{enumerate}

\end{proof}

\begin{lemma} \label{Item05-05}

Sei $d \in \mathbb{N}$ quadratfrei, und sei $k = \mathbb{Q}(i\sqrt{d})$.
\\
Sei $F$ eine $\mathbb{Q}$-Quaternionenalgebra,
sei $M$ eine Erweiterung von $F$ zur $k$-Quaternionen\-algebra,
und seien $\mathfrak{N}, \mathfrak{N}'$ zwei $M$-Maximalordnungen.
Dann gilt:
\begin{enumerate}[(i)]

\item \label{Item05-05.1}
$\Lambda(\mathfrak{N} \cap \mathfrak{N}') = N_{abs}(N(\mathfrak{N}\mathfrak{N}'))^{-1}$.

\item \label{Item05-05.2}
Sei  $\mathfrak{M}$ eine $M$-Maximalordnung. Dann gilt
$\Lambda(\mathfrak{N} \cap \mathfrak{N}') \Lambda(\mathfrak{N} \cap \mathfrak{M}) \Lambda(\mathfrak{M} \cap \mathfrak{N}') \in \mathbb{Q}^{\times(2)}$.

\item \label{Item05-05.3}
$\mathfrak{N}$ ist genau dann isomorph zu $\mathfrak{N}'$, wenn es ein $f \in \mathbb{N}$ gibt mit:
\\
$f \mid \Sigma_k(F)$ und $\left( \dfrac{f \Lambda(\mathfrak{N} \cap \mathfrak{N}'),-d} {p} \right) = 1$ f\"ur alle Stellen $p$ von $\mathbb{Q}$.

\end{enumerate}

\end{lemma}

\begin{proof} ${}$
\begin{enumerate}[(i)]

\item[(\ref{Item05-05.1})]
Beweis durch Zusammenf\"ugen des folgenden Ergebnisses f\"ur alle endlichen Stellen $p$ von $\mathbb{Q}$.
Seien dazu $\mathfrak{p} \mid p$ eine Stelle von $k$ und $\pi \in \mathfrak{o}_\mathfrak{p}$ ein Primelement.
Falls $\mathfrak{p} \mid \Sigma_k(F)$, gilt $\mathfrak{N}_\mathfrak{p} = \mathfrak{N}'_\mathfrak{p}$,
also $\Lambda(\mathfrak{N} \cap \mathfrak{N}')_\mathfrak{p} = 1$
und $N(\mathfrak{N}\mathfrak{N}')_\mathfrak{p} = \mathfrak{o}_\mathfrak{p}$.
Falls $\mathfrak{p} \nmid \Sigma_k(F)$, gibt es nach Lemma \ref{Item02-02} ein $r \in \mathbb{N}_0$ mit
$[\mathfrak{N}_\mathfrak{p} : (\mathfrak{N} \cap \mathfrak{N}')_\mathfrak{p}] =
N_{abs}(\pi\mathfrak{o}_\mathfrak{p})^r = N_{abs}(N(\mathfrak{N}\mathfrak{N}')_\mathfrak{p})^{-1}$.

\item[(\ref{Item05-05.2})]
$\mathfrak{A} := (\mathfrak{N} \mathfrak{M} \mathfrak{N}'\mathfrak{N})^{-1} =
(\mathfrak{N} \mathfrak{M} \mathfrak{M} \mathfrak{N}' \mathfrak{N}' \mathfrak{N})^{-1}$
ist ein zweiseitiges ganzes $\mathfrak{N}$-Ideal und teilerfremd zu $\Sigma_k(F)$.
Daher gibt es ein ganzes $\mathfrak{o}$-Ideal $\mathfrak{a}$ mit $\mathfrak{A} = \mathfrak{N}\mathfrak{a}$, siehe dazu etwa \cite[Lemma 6.7.5]{Ma}.
Dann gilt $\mathfrak{a}^2 = N(\mathfrak{A}) = N(\mathfrak{N}'\mathfrak{N})^{-1} N(\mathfrak{M}\mathfrak{N}')^{-1} N(\mathfrak{N}\mathfrak{M})^{-1}$, nach (\ref{Item05-05.1}) also
$\Lambda(\mathfrak{N}' \cap \mathfrak{N}) \Lambda(\mathfrak{M} \cap \mathfrak{N}') \Lambda(\mathfrak{N} \cap \mathfrak{M}) =
N_{abs}(\mathfrak{a})^2 \in \mathbb{Q}^{\times(2)}$.

\item[(\ref{Item05-05.3})]
Bezeichne $R$ die von den Verzweigungsstellen von $M$ erzeugte $\mathfrak{o}$-Idealgruppe.
\\
Nach \cite[Chapter 6.7]{Ma} ist $\mathfrak{N}$ genau dann isomorph zu $\mathfrak{N}'$,
wenn $N(\mathfrak{N}\mathfrak{N}') \in RI^{(2)}H$,
also wenn es ein o.B.d.A. ganzes $\mathfrak{f} \in R$ gibt mit $\mathfrak{f}N(\mathfrak{N}\mathfrak{N}')^{-1} \in I^{(2)}H$.
\\
Nach \cite[Kapitel III, § 8, S\"atze 6 und 7]{Bo} gilt $\mathfrak{f}N(\mathfrak{N}\mathfrak{N}')^{-1} \in I^{(2)}H$ genau dann, wenn
$N_{abs}(\mathfrak{f}) \Lambda(\mathfrak{N} \cap \mathfrak{N}') = N_{abs}(\mathfrak{f}N(\mathfrak{N}\mathfrak{N}')^{-1})
\in N_{abs}(k^\times)$,
nach \cite[Kapitel I, § 7, Satz 1]{Bo} (Minkowski-Hasse) also genau dann, wenn
$\left(\dfrac{N_{abs}(\mathfrak{f}) \Lambda(\mathfrak{N} \cap \mathfrak{N}'),-d}{p}\right) = 1$ f\"ur alle $p$.
\\
Sei $f$ der quadratfreie Anteil von $N_{abs}(\mathfrak{f})$. Nach Lemma \ref{Item03-02}.(\ref{Item03-02.1}) gilt dann $f \mid \Sigma_k(F)$.
Ist umgekehrt $f \in \mathbb{N}$ mit $f \mid \Sigma_k(F)$, so gibt es ein $\mathfrak{f} \in R$ mit $f = N_{abs}(\mathfrak{f})$.

\end{enumerate}

\end{proof}

\begin{theorem} \label{Item05-06}

Sei $d \in \mathbb{N}$ quadratfrei, und sei $k = \mathbb{Q}(i\sqrt{d})$ mit Hauptordnung $\mathfrak{o}$.
Sei $F$ eine $\mathbb{Q}$-Quaternionenalgebra, und sei $M$ eine Erweiterung von $F$ zur $k$-Quaternionenalgebra.
\begin{enumerate}[(i)]

\item \label{Item05-06.1}
Sei $\mathfrak{N}$ eine $M$-Maximalordnung.
\\
Dann ist $F \cap \mathfrak{N}$ eine $\mathfrak{o}$-kompatible $F$-Ordnung.

\item \label{Item05-06.2}
Sei umgekehrt $\mathfrak{G}$ eine $\mathfrak{o}$-kompatible $F$-Ordnung.
\\
Dann gibt es eine $M$-Maximalordnung $\mathfrak{N}$ mit $\mathfrak{G} = F \cap \mathfrak{N}$.

\item \label{Item05-06.3}
Sind $\mathfrak{N}$ eine $M$-Maximalordnung und $\mathfrak{F}$ eine $F$-Maximalordnung mit $F \cap \mathfrak{N} \subset \mathfrak{F}$,
so gibt es genau eine $M$-Maximalordnung $\mathfrak{M}$ mit $\mathfrak{F} \subset \mathfrak{M}$
und $\Lambda(F \cap \mathfrak{N}) = \Lambda(\mathfrak{M} \cap \mathfrak{N})$.
Die Inklusion $\mathfrak{F} \hookrightarrow \mathfrak{M}$ induziert einen Isomorphismus
$\mathfrak{F} / (F \cap \mathfrak{N}) \rightarrow \mathfrak{M} / (\mathfrak{M} \cap \mathfrak{N})$.

\item \label{Item05-06.4}
Seien $\mathfrak{N}, \mathfrak{N}'$ zwei $M$-Maximalordnungen.
\\
Dann gilt:
$\Delta(F \cap \mathfrak{N}) \Lambda(\mathfrak{N} \cap \mathfrak{N}') \Delta(F \cap \mathfrak{N}') \in \mathbb{Q}^{\times (2)}$.
\\
Speziell ist $\mathfrak{N}$ isomorph zu $\mathfrak{N}'$, falls $\Delta(F \cap \mathfrak{N}) = \Delta(F \cap \mathfrak{N}')$ gilt.

\end{enumerate}

\end{theorem}

\begin{proof}

jeweils durch Zusammenf\"ugen der Ergebnisse f\"ur alle endlichen Stellen $p$ von $\mathbb{Q}$.
\begin{enumerate}[(i)]

\item[(\ref{Item05-06.1})]
Falls $p \mid \Lambda(F \cap \mathfrak{N})$,
enth\"alt $F_p \cap \mathfrak{N}_p$ nach \ref{Item04-08}.(\ref{Item04-08.1}) eine zu $\mathfrak{o}_p$ isomorphe Ordnung.

\item[(\ref{Item05-06.2})]
Falls $p \mid \Sigma_k(F)$, ist $\mathfrak{G}_p$ die $F_p$-Maximalordnung.
Sei $\mathfrak{N}_p$ die $M_p$-Maximalordnung.
Falls $p \nmid \Sigma_k(F)$, gibt es nach \ref{Item04-09} eine $M_p$-Maximalordnung $\mathfrak{N}_p$
mit $\mathfrak{G}_p = F_p \cap \mathfrak{N}_p$.

\item[(\ref{Item05-06.3})]
Falls $p \mid \Sigma_k(F)$, ist $F_p \cap \mathfrak{N}_p = \mathfrak{F}_p$ die $F_p$-Maximalordnung,
und mit $\mathfrak{M}_p = \mathfrak{N}_p$ gilt 
$[\mathfrak{F}_p : (F \cap \mathfrak{N})_p] = [\mathfrak{M}_p : (\mathfrak{M}_p \cap \mathfrak{N}_p)] = 1$.
Falls $p \nmid \Sigma_k(F)$, gibt es nach \ref{Item04-08}.(\ref{Item04-08.2})
genau eine $M_p$-Maximalordnung $\mathfrak{M}_p$ mit $\mathfrak{F}_p \subset \mathfrak{M}_p$
und $[\mathfrak{F}_p : (F_p \cap \mathfrak{N}_p)] = [\mathfrak{M}_p : (\mathfrak{M}_p \cap \mathfrak{N}_p)]$.
\\
$\mathfrak{F} / (F \cap \mathfrak{N}) \rightarrow \mathfrak{M} / (\mathfrak{M} \cap \mathfrak{N})$ ist Isomorphismus,
da f\"ur alle $p \ne \infty$ ein Isomorphismus.

\item[(\ref{Item05-06.4})]
Seien $\mathfrak{F}, \mathfrak{F}'$ zwei $F$-Maximalordnungen
mit $F \cap \mathfrak{N} \subset \mathfrak{F}$ und $F \cap \mathfrak{N}' \subset \mathfrak{F}'$.
Nach (\ref{Item05-06.3}) gibt es $M$-Maximalordnungen $\mathfrak{M}, \mathfrak{M}'$
mit $\mathfrak{F} \subset \mathfrak{M}$ und $\mathfrak{F}' \subset \mathfrak{M}'$
sowie mit $\Lambda(F \cap \mathfrak{N}) = \Lambda(\mathfrak{M} \cap \mathfrak{N})$ und 
$\Lambda(F \cap \mathfrak{N}') = \Lambda(\mathfrak{M}' \cap \mathfrak{N}')$.
Mit Lemma \ref{Item05-05}.(\ref{Item05-05.2}) folgt
$\Lambda(\mathfrak{N} \cap \mathfrak{N}') \Lambda(F \cap \mathfrak{N}) \Lambda(\mathfrak{M} \cap \mathfrak{N}')
\in \mathbb{Q}^{\times(2)}$ und
$\Lambda(\mathfrak{M} \cap \mathfrak{N}') \Lambda(\mathfrak{M} \cap \mathfrak{M}') \Lambda(F \cap \mathfrak{N}')
\in \mathbb{Q}^{\times(2)}$, also
$\Delta(F \cap \mathfrak{N}) \Lambda(\mathfrak{N} \cap \mathfrak{N}') \Delta(F \cap \mathfrak{N}') \Lambda(\mathfrak{M} \cap \mathfrak{M}') \in \mathbb{Q}^{\times(2)}$.
Wir zeigen $\Lambda(\mathfrak{M} \cap \mathfrak{M}') \in \mathbb{Q}^{\times (2)}$:
\\
Falls $F$ an der Stelle $p$ verzweigt ist, gilt $\mathfrak{F}_p = \mathfrak{F}'_p$,
und nach Lemma \ref{Item04-06} bzw. \ref{Item04-07} gilt $\mathfrak{M}_p = \mathfrak{M}'_p$
oder $[\mathfrak{M}_p : (\mathfrak{M}_p \cap \mathfrak{M}'_p)] = p^2$.
\\
Falls $F$ an der Stelle p zerf\"allt, gibt es nach Lemma \ref{Item02-02}.(\ref{Item02-02.1}) ein $r \in \mathbb{N}_0$ und Matrix\-einheiten in $F_p$,
bez\"uglich derer $\mathfrak{F}_p = \begin{pmatrix} \mathbb{Z}_p & \mathbb{Z}_p \\ \mathbb{Z}_p & \mathbb{Z}_p \end{pmatrix}$ und
$\mathfrak{F}'_p = \begin{pmatrix} \mathbb{Z}_p & p^{-r}\mathbb{Z}_p \\ p^{r}\mathbb{Z}_p & \mathbb{Z}_p \end{pmatrix}$ gilt.
Damit ist $\mathfrak{M}_p = \begin{pmatrix} \mathfrak{o}_p & \mathfrak{o}_p \\ \mathfrak{o}_p & \mathfrak{o}_p \end{pmatrix}$ und
$\mathfrak{M}'_p = \begin{pmatrix} \mathfrak{o}_p & p^{-r}\mathfrak{o}_p \\ p^{r}\mathfrak{o}_p & \mathfrak{o}_p \end{pmatrix}$,
also $[\mathfrak{M}_p : (\mathfrak{M}_p \cap \mathfrak{M}'_p)] = p^{2r}$.

\end{enumerate}

\end{proof}

\begin{lemma} \label{Item05-07}

Sei $d \in \mathbb{N}$ quadratfrei, und sei $k = \mathbb{Q}(i\sqrt{d})$ mit Diskriminante $D$.
\\
Seien $E$ eine $\mathbb{Q}$-Quaternionenalgebra und $M$ eine Erweiterung von $E$ zur $k$-Quaternionen\-algebra.
Sei $T \in M^\times$ mit $\Phi_E(T) = T^*$, $N(T) \in \mathbb{Z}$, $S(T) \in 4\mathbb{Z}$ und $S(T)^2 \ne 4N(T)$.
\\
Sei  $p$ eine endliche Stelle von $\mathbb{Q}$, an der $E$ zerf\"allt. Dann gilt:
\\
Es gibt eine $E_p$-Maximalordnung $\mathfrak{E}_p$ und eine $M_p$-Maximalordnung $\mathfrak{N}_p$
mit $\mathfrak{E}_p \subset \mathfrak{N}_p$,
\\
so dass die $F(E,T)_p$-Ordnung $F(E,T)_p \cap \mathfrak{N}_p$ die Diskriminante $D^2d^2N(T)^2\mathbb{Z}_p$ hat.

\end{lemma}

\begin{proof}

Es gibt Matrixeinheiten in $E_p$, bez\"uglich derer
$T = \begin{pmatrix} \sigma + \alpha i\sqrt{d} & \beta i\sqrt{d} \\ \gamma i\sqrt{d} & \sigma - \alpha i\sqrt{d} \end{pmatrix}$
mit $\sigma, \alpha, \beta, \gamma \in \mathbb{Q}_p$.
Sei $\tau := N(T) - S(T)^2/4 \in \mathbb{Z}$.
Wegen $\tau \ne 0$ gibt es $a, b \in \mathbb{Q}_p$ mit $j := a^2\beta - 2ab\alpha - b^2\gamma \ne 0$.
Sei $J := \begin{pmatrix} a/d & b/d \\ a \alpha + b \gamma & a \beta - b \alpha \end{pmatrix}$.
Dann ist $N(J) = j/d$ und $JTJ^{-1} = \begin{pmatrix} \sigma & i\sqrt{d}/d \\ \tau i\sqrt{d} & \sigma \end{pmatrix}$.
Also gibt es Matrixeinheiten in $E_p$ bez\"uglich derer $T = \begin{pmatrix} \sigma & i\sqrt{d}/d \\ \tau i\sqrt{d} & \sigma \end{pmatrix}$ gilt.
Diesbez\"uglich sei $\mathfrak{E}_p = M_2(\mathbb{Z}_p)$ und $\mathfrak{N}_p = M_2(\mathfrak{o}_p)$. Dann ist
$\left\{ 1, U := \begin{pmatrix} 0 & 1 \\ \tau d & 0 \end{pmatrix},
V := \begin{pmatrix} \sigma d & -i\sqrt{d} \\ \tau d i\sqrt{d} & -\sigma d \end{pmatrix},
W := \begin{pmatrix} i\sqrt{d} & 0 \\ 2\sigma d & -i\sqrt{d} \end{pmatrix} = (\sigma U + UV/d)/\tau \right\}$
\\
$\mathbb{Z}_p$-Basis einer $F(E,T)_p$-Ordnung $\mathfrak{G}_p \subset \mathfrak{N}_p$
mit Diskriminante $16 d^4 N(T)^2\mathbb{Z}_p$.
Ist $p=2$ in $k$ verzweigt oder $p \ne 2$, so gilt $\mathfrak{G}_p = F(E,T)_p \cap \mathfrak{N}_p$.
Ist $p = 2$ in $k$ tr\"age oder zerlegt,
so ist $\{ 1, (U + V)/2, V, (1 + W)/2 \}$
$\mathbb{Z}_p$-Basis von $F(E,T)_p \cap \mathfrak{N}_p$
mit Diskriminante $d^4 N(T)^2\mathbb{Z}_p$.

\end{proof}

\begin{lemma} \label{Item05-08}

Sei $d \in \mathbb{N}$ quadratfrei, und sei $k = \mathbb{Q}(i\sqrt{d})$.
\\
Sei $E$ eine $\mathbb{Q}$-Quaternionenalgebra mit $|\Sigma(E)| = \Sigma_k(E)$,
sei $M$ eine Erweiterung von $E$ zur $k$-Quaternionenalgebra,
und sei $F$ eine weitere $\mathbb{Q}$-Quaternionenalgebra mit $F \subset M$.
\\
Dann gibt es ein $e \in \mathbb{N}$ mit $e \mid \Sigma_k(E)$,
eine $E$-Maximalordnung $\mathfrak{E}$ und eine $F$-Maximal\-ordnung $\mathfrak{F}$,
sowie $M$-Maximalordnungen $\mathfrak{N}$ und $\mathfrak{L}$ mit $\mathfrak{E} \subset \mathfrak{N}$ und $\mathfrak{F} \subset \mathfrak{L}$,
so dass gilt:
\\
$\left(\dfrac{\Sigma(F) e \Lambda(\mathfrak{L} \cap \mathfrak{N}) \Sigma(E),-d}{p}\right) =
\left(\dfrac{E}{p}\right) \left(\dfrac{F}{p}\right)$ f\"ur alle Stellen $p$ von $\mathbb{Q}$.

\end{lemma}

\begin{proof}

Nach Satz \ref{Item03-03}.(\ref{Item03-03.3}) gibt es ein $T \in M^\times$ mit $\Phi_E(T) = T^*$ und $F = F(E,T)$.
\\
Falls $S(T)^2 = 4N(T)$, gilt $N(T) \in \mathbb{Q}^{\times(2)}$.
Nach Satz \ref{Item03-03}.(\ref{Item03-03.2}) ist $F$ isomorph zu $E$, also $\Sigma(F) = \Sigma(E)$.
Seien $\mathfrak{E},\mathfrak{F}$ Maximalordnungen in $E,F$,
und seien $\mathfrak{N},\mathfrak{L}$ zwei $M$-Maximal\-ordnungen 
mit $\mathfrak{E} \subset \mathfrak{N}$ und $\mathfrak{F} \subset \mathfrak{N}$,
so dass $\mathfrak{E}$ zu $\mathfrak{F}$ und $\mathfrak{N}$ zu $\mathfrak{L}$ isomorph ist.
Nach Lemma \ref{Item05-05}.(\ref{Item05-05.3}) gibt es $e \in \mathbb{N}$ mit $e \mid \Sigma_k(E)$ und
$\left( \dfrac{e \Lambda(\mathfrak{L} \cap \mathfrak{N}),-d} {p} \right) = 1 = \left(\dfrac{E}{p}\right) \left(\dfrac{F}{p}\right)$
f\"ur alle $p$.
\\
Wir k\"onnen nun also $S(T)^2 \ne 4N(T)$ und o.B.d.A. $N(T) \in \mathbb{Z}$, $S(T) \in 4 \mathbb{Z}$ annehmen.
Falls $p \mid \Sigma_k(E)$, ist $F_p \cap \mathfrak{N}_p$ die $F_p$-Maximalordnung, hat also die Diskriminante $p^2\mathbb{Z}_p$.
Nach Lemma \ref{Item05-07} gibt es eine $E$-Maximalordnung $\mathfrak{E}$ und eine $M$-Maximalordnung $\mathfrak{N}$ mit $\mathfrak{E} \subset \mathfrak{N}$, so dass die $F_p$-Ordnung $F_p \cap \mathfrak{N}_p$
an allen endlichen Stellen $p$ mit $p \nmid \Sigma_k(E)$ die Diskriminante $D^2 d^2 N(T)^2\mathbb{Z}_p$ hat.
Sei $D d N(T) = tg$ mit $t \in \mathbb{Z}$ und $g \in \mathbb{N}$,
so dass $t$ teilerfremd zu $\Sigma_k(E)$ ist und alle Primteiler von $g$ auch Teiler von $\Sigma_k(E)$ sind.
Dann hat $F \cap \mathfrak{N}$ die Diskriminante $(D d N(T) \Sigma_k(E) / g)^2 \mathbb{Z}$.
Andererseits hat $F \cap \mathfrak{N}$ die Diskriminante $(\Sigma(F) \Lambda(F \cap  \mathfrak{N}))^2 \mathbb{Z}$.
Unter Beachtung des Vorzeichens (Satz \ref{Item03-03}.(\ref{Item03-03.2}) f\"ur $p = \infty$) gilt also
$|D| d N(T)\Sigma(E) = g  \Sigma(F) \Lambda(F \cap \mathfrak{N})$.
Sei nun $\mathfrak{F}$ eine $F$-Maximalordnung mit $F \cap \mathfrak{N} \subset \mathfrak{F}$.
Nach Satz \ref{Item05-06}.(\ref{Item05-06.3}) gibt es eine $M$-Maximalordnung $\mathfrak{L}$ mit $\mathfrak{F} \subset \mathfrak{L}$
und $\Lambda(F \cap \mathfrak{N}) = \Lambda(\mathfrak{L} \cap \mathfrak{N})$.
Sei $e \in \mathbb{N}$ der quadratfreie Anteil von $g$.
Die Behauptung folgt jetzt mit Satz \ref{Item03-03}.(\ref{Item03-03.2}).

\end{proof}

\begin{theorem} \label{Item05-09}

Sei $d \in \mathbb{N}$ quadratfrei, und sei $k = \mathbb{Q}(i\sqrt{d})$.
\\
Seien $F, F'$ zwei $\mathbb{Q}$-Quaternionenalgebren,
sei $M$ eine gemeinsame Erweiterung von $F$ und $F'$ zur $k$-Quaternionenalgebra,
und seien $\mathfrak{M}, \mathfrak{M}'$ zwei $M$-Maximalordnungen.
\\
Dann gibt es ein $f \in \mathbb{N}$ mit $f \mid \Sigma_k(F)$,
so dasss f\"ur alle Stellen $p$ von $\mathbb{Q}$ gilt:
\\
$\left(\dfrac{\Delta(F \cap \mathfrak{M}) f \Lambda(\mathfrak{M} \cap \mathfrak{M}') \Delta(F' \cap \mathfrak{M}'),-d}{p}\right)
= \left(\dfrac{F}{p}\right) \left(\dfrac{F'}{p}\right)$.

\end{theorem}

\begin{proof}

Sei $\mathfrak{F}'$ eine $F$-Maximalordnung mit $F \cap \mathfrak{M} \subset \mathfrak{F}'$.
Nach Satz \ref{Item05-06}.(\ref{Item05-06.3}) gibt es eine $M$-Maximalordnung $\mathfrak{L}'$
mit $\mathfrak{F}' \subset \mathfrak{L}'$ und $\Lambda(F \cap \mathfrak{M}) = \Lambda(\mathfrak{L}' \cap \mathfrak{M})$,
und nach Lemma \ref{Item03-02} gibt es eine $\mathbb{Q}$-Quaternionenalgebra $E \subset M$ mit $|\Sigma(E)| = \Sigma_k(E)$.
Seien $e, \mathfrak{E}, \mathfrak{F}, \mathfrak{N}$ und $\mathfrak{L}$ wie in Lemma \ref{Item05-08}.
Mit Lemma \ref{Item05-05}.(\ref{Item05-05.2}) folgt
$\Lambda(\mathfrak{L} \cap \mathfrak{N}) \Lambda(\mathfrak{L} \cap \mathfrak{L}') \Lambda(\mathfrak{L}' \cap \mathfrak{M})
\Lambda(\mathfrak{M} \cap \mathfrak{N}) \in \mathbb{Q}^{\times(2)}$,
und nach Satz \ref{Item05-06}.(\ref{Item05-06.4}) gilt $\Lambda(\mathfrak{L} \cap \mathfrak{L}') =
\Lambda(F \cap \mathfrak{L}) \Lambda(\mathfrak{L} \cap \mathfrak{L}') \Lambda(F \cap \mathfrak{L}') \in \mathbb{Q}^{\times (2)}$.
Mit Lemma \ref{Item05-08} und $\Sigma(F) \Lambda(F \cap \mathfrak{M}) = \Delta(F \cap \mathfrak{M})$ folgt
$\left(\dfrac{\Delta(F \cap \mathfrak{M}) e \Lambda(\mathfrak{M} \cap \mathfrak{N}) \Sigma(E),-d}{p}\right) = \left(\dfrac{E}{p}\right) \left(\dfrac{F}{p}\right)$
\\
f\"ur alle Stellen $p$ von $\mathbb{Q}$.
Entsprechend gibt es ein $e' \in \mathbb{N}$ mit $e' \mid \Sigma_k(E)$, eine $E$-Maximalordnung $\mathfrak{E}'$
und eine $M$-Maximalordnung $\mathfrak{N}'$ mit $\mathfrak{E}' \subset \mathfrak{N}'$,
so dass f\"ur alle Stellen $p$ von $\mathbb{Q}$ gilt:
$\left(\dfrac{\Delta(F' \cap \mathfrak{M}') e' \Lambda(\mathfrak{M}' \cap \mathfrak{N}') \Sigma(E),-d}{p}\right) =
\left(\dfrac{E}{p}\right) \left(\dfrac{F'}{p}\right)$.
Bezeichne $f$ den quadratfreien Anteil von $ee'$.
Die Behauptung folgt durch Multiplikation der beiden Gleichungen,
denn nach Satz \ref{Item05-06}.(\ref{Item05-06.4}) gilt $\Lambda(\mathfrak{N} \cap \mathfrak{N}') =
\Lambda(E \cap \mathfrak{N}) \Lambda(\mathfrak{N} \cap \mathfrak{N}') \Lambda(E \cap \mathfrak{N}') \in \mathbb{Q}^{\times (2)}$,
und mit Lemma \ref{Item05-05}.(\ref{Item05-05.2}) folgt
$\Lambda(\mathfrak{M} \cap \mathfrak{M}') \Lambda(\mathfrak{M} \cap \mathfrak{N}) \Lambda(\mathfrak{N} \cap \mathfrak{N}') \Lambda(\mathfrak{N}' \cap \mathfrak{M}') \in \mathbb{Q}^{\times(2)}$.

\end{proof}

\begin{corollary} \label{Item05-10}

Sei $d \in \mathbb{N}$ quadratfrei, und sei $k = \mathbb{Q}(i\sqrt{d})$ mit Hauptordnung $\mathfrak{o}$.
\\
Sei $F \subset M_2(k)$ eine $\mathbb{Q}$-Quaternionenalgebra,
und sei $\mathfrak{M}$ eine $M_2(k)$-Maximalordnung.
F\"ur alle Stellen $p$ von $\mathbb{Q}$ gilt dann
$\left(\dfrac{\Delta(F \cap \mathfrak{M}) \Lambda(\mathfrak{M} \cap M_2(\mathfrak{o})),-d}{p}\right) = \left(\dfrac{F}{p}\right)$.
\\
Speziell gilt $\left(\dfrac{\Delta(F \cap M_2(\mathfrak{o})),-d}{p}\right) = \left(\dfrac{F}{p}\right)$ f\"ur alle Stellen $p$ von $\mathbb{Q}$.

\end{corollary}

\begin{proof}

Seien $F' = M_2(\mathbb{Q})$, $M = M_2(k)$ und $\mathfrak{M}' = M_2(\mathfrak{o})$.
Damit gilt $\Sigma_k(F) = 1$ und $\Delta(F' \cap \mathfrak{M'}) = 1$ sowie $\left(\dfrac{F'}{p}\right) = 1$ f\"ur alle Stellen $p$.
Die Behauptung folgt mit Satz \ref{Item05-09}.

\end{proof}


\section{Die nichtzyklischen endlichen Untergruppen} \label{Chapter06}

\begin{definition} \label{Item06-01} ${}$
\begin{enumerate}[(i)]

\item \label{Item06-01.1}
Bezeichne $P: SL_2(\mathbb{C}) \rightarrow PSL_2(\mathbb{C})$ die kanonische Projektion.
\\
F\"ur Untergruppen $\Gamma \subset SL_2(\mathbb{C})$ k\"urzen wir
$P\Gamma := P(\Gamma) \cong \Gamma / (\Gamma \cap \{ \pm 1 \})$ ab.

\item \label{Item06-01.2}
F\"ur einen Ring $R \subset M_2(\mathbb{C})$ mit Eins sei $\Gamma(R) := \{ A \in R \mid N(A) = 1 \}$

\item \label{Item06-01.3}
F\"ur eine endliche Menge $\mathcal{G} \subset PSL_2(\mathbb{C})$
seien $F(\mathcal{G}) \subset M_2(\mathbb{C})$ der von $P^{-1}(\mathcal{G})$ erzeugte $\mathbb{Q}$-Vektorraum
und $\mathfrak{F}(\mathcal{G})$ der von $P^{-1}(\mathcal{G})$ erzeugte $\mathbb{Z}$-Modul.

\item \label{Item06-01.4}
Sei $d \in \mathbb{N}$ quadratfrei, und sei $k = \mathbb{Q}(i\sqrt{d})$ mit Hauptordnung $\mathfrak{o}$.
\\
Seien $M$ eine $k$-Quaternionenalgebra und $\mathfrak{M}$ eine $M$-Maximalordnung.
\\
Dann bezeichnen wir $P\Gamma(\mathfrak{M})$ als verallgemeinerte Bianchi-Gruppe,
\\
und speziell bezeichnen wir $P\Gamma(M_2(\mathfrak{o})) = PSL_2(\mathfrak{o})$ als Bianchi-Gruppe.

\end{enumerate}

\end{definition}

\noindent
F\"ur $2 \le m \in \mathbb{N}$ bezeichne $\mathcal{D}_m \subset PSL_2(\mathbb{C})$ stets eine (projektive) $m$-Diedergruppe,
und $\mathcal{T} \subset PSL_2(\mathbb{C})$ bezeichne stets eine (projektive) Tetraedergruppe.
Dann ist $P^{-1}(\mathcal{D}_m)$ bzw. $P^{-1}(\mathcal{T})$ eine bin\"are $m$-Diedergruppe bzw. bin\"are Tetraedergruppe.
$\mathcal{D}_3$, $\mathcal{T}$ bzw. $\mathcal{D}_2$ sind isomorph zur
symmetrischen Gruppe $\mathcal{S}_3$, alternierenden Gruppe $\mathcal{A}_4$ bzw. Kleinschen Vierergruppe $\mathcal{V}_4$.
Eine Tetraedergruppe $\mathcal{T}$ enth\"alt genau eine $2$-Diedergruppe $\mathcal{D}_2$,
und umgekehrt ist eine $2$-Diedergruppe $\mathcal{D}_2$ Normalteiler in genau einer Tetraedergruppe $\mathcal{T}$.
\\
Wir fassen die bekannten Grundlagen in den Lemmata \ref{Item06-02} und \ref{Item06-03} zusammen:

\begin{lemma} \label{Item06-02}

Sei $d \in \mathbb{N}$ quadratfrei, und sei $k = \mathbb{Q}(i\sqrt{d})$.
\\
Sei $M$ eine $ \mathbb{Q}$-Quaternionenalgebra oder eine $k$-Quaternionenalgebra.
\\
Ist $\mathcal{G} \subset P\Gamma(M)$ eine nichttriviale endliche Gruppe,
dann hat $\mathcal{G}$ die Ordnung $2$ oder $3$,
oder $\mathcal{G}$ ist eine $3$-Diedergruppe, Tetraedergruppe oder $2$-Diedergruppe.

\end{lemma}

\begin{proof}

Sei $\mathcal{G}$ eine nichttriviale endliche Untergruppe von $P\Gamma(M)$,
und sei $\mathcal{C}_m \subset \mathcal{G}$ eine zyklische Untergruppe der Ordnung $m > 1$.
Dann ist $P^{-1}(\mathcal{C}_m)$ zyklisch und hat die Ordnung $2m$.
Sei $U$ ein erzeugendes Element von $P^{-1}(\mathcal{C}_m)$,
und sei $\zeta \in \mathbb{C}$ eine primitive $2m$-te Einheitswurzel.
Dann ist $S(U) = \zeta + \overline{\zeta} \in \mathfrak{o} \cap \mathbb{R} = \mathbb{Z}$.
Daher muss $m = 2$ oder $m = 3$ sein.
Nach \cite[Chapter 4.4]{Sp} ist $P^{-1}(\mathcal{G})$ also
eine zyklische Gruppe der Ordnung $4$ oder $6$, oder eine bin\"are $3$-Diedergruppe, bin\"are Tetraedergruppe oder bin\"are $2$-Diedergruppe.

\end{proof}

\begin{lemma} \label{Item06-03}

Sei $\mathcal{G}$ eine $3$-Diedergruppe $\mathcal{D}_3$, Tetraedergruppe $\mathcal{T}$ oder $2$-Diedergruppe $\mathcal{D}_2$.
\\
Dann ist $F(\mathcal{G})$ eine $\mathbb{Q}$-Quaternionenalgebra,
$\mathfrak{F}(\mathcal{G})$ ist eine $F(\mathcal{G})$-Ordnung,
und es gilt $P^{-1}(\mathcal{G}) = \Gamma(\mathfrak{F}(\mathcal{G})) = \mathfrak{F}(\mathcal{G})^\times$.
Im Einzelnen gilt:
\begin{enumerate}[(i)]

\item \label{Item06-03.1}
$P^{-1}(\mathcal{D}_3)$ wird von zwei Elementen $U, V$ erzeugt, mit $U^3 = -1$ und $VUV^{-1} = U^{-1}$.
\\
$\Sigma(F(\mathcal{D}_3)) = -3$ und $\Lambda(\mathfrak{F}(\mathcal{D}_3)) = 1$, also $\Delta(\mathfrak{F}(\mathcal{D}_3)) = -3$.

\item \label{Item06-03.2}
$P^{-1}(\mathcal{D}_2)$ wird von zwei Elementen $U, V$ erzeugt, mit $U^2 = -1$ und $VUV^{-1} = U^{-1}$.
\\
Falls $\mathcal{D}_2 \subset \mathcal{T}$, wird $P^{-1}(\mathcal{T})$ von $U, V$ und $W := (1 - U - V - UV)/2$ erzeugt.
\\
Dann gilt $F(\mathcal{D}_2) = F(\mathcal{T})$, $\Sigma(F(\mathcal{D}_2)) = \Sigma(F(\mathcal{T})) = -2$
und $\mathfrak{F}(\mathcal{D}_2) \subset \mathfrak{F}(\mathcal{T})$,
\\
$\Lambda(\mathfrak{F}(\mathcal{T})) = 1$, also  $\Delta(\mathfrak{F}(\mathcal{T})) = -2$,
und $\Lambda(\mathfrak{F}(\mathcal{D}_2)) = 2$, also $\Delta(\mathfrak{F}(\mathcal{D}_2)) = -4$.

\end{enumerate}

\end{lemma}

\begin{proof} ${}$
\begin{enumerate}[(i)]

\item [(\ref{Item06-03.1})]
$P^{-1}(\mathcal{D}_3)$ wird von zwei Elementen $U$, $V$ erzeugt,
mit $U^3 = -1$ ($U^2 - U + 1 = 0$), $V^2 = -1$, $VUV^{-1} = U^{-1}$, siehe \cite[4.4.7]{Sp}.
 $U$, $V$ erzeugen \"uber $\mathbb{Q}$ eine Quaternionenalgebra $F(\mathcal{D}_3)$
und \"uber $\mathbb{Z}$ eine $F(\mathcal{D}_3)$-Ordnung $\mathfrak{F}(\mathcal{D}_3)$ mit Diskriminante $9\mathbb{Z}$.
\\
Also zerf\"allt $F(\mathcal{D}_3)$ an allen Stellen $p \not= 3, \infty$ von $\mathbb{Q}$.
Jedes $A \in F(\mathcal{D}_3)$ ist darstellbar als $A = \alpha + \beta U + \gamma V + \delta UV$ mit $\alpha, \beta, \gamma, \delta \in \mathbb{Q}$.
Falls $A \not= 0$ ist, ist $N(A) = \alpha ^2 + \alpha \beta + \beta ^2 + \gamma ^2 + \gamma \delta + \delta ^2 > 0$.
Daher ist $F(\mathcal{D}_3)$ Divisionsalgebra.
\\
Die Anzahl der Verzweigungsstellen von $F(\mathcal{D}_3)$ ist gerade.
Also ist $F(\mathcal{D}_3)$ genau an den Stellen $3$ und $\infty$ verzweigt,
und $\mathfrak{F}(\mathcal{D}_3)$ ist eine $F(\mathcal{D}_3)$-Maximalordnung.

\item [(\ref{Item06-03.2})]
Sei $\mathcal{D}_2 \subset \mathcal{T}$.
Dann wird $P^{-1}(\mathcal{D}_2)$ von zwei Elementen $U$, $V$ erzeugt, mit $U^2 = -1$, $V^2 = -1$, $VUV^{-1} = U^{-1}$, siehe \cite[4.4.7]{Sp}.
$P^{-1}(\mathcal{T})$ wird von $U$ und $V$ sowie von einem Element $W$ erzeugt, mit $W^3 = -1$, $UWU^{-1} = VW$ und $WV = UW$.
Eine elementare Rechnung zeigt, dass $W = (1 - U - V - UV)/2$ gilt, siehe \cite[4.4.10]{Sp}.
\\
$U$, $V$, $W$ erzeugen \"uber $\mathbb{Q}$ eine Quaternionenalgebra $F(\mathcal{T}) = F(\mathcal{D}_2)$
und \"uber $\mathbb{Z}$ eine Ordnung $\mathfrak{F}(\mathcal{T})$ mit Diskriminante $4\mathbb{Z}$.
Also zerf\"allt $F(\mathcal{T})$ an allen Stellen $p \not= 2, \infty$.
\\
Jedes $A \in F(\mathcal{T})$ ist darstellbar als $A = \alpha + \beta U + \gamma V + \delta UV$ mit $\alpha, \beta, \gamma, \delta \in \mathbb{Q}$.
Falls $A \not= 0$, ist $N(A) = \alpha ^2 + \beta ^2 + \gamma ^2 + \delta ^2 > 0$.
Daher ist $F(\mathcal{T})$ Divisionsalgebra.
Die Anzahl der Verzweigungsstellen von $F(\mathcal{T})$ ist gerade.
Also ist $F(\mathcal{T})$ genau an den Stellen $2$ und $\infty$ verzweigt,
und $\mathfrak{F}(\mathcal{T})$ ist eine $F(\mathcal{T})$-Maximalordnung.
Man rechnet nun leicht nach, dass $\mathfrak{F}(\mathcal{D}_2)$ die Diskriminante $16\mathbb{Z}$ hat,
also $[\mathfrak{F}(\mathcal{T}) : \mathfrak{F}(\mathcal{D}_2)] = 2$ gilt.

\end{enumerate}

\noindent
Bei den Erzeugendenrelationen f\"ur $P^{-1}(\mathcal{D}_3)$ und $P^{-1}(\mathcal{D}_2)$ k\"onnen wir die Gleichung $V^2 = -1$ weglassen,
denn aus $V^2UV^{-2} = VU^{-1}V^{-1} = (VUV^{-1})^{-1} = U$ erhalten wir $V^2 \in \mathbb{C}[U] \cap \mathbb{C}[V] = \mathbb{C}$,
wegen $V \not\in \mathbb{C}$ also $S(V) = 0$ und $V^2 = -N(V) = -1$.

\noindent
Wegen $N (F(\mathcal{G})^\times) \subset \mathbb{Q}^+$
ist $\mathfrak{F}(\mathcal{G})^\times = \Gamma(\mathfrak{F}(\mathcal{G}))$ endlich.
Wegen $\mathcal{G} \subset P\Gamma(\mathfrak{F}(\mathcal{G}))$ folgt mit Lemma \ref{Item06-02},
dass $P\Gamma(\mathfrak{F}(\mathcal{G})) \in \{ \mathcal{D}_3, \mathcal{T}, \mathcal{D}_2 \}$.
Falls $\mathcal{G} \in \{ \mathcal{D}_3, \mathcal{T} \}$, ist also $\mathcal{G} = P\Gamma(\mathfrak{F}(\mathcal{G}))$.
Falls $\mathcal{G} = \mathcal{D}_2$, gilt $W \not\in \mathfrak{F}(\mathcal{D}_2)$,
also $P\Gamma(\mathfrak{F}(\mathcal{D}_2)) \subsetneqq P\Gamma(\mathfrak{F}(\mathcal{T})) = \mathcal{T}$
und $P\Gamma(\mathfrak{F}(\mathcal{D}_2)) = \mathcal{D}_2$.

\end{proof}

\begin{theorem} \label{Item06-04}

Sei $d \in \mathbb{N}$ quadratfrei, und sei $k = \mathbb{Q}(i\sqrt{d})$.
\\
Seien $\mathcal{D}_3$ eine $3$-Diedergruppe und $\mathcal{T}$ eine Tetraedergruppe. Dann gilt:
\begin{enumerate}[(i)]

\item \label{Item06-04.1}
Es gibt bis auf Isomorphie genau eine $k$-Quaternionenalgebra $M$ mit $\mathcal{D}_3 \subset P\Gamma(M)$
und bis auf Isomorphie genau eine $M$-Maximalordnung $\mathfrak{M}$ mit $\mathcal{D}_3 \subset P\Gamma(\mathfrak{M})$.
\\
Genau dann gilt $M \cong M_2(k)$, wenn $d \not\equiv 2 \pmod{3}$.

\item \label{Item06-04.2}
Es gibt bis auf Isomorphie genau eine $k$-Quaternionenalgebra $M$ mit $\mathcal{T} \subset P\Gamma(M)$
und bis auf Isomorphie genau eine $M$-Maximalordnung $\mathfrak{M}$ mit $\mathcal{T} \subset P\Gamma(\mathfrak{M})$.
\\
Genau dann gilt $M \cong M_2(k)$, wenn $d \not\equiv 7 \pmod{8}$.

\end{enumerate}

\end{theorem}

\begin{proof}

Sei $M$ eine $k$-Quaternionenalgebra, und sei $\mathcal{G} = \mathcal{D}_3$ oder $\mathcal{G} = \mathcal{T}$.
\\
$\mathcal{G} \subset P\Gamma(M)$ gilt genau dann, wenn $M$ eine Erweiterung von $F(\mathcal{G})$ zur $k$-Algebra ist,
was den Isomorphietyp von $M$ festlegt. Der Isomorphietyp von $\mathfrak{M}$ ist eindeutig nach Satz \ref{Item05-06}.
\\
Nach Lemma \ref{Item03-02}.(\ref{Item03-02.2}) ist $M \cong M_2(k)$ genau dann,
wenn $\Sigma_k(F(\mathcal{G})) = \Sigma_k(M_2(\mathbb{Q})) = 1$ gilt,
also wenn die einzige endliche Verzweigungsstelle $p$ von $F(\mathcal{G})$ in $k$ nicht zerlegt ist,
also $3$ bzw. $2$ in $k$ nicht zerlegt ist, also wenn $-d \not\equiv 1 \pmod{3}$ bzw. $-d \not\equiv 1 \pmod{8}$ ist.

\end{proof}

\noindent
$\mathcal{D}_3 \subset P\Gamma(\mathfrak{M})$ bzw. $\mathcal{T} \subset P\Gamma(\mathfrak{M})$
ist stets eine maximalendliche Untergruppe von $P\Gamma(\mathfrak{M})$.
Ist aber $\mathcal{D}_2 \subset P\Gamma(\mathfrak{M})$,
so  gibt es (genau) eine Tetraedergruppe $\mathcal{T}$ mit $\mathcal{D}_2 \subset \mathcal{T} \subset P\Gamma(M)$,
und $\mathcal{D}_2$ ist genau dann maximalendliche Untergruppe von $P\Gamma(\mathfrak{M})$, wenn $\mathcal{T} \not\subset P\Gamma(\mathfrak{M})$.

\begin{lemma} \label{Item06-05}

Sei $d \in \mathbb{N}$ quadratfrei, und sei $k = \mathbb{Q}(i\sqrt{d})$ mit Hauptordnung $\mathfrak{o}$.
\\
Sei $\mathcal{D}_2$ eine $2$-Diedergruppe, und sei $M$ eine $k$-Quaternionenalgebra mit $\mathcal{D}_2 \subset P\Gamma(M)$.
\begin{enumerate}[(i)]

\item \label{Item06-05.1}
Ist $\mathfrak{M}$ eine $M$-Maximalordnung, so gilt:
Genau dann ist $\mathcal{D}_2$ eine maximalendliche Untergruppe von $P\Gamma(\mathfrak{M})$,
wenn $F(\mathcal{D}_2) \cap \mathfrak{M} = \mathfrak{F}(\mathcal{D}_2)$.

\item \label{Item06-05.2}
Genau dann ist $\mathfrak{F}(\mathcal{D}_2)$ eine $\mathfrak{o}$-kompatible $M$-Ordnung,
wenn $d \not\equiv 3 \pmod{4}$.

\end{enumerate}

\end{lemma}

\begin{proof}

Wir k\"urzen $F :=  F(\mathcal{D}_2) =  F(\mathcal{T})$ ab.
\begin{enumerate}[(i)]

\item [(\ref{Item06-05.1})]
Genau dann gilt $\mathcal{D}_2 \subset P\Gamma(\mathfrak{M})$ und $\mathcal{T} \not\subset P\Gamma(\mathfrak{M})$,
wenn $\mathfrak{F}(\mathcal{D}_2) \subset \mathfrak{M}$ und $\mathfrak{F}(\mathcal{T}) \not\subset \mathfrak{M}$.
\\
Wir zeigen die \"Aquivalenz zu $F \cap \mathfrak{M} = \mathfrak{F}(\mathcal{D}_2)$ durch Widerspruch.
Wir nehmen also $\mathfrak{F}(\mathcal{D}_2) \subsetneqq F \cap \mathfrak{M}$ an.
Wegen $\Lambda(\mathfrak{F}(\mathcal{D}_2)) = 2$ ist dann $F \cap \mathfrak{M}$ eine $F$-Maximalordnung.
Da $F$ an der Stelle $2$ verzweigt ist, ist $\mathfrak{F}(\mathcal{T})_2$ die $F_2$-Maximalordnung
und $\mathfrak{F}(\mathcal{T})$ die einzige $F$-Maximalordnung, die $\mathfrak{F}(\mathcal{D}_2)$ enth\"alt.
Also ist $F \cap \mathfrak{M} = \mathfrak{F}(\mathcal{T})$, Widerspruch.

\item [(\ref{Item06-05.2})]
Nach Satz \ref{Item04-05}.(\ref{Item04-05.1}) gibt es genau eine $F_2$-Ordnung $\mathfrak{G}_2 \subset \mathfrak{F}(\mathcal{T})_2$ mit $[\mathfrak{F}(\mathcal{T})_2 : \mathfrak{G}_2] = 2$,
also genau eine $F$-Ordnung $\mathfrak{G} \subset \mathfrak{F}(\mathcal{T})$ mit $\Lambda(\mathfrak{G}) = 2$.
Nach Satz \ref{Item05-04} ist $\mathfrak{F}(\mathcal{D}_2)$ also genau dann eine $\mathfrak{o}$-kompatible $F$-Ordnung,
wenn $2 \nmid \Sigma_k(F)$ und $2 \in N_{abs}(I)$ gilt,
also wenn $2$ in $k$ nicht zerlegt und nicht tr\"age ist,
also wenn $-d \not\equiv 1 \pmod{4}$ gilt.

\end{enumerate}

\end{proof}

\begin{theorem} \label{Item06-06}

Sei $d \in \mathbb{N}$ quadratfrei, und sei $k = \mathbb{Q}(i\sqrt{d})$.
\\
Genau dann gibt es eine $k$-Quaternionenalgebra $M$ und eine $M$-Maximalordnung $\mathfrak{M}$
mit einer maximalendlichen $2$-Diedergruppe $\mathcal{D}_2 \subset P\Gamma(\mathfrak{M})$, wenn $d \not\equiv 3 \pmod{4}$.
\\
In diesem Fall gilt $M \cong M_2(k)$, und $\mathfrak{M}$ ist bis auf Isomorphie eindeutig.

\end{theorem}

\begin{proof}

Die erste Behauptung folgt mit Lemma \ref{Item06-05} und Satz \ref{Item05-06}.(\ref{Item05-06.1})(\ref{Item06-05.2}).
\\
Die zweite Behauptung folgt wegen $d \not\equiv 7 \pmod{8}$
mit Satz \ref{Item06-04}.(\ref{Item06-04.2}) und Satz \ref{Item05-06}.(\ref{Item05-06.4}).

\end{proof}

\begin{theorem} \label{Item06-07}

Sei $d \in \mathbb{N}$ quadratfrei, und sei $k = \mathbb{Q}(i\sqrt{d})$.
\\
Sei $\mathfrak{M}$ eine $M_2(k)$-Maximalordnung. Dann gilt:
\begin{enumerate}[(i)]

\item \label{Item06-07.1}
$P\Gamma(\mathfrak{M})$ enth\"alt genau dann eine $3$-Diedergruppe $\mathcal{D}_3$,
\\
wenn f\"ur alle Stellen $p \ne 3, \infty$ von $\mathbb{Q}$ gilt:
$\left(\dfrac{-3 \Lambda(\mathfrak{M} \cap M_2(\mathfrak{o})),-d}{p}\right) = 1$.

\item \label{Item06-07.2}
$P\Gamma(\mathfrak{M})$ enth\"alt genau dann eine Tetraedergruppe $\mathcal{T}$,
\\
wenn f\"ur alle Stellen $p \ne 2, \infty$ von $\mathbb{Q}$ gilt:
$\left(\dfrac{-2 \Lambda(\mathfrak{M} \cap M_2(\mathfrak{o})),-d}{p}\right) = 1$.

\item \label{Item06-07.3}
$P\Gamma(\mathfrak{M})$ enth\"alt genau dann eine maximalendliche $2$-Diedergruppe $\mathcal{D}_2$,
\\
wenn f\"ur alle Stellen $p \ne 2, \infty$ von $\mathbb{Q}$ gilt:
$\left(\dfrac{-\Lambda(\mathfrak{M} \cap M_2(\mathfrak{o})),-d}{p}\right) = 1$.

\end{enumerate}

\end{theorem}

\begin{proof}

Wir k\"urzen $\lambda := \Lambda(\mathfrak{M} \cap M_2(\mathfrak{o}))$ ab.
\begin{enumerate}[(i)]

\item [(\ref{Item06-07.1})]
Sei $\mathcal{D}_3 \subset P\Gamma(\mathfrak{M})$,
oder \"aquivalent $\mathfrak{F}(\mathcal{D}_3) \subset \mathfrak{M}$,
oder \"aquivalent $\mathfrak{F}(\mathcal{D}_3) = F(\mathcal{D}_3) \cap \mathfrak{M}$.
\\
Nach Korollar \ref{Item05-10} und Lemma \ref{Item06-03}.(\ref{Item06-03.1}) ist dann
$\left(\dfrac{-3 \lambda,-d}{p}\right) = 1$ f\"ur alle $p \ne 3, \infty$.

Sei nun umgekehrt $\left(\dfrac{-3 \lambda,-d}{p}\right) = 1$ f\"ur alle $p \ne 3, \infty$.
Die Anzahl der Stellen $p$ mit $\left(\dfrac{-3 \lambda,-d}{p}\right) = -1$ ist gerade.
Also ist $\left(\dfrac{-3 \lambda,-d}{3}\right) = -1$,
speziell $d \not\equiv 2 \pmod{3}$.
\\
Sei $\mathcal{D}'_3$ eine $3$-Diedergruppe.
Nach Satz \ref{Item06-04}.(\ref{Item06-04.1}) gibt es eine $M_2(k)$-Maximalordnung $\mathfrak{M}'$
mit $\mathcal{D}'_3 \subset P\Gamma(\mathfrak{M}')$.
Dann ist $\left(\dfrac{-3 \Lambda(\mathfrak{M}' \cap M_2(\mathfrak{o})),-d}{p}\right) = \left(\dfrac{-3 \lambda,-d}{p}\right)$
f\"ur alle Stellen $p$.
Nach Lemma \ref{Item05-05}.(\ref{Item05-05.2}) ist
$\Lambda(\mathfrak{M} \cap \mathfrak{M}') \lambda \Lambda(M_2(\mathfrak{o}) \cap \mathfrak{M}') \in \mathbb{Q}^{\times(2)}$.
Nach Lemma \ref{Item05-05}.(\ref{Item05-05.3}) ist also $\mathfrak{M}$ isomorph zu $\mathfrak{M}'$.
Daher gibt es eine $3$-Diedergruppe $\mathcal{D}_3 \subset P\Gamma(\mathfrak{M})$.

\item [(\ref{Item06-07.2})]
folgt mit Lemma \ref{Item06-03}.(\ref{Item06-03.2}) und umgekehrt mit Satz \ref{Item06-04}.(\ref{Item06-04.2}) wie im Beweis von (\ref{Item06-07.1}).

\item [(\ref{Item06-07.3})]
Sei $\mathcal{D}_2 \subset P\Gamma(\mathfrak{M})$ und $\mathcal{T} \not\subset P\Gamma(\mathfrak{M})$.
Nach Lemma \ref{Item06-05}.(\ref{Item06-05.1}), Korollar \ref{Item05-10} und Lemma \ref{Item06-03}.(\ref{Item06-03.2})
ist dann $\left(\dfrac{-\lambda,-d}{p}\right) = 1$ f\"ur alle $p \ne 2, \infty$.

Sei nun umgekehrt $\left(\dfrac{-\lambda,-d}{p}\right) = 1$ f\"ur alle $p \ne 2, \infty$.
Die Anzahl der Stellen $p$ mit $\left(\dfrac{-\lambda,-d}{p}\right) = -1$ ist gerade.
Also ist $\left(\dfrac{-\lambda,-d}{2}\right) = -1$ und speziell $d \not\equiv 7 \pmod{8}$.
\\
Falls $d \equiv 3 \pmod{8}$, w\"are $2$ in $k$ tr\"age.
Nach Satz \ref{Item05-06}.(\ref{Item05-06.1}) und \ref{Item05-04}.(\ref{Item05-04.1}) gilt $\lambda \in N_{abs}(I)$,
der quadratfreie Anteil von $\lambda$ w\"are ungerade, Widerspruch.
Also ist $d \not\equiv 3 \pmod{4}$.
\\
Nach Satz \ref{Item06-06} gibt es eine $M_2(k)$-Maximalordnung $\mathfrak{M}'$
mit $\mathcal{D}'_2 \subset P\Gamma(\mathfrak{M}')$ und $\mathcal{T}' \not\subset P\Gamma(\mathfrak{M}')$.
Dann gilt $\left(\dfrac{-\Lambda(\mathfrak{M}' \cap M_2(\mathfrak{o})),-d}{p}\right) = \left(\dfrac{-\lambda,-d}{p}\right)$
f\"ur alle Stellen $p$,
und nach den Lemmata \ref{Item05-05}.(\ref{Item05-05.2}), \ref{Item05-05}.(\ref{Item05-05.3}) ist $\mathfrak{M}'$ isomorph zu $\mathfrak{M}$.

\end{enumerate}

\end{proof}

\begin{theorem} \label{Item06-08}

Sei $d \in \mathbb{N}$ quadratfrei, und sei $k = \mathbb{Q}(i\sqrt{d})$ mit Hauptordnung $\mathfrak{o}$.
Dann gilt:
\begin{enumerate}[(i)]

\item \label{Item06-08.1}
$PSL_2(\mathfrak{o})$ enth\"alt genau dann eine $3$-Diedergruppe $\mathcal{D}_3$,
\\
wenn $p \equiv 1 \pmod{3}$ f\"ur alle Primteiler $p \not= 3$ von $d$.

\item \label{Item06-08.2}
$PSL_2(\mathfrak{o})$ enth\"alt genau dann eine Tetraedergruppe $\mathcal{T}$,
\\
wenn $p \equiv 1$ oder $p \equiv 3 \pmod{8}$ f\"ur alle Primteiler $p \not= 2$ von $d$.

\item \label{Item06-08.3}
$PSL_2(\mathfrak{o})$ enth\"alt genau dann eine maximalendliche $2$-Diedergruppe $\mathcal{D}_2$,
\\
wenn $p \equiv 1 \pmod{4}$ f\"ur alle Primteiler $p \not= 2$ von $d$.

\end{enumerate}

\end{theorem}

\begin{proof}

Nach Satz \ref{Item06-07} gilt:
\begin{enumerate}[(i)]

\item[(\ref{Item06-08.1})]
$PSL_2(\mathfrak{o})$ enth\"alt genau dann eine $3$-Diedergruppe $\mathcal{D}_3$,
\\
wenn f\"ur alle Stellen $p \ne 3, \infty$ von $\mathbb{Q}$ gilt: $\left(\dfrac{-3,-d}{p}\right) = 1$.

\item[(\ref{Item06-08.2})]
$PSL_2(\mathfrak{o})$ enth\"alt genau dann eine Tetraedergruppe $\mathcal{T}$,
\\
wenn f\"ur alle Stellen $p \ne 2, \infty$ von $\mathbb{Q}$ gilt: $\left(\dfrac{-2,-d}{p}\right) = 1$.

\item[(\ref{Item06-08.3})]
$PSL_2(\mathfrak{o})$ enth\"alt genau dann eine maximalendliche $2$-Diedergruppe $\mathcal{D}_2$,
\\
wenn f\"ur alle $p \ne 2, \infty$ gilt: $\left(\dfrac{-1,-d}{p}\right) = 1$.

\end{enumerate}

\noindent
Die Behauptung folgt dann jeweils durch Auswertung der Hilbertsymbole.

\end{proof}


\section{Konjugationsklassenanzahlen} \label{Chapter07}

\begin{definition} \label{Item07-01}
Sei $d \in \mathbb{N}$ quadratfrei, und sei $k = \mathbb{Q}(i\sqrt{d})$.
Sei $\mathcal{G}$ eine endliche Gruppe.
\\
Seien $F$ eine $\mathbb{Q}$-Quaternionenalgebra und $M$ eine Erweiterung von $F$ zur $k$-Quater\-nionenalgebra.
Seien $\mathfrak{G}$ eine $F$-Ordnung und $\mathfrak{N}$ eine $M$-Ordnung.

\begin{enumerate}[(i)]

\item \label{Item07-01.1}
Bezeichne $\mu(\mathcal{G}, \mathfrak{N})$ die Anzahl der $P\Gamma(\mathfrak{N})$-Konjugationsklassen
maximalendlicher Untergruppen $\mathcal{G}' \subset P\Gamma(\mathfrak{N})$ vom Isomorphietyp $\mathcal{G}$.

\item \label{Item07-01.2}
Bezeichne $B(\mathfrak{G}, \mathfrak{N})$ die Anzahl der $\mathfrak{N}^\times$-Konjugationsklassen
optimaler Einbettungen $j: \mathfrak{G} \hookrightarrow \mathfrak{N}$,
und $B_1(\mathfrak{G}, \mathfrak{N})$ die entsprechende Anzahl der $\Gamma(\mathfrak{N})$-Konjugationsklassen.

\item \label{Item07-01.3}
F\"ur $\mathfrak{O} = \mathfrak{G}$ oder $\mathfrak{O} = \mathfrak{N}$
bezeichne $\mathcal{A}(\mathfrak{O})$ die Gruppe der Automorphismen von $\mathfrak{O}$
und $\mathcal{I}(\mathfrak{O}) \subset \mathcal{A}(\mathfrak{O})$ die Unterguppe der Konjugationen mit Elementen aus $\mathfrak{O}^\times$.

\end{enumerate}

\end{definition}

\noindent
$Bemerkung$: Die Bezeichnungen $\mu(\mathcal{D}_3, \mathfrak{N})$, $\mu(\mathcal{T}, \mathfrak{N})$ und $\mu(\mathcal{D}_2, \mathfrak{N})$
entsprechen den in \cite[§ 22 und § 24]{Kr} definierten Bezeichnungen
$\mu_3(\mathfrak{N})$, $\mu_{\mathcal{T}}(\mathfrak{N}) = \mu^{\mathcal{T}}_2(\mathfrak{N})$ und $\mu^-_2(\mathfrak{N})$.

\begin{lemma} \label{Item07-02}

Sei $d \in \mathbb{N}$ quadratfrei, sei $k = \mathbb{Q}(i\sqrt{d})$ mit Diskriminante $D$,
und bezeichne $t$ die Anzahl der verschiedenen Primteiler von $D$.
Sei $F$ eine $\mathbb{Q}$-Quaternionenalgebra,
und bezeichne $r$ die Anzahl der Primteiler von $\Sigma_k(F)$.
Sei $s$ der Rang der $t$$\times$$r$-Matrix $(h_{pq})$
mit Koeffizienten in $\mathbb{Z}/2\mathbb{Z}$
und $h_{pq} := 0$ oder $1$, je nachdem, ob $\left(\dfrac{q,-d}{p}\right) = 1$ oder $-1$ gilt,
wo $p$ die verschiedenen Primteiler von $D$ und $q$ die Primteiler von $\Sigma_k(F)$ durchl\"auft.
Seien $M$ eine Erweiterung von $F$ zur $k$-Quaternionenalgebra
und $\mathfrak{M}$ eine $M$-Maximalordnung.
Dann gilt $[\mathcal{A}(\mathfrak{M}) : \mathcal{I}(\mathfrak{M})] = 2^{t+2r-s-1}$.

\end{lemma}

\begin{proof}

Sei $R$ sei die von den Verzweigungsstellen von $M$ erzeugte $\mathfrak{o}$-Idealgruppe.
Ist $j \in \mathcal{A}(\mathfrak{M})$, so
gibt es ein $J \in M^\times$ mit $j(A) = JAJ^{-1}$ f\"ur alle $A \in \mathfrak{M}$,
und $\mathfrak{J} := \mathfrak{M} J$ ist ein zweiseitiges $\mathfrak{M}$-Ideal
mit $\mathfrak{j} := N(\mathfrak{J}) \in RI^{(2)} \cap H$.
Da $J$ bis auf einen Faktor $a \in k^\times$ eindeutig ist,
induziert die Zuordnung $j \mapsto \mathfrak{j}$ wegen $N(a) = a^2$ einen
surjektiven Homomorphismus $\mathcal{A}(\mathfrak{M}) \rightarrow (RI^{(2)} \cap H) / H^{(2)}$
mit Kern $\mathcal{I}(\mathfrak{M})$,
siehe \cite[Lemma 6.7.5, Theorem 7.7.7]{Ma}.
Also gilt $[\mathcal{A}(\mathfrak{M}) : \mathcal{I}(\mathfrak{M})] = [(RI^{(2)} \cap H) : (I^{(2)} \cap H)] [(I^{(2)} \cap H) : H^{(2)}]$.
\\[4pt]
Die Zuordnungsvorschrift $\mathfrak{a} \mapsto \mathfrak{a}^2$ 
induziert einen Isomorphismus $I/H \rightarrow I^{(2)}/H^{(2)}$.
Also gilt $[I : H] = [I^{(2)} : H^{(2)}]$
und $[I : I^{(2)}H] [ I^{(2)}H : H] = [I^{(2)} : (I^{(2)} \cap H)] [ (I^{(2)} \cap H) : H^{(2)}]$.
Wegen $[I^{(2)} : (I^{(2)} \cap H)] = [I^{(2)}H : H]$ gilt also $[(I^{(2)} \cap H) : H^{(2)}] = [I : I^{(2)}H] = 2^{t-1}$.
\\
Weiter gilt $[RI^{(2)} : (RI^{(2)} \cap H)] = [RI^{(2)}H : H] = [RI^{(2)}H : I^{(2)}H] [(I^{(2)} : (I^{(2)} \cap H)]$ und
$[RI^{(2)} : I^{(2)}] [(I^{(2)} : (I^{(2)} \cap H)] = [RI^{(2)} :  (RI^{(2)} \cap H)] [(RI^{(2)} \cap H) : (I^{(2)} \cap H)]$, also 
\\
$[(RI^{(2)} \cap H) : (I^{(2)} \cap H)] = [RI^{(2)} : I^{(2)}]  [RI^{(2)}H : I^{(2)}H]^{-1} = 2^{2r} [R: (R \cap  I^{(2)}H)]^{-1}$.
\\[4pt]
Wir zeigen noch $[R: (R \cap  I^{(2)}H)] = 2^s$:
Mit der Multiplikation als Verkn\"upfung ist $R / (R \cap  I^{(2)}H)$ ein $\mathbb{Z}/2\mathbb{Z}$-Modul.
F\"ur $f \in \mathbb{N}$ sei $(f_p)$ das $t$-Tupel mit Koeffizienten aus $\mathbb{Z}/2\mathbb{Z}$
und $f_p := 0$ oder $1$, je nachdem, ob $\left(\dfrac{f,-d}{p}\right) = 1$ oder $-1$ gilt,
wo $p$ die verschiedenen Primteiler von $D$ durchl\"auft.
Die Zuordnungsvorschrift $\mathfrak{f} \mapsto (N_{abs}(\mathfrak{f})_p)$ induziert einen
$\mathbb{Z}/2\mathbb{Z}$-Modulhomomphismus $\Psi: R / (R \cap  I^{(2)}H) \rightarrow (\mathbb{Z}/2\mathbb{Z})^t$.
Das Bild von $\Psi$ wird von den Tupeln $(q_p)$ erzeugt,
wo $q$ die Primteiler von $\Sigma_k(F)$ durchl\"auft, hat also die \mbox{Dimension $s$}.
Nach \cite[Kapitel III, § 8, S\"atze 6, 7 und Kapitel I, § 7, Satz 1]{Bo} gilt $\mathfrak{f} \in I^{(2)}H$ genau dann,
wenn $\left(\dfrac{N_{abs}(\mathfrak{f}),-d}{p}\right) = 1$ an allen Stellen $p$ von $\mathbb{Q}$ (mit $p \mid D$).
Also ist $\Psi$ Isomorphismus.

\end{proof}

\begin{theorem} \label{Item07-03}

Sei $d \in \mathbb{N}$ quadratfrei, und sei $k = \mathbb{Q}(i\sqrt{d})$.
\\
Seien $F$ eine $\mathbb{Q}$-Quaternionenalgebra und $M$ eine Erweiterung von $F$ zur $k$-Quater\-nionenalgebra.
Seien $\mathfrak{G}$ eine $F$-Ordnung und $\mathfrak{N}$ eine $M$-Maximalordnung mit $\mathfrak{G} = F \cap \mathfrak{N}$.
\\
Sei $C(\mathfrak{G})$ gleich dem Produkt der in Satz \ref{Item04-09} angegebenen
$C_p(\mathfrak{G}_p)$ \"uber alle endlichen Stellen $p$ von $\mathbb{Q}$.
Dann ist $B(\mathfrak{G}, \mathfrak{N}) = C(\mathfrak{G}) [\mathcal{A}(\mathfrak{N}) : \mathcal{I}(\mathfrak{N})]$
und $B_1(\mathfrak{G}, \mathfrak{N}) = 2 B(\mathfrak{G}, \mathfrak{N})$.

\end{theorem}

\begin{proof}

$C(\mathfrak{G})$ ist gleich der Anzahl der $M$-Maximalordnungen $\mathfrak{L}$ mit $\mathfrak{G} = F \cap \mathfrak{L}$.
Ist $\mathfrak{L}$ eine $M$-Maximalordnung mit $\mathfrak{G} = F \cap \mathfrak{L}$,
so gibt es nach Satz \ref{Item05-06}.(\ref{Item05-06.4})
einen Automorphismus $j$ von $M$ mit $ j(\mathfrak{L}) = \mathfrak{N}$ und $j(\mathfrak{G}) = j(F) \cap \mathfrak{N}$.
F\"ur jedes $\mathfrak{L}$ mit $\mathfrak{G} = F \cap \mathfrak{L}$
w\"ahlen wir einen Automorphismus $j_\mathfrak{L}$ von $M$ mit $\mathfrak{N} = j_\mathfrak{L}(\mathfrak{L})$.
Eine Einbettung $j:  \mathfrak{G} \hookrightarrow \mathfrak{N}$ mit $j(\mathfrak{G}) = j(F) \cap \mathfrak{N}$
ist zu einem Automorphismus von $M$ fortsetzbar.
Mit $\mathfrak{L} := j^{-1}(\mathfrak{N})$ gilt $\mathfrak{G} = F \cap \mathfrak{L}$
und $a(j) := j \circ j_\mathfrak{L}^{-1} \in \mathcal{A}(\mathfrak{N})$.
Die Abbildung mit Zuordnungsvorschrift $j \mapsto (\mathfrak{L}, a(j))$ ist injektiv.
Sind umgekehrt eine $M$-Maximalordnung $\mathfrak{L}$ mit $\mathfrak{G} = F \cap \mathfrak{L}$
und ein $a \in \mathcal{A}(\mathfrak{N})$ gegeben,
dann ist $j: = a \circ j_\mathfrak{L}$ ein Automorphismus von $M$,
und seine Einschr\"ankung auf $\mathfrak{G}$ induziert eine Einbettung
$j: \mathfrak{G} \hookrightarrow \mathfrak{N}$ mit $j(\mathfrak{G}) = j(F) \cap \mathfrak{N}$.
\\
Zwei Einbettungen $j, j':  \mathfrak{G} \hookrightarrow \mathfrak{N}$
mit $j(\mathfrak{G}) = j(F) \cap \mathfrak{N}$ und $j'(\mathfrak{G}) = j'(F) \cap \mathfrak{N}$
sind genau dann $\mathfrak{N}^\times$-konjugiert, wenn $j' = a \circ j$ f\"ur ein $a \in \mathcal{I}(\mathfrak{N})$,
also genau dann, wenn $j'^{-1}(\mathfrak{N}) = j^{-1}(\mathfrak{N})$ und $a(j') \circ a(j)^{-1} \in \mathcal{I}(\mathfrak{N})$.
Damit folgt die erste Behauptung.
\\
Sei $\mathcal{I}_1(\mathfrak{N})$ die Gruppe der Konjugationen von $\mathfrak{N}$ mit Elementen aus $\Gamma(\mathfrak{N})$.
Offenbar gilt $B(\mathfrak{G}, \mathfrak{N}) = C(\mathfrak{G}) [\mathcal{A}(\mathfrak{N}) : \mathcal{I}_1(\mathfrak{N})] = 
B(\mathfrak{G}, \mathfrak{N}) [\mathcal{I}(\mathfrak{N}) : \mathcal{I}_1(\mathfrak{N})]$
und $\mathcal{I}(\mathfrak{N}) \cong \mathfrak{N}^\times/\mathfrak{o}^\times$
sowie $\mathcal{I}_1(\mathfrak{N}) \cong \Gamma(\mathfrak{N})/\{\pm{1}\}$.
Nach \cite[Theorem 11.6.1]{Ma} gilt $N(\mathfrak{N}^\times) = \mathfrak{o}^\times$,
also induziert die Norm eine exakte Sequenz
$\{1\} \rightarrow \Gamma(\mathfrak{N})/\{\pm{1}\} \rightarrow \mathfrak{N}^\times/\mathfrak{o}^\times
\rightarrow \mathfrak{o}^\times/\mathfrak{o}^{\times(2)} \rightarrow \{1\}$.

\end{proof}

\begin{lemma} \label{Item07-04}

Sei $d \in \mathbb{N}$ quadratfrei, und sei $k = \mathbb{Q}(i\sqrt{d})$.
\\
Sei $M$ eine $k$-Quaternionenalgebra, und sei $\mathfrak{N}$ eine $M$-Ordnung.
\\
Sei $\mathcal{G} \subset P\Gamma(\mathfrak{N})$ eine nichtzyklische maximalendliche Untergruppe.
\begin{enumerate}[(i)]

\item \label{Item07-04.1}
Dann ist $\mu(\mathcal{G}, \mathfrak{N}) =
B_1(\mathfrak{F}(\mathcal{G}), \mathfrak{N}) [\mathcal{A}(\mathfrak{F}(\mathcal{G})) : \mathcal{I}(\mathfrak{F}(\mathcal{G}))]^{-1}$

\item \label{Item07-04.2}
Ist $\mathcal{G}$ eine $3$-Diedergruppe oder Tetraedergruppe, dann ist $[\mathcal{A}(\mathfrak{F}(\mathcal{G})) : \mathcal{I}(\mathfrak{F}(\mathcal{G}))] = 2$.
\\
Ist $\mathcal{G}$ eine $2$-Diedergruppe, dann ist $[\mathcal{A}(\mathfrak{F}(\mathcal{G})) : \mathcal{I}(\mathfrak{F}(\mathcal{G}))] = 6$.

\end{enumerate}

\end{lemma}

\begin{proof}

Nach Lemma \ref{Item06-03} gilt $\mathfrak{F}(\mathcal{G})^\times = \Gamma(\mathfrak{F}(\mathcal{G})) = P^{-1}(\mathcal{G})$.
\begin{enumerate}[(i)]

\item[(\ref{Item07-04.1})]
Wegen $P^{-1}(\mathcal{G}) = \mathfrak{F}(\mathcal{G})^\times$ induziert die Zuordnungsvorschrift $\mathcal{G}' \mapsto \mathfrak{F}(\mathcal{G}')$ eine Bijektion zwischen den maximalendlichen Untergruppen $\mathcal{G}' \subset P\Gamma(\mathfrak{N})$ vom Isomorphietyp $\mathcal{G}$
und den optimal eingebetteten Ordnungen $\mathfrak{F}' \subset \mathfrak{N}$ vom Isomorphietyp $\mathfrak{F}(\mathcal{G})$.
\\
Zwei solche Gruppen $\mathcal{G}'$, $\mathcal{G}''$ sind genau dann $P\Gamma(\mathfrak{N})$-konjugiert,
wenn $P^{-1}(\mathcal{G}')$ und $P^{-1}(\mathcal{G}'')$ $\Gamma(\mathfrak{N})$-konjugiert sind,
also wenn $\mathfrak{F}(\mathcal{G}')$ und $\mathfrak{F}(\mathcal{G}'')$ $\Gamma(\mathfrak{N})$-konjugiert sind.
\\
Also ist $\mu(\mathcal{G}, \mathfrak{N}) =
B_1(\mathfrak{F}(\mathcal{G}), \mathfrak{N})[\mathcal{A}'(\mathcal{G}) : \mathcal{I}'(\mathcal{G})]^{-1}$,
wo $\mathcal{A}'(\mathcal{G})$ die Gruppe der Einbettungen $j: \mathfrak{F}(\mathcal{G}) \hookrightarrow \mathfrak{N}$
mit $j(\mathfrak{F}(\mathcal{G})) = \mathfrak{F}(\mathcal{G})$ ist,
und $\mathcal{I}'(\mathcal{G})$ die Untergruppe der $j \in \mathcal{A}'(\mathcal{G})$,
die durch Konjugation mit einem Element aus $\Gamma(\mathfrak{N})$ erzeugt werden.
\\
Offenbar gilt $\mathcal{A}'(\mathcal{G}) \subset \mathcal{A}(\mathfrak{F}(\mathcal{G}))$
und $\mathcal{I}(\mathfrak{F}(\mathcal{G})) \subset \mathcal{I}'(\mathcal{G})$.
Ein $j \in \mathcal{A}(\mathfrak{F}(\mathcal{G}))$ induziert genau
einen Automorphismus von $M$ mit $j(\mathfrak{F}(\mathcal{G})) = \mathfrak{F}(\mathcal{G})$,
also gilt $\mathcal{A}'(\mathcal{G}) = \mathcal{A}(\mathfrak{F}(\mathcal{G}))$.
Falls $J \in \Gamma(\mathfrak{N})$ mit $J\mathfrak{F}(\mathcal{G})J^{-1} = \mathfrak{F}(\mathcal{G})$,
so ist wegen $P^{-1}(\mathcal{G}) = \Gamma(\mathfrak{F}(\mathcal{G}))$ auch $JP^{-1}(\mathcal{G})J^{-1} = P^{-1}(\mathcal{G})$.
Weil $P^{-1}(\mathcal{G}) \subset \Gamma(\mathfrak{N})$ maximalendlich, also sein eigener Normalisator in $\Gamma(\mathfrak{N})$ ist,
folgt $J \in P^{-1}(\mathcal{G}) = \mathfrak{F}(\mathcal{G})^\times$.
Daher gilt $\mathcal{I}'(\mathcal{G}) = \mathcal{I}(\mathfrak{F}(\mathcal{G}))$.

\item[(\ref{Item07-04.2})]
$\mathfrak{F}(\mathcal{G})^\times$ hat das Zentrum $\{\pm 1\}$,
also ist $\mathcal{I}(\mathfrak{F}(\mathcal{G})) \cong \mathfrak{F}(\mathcal{G})^\times / \{\pm 1\} \cong \mathcal{G}$.
\\
$\mathfrak{F}(\mathcal{D}_3)^\times = P^{-1}(\mathcal{D}_3)$ wird erzeugt von $U$ und $V$ mit $U^3 = -1$ und $VUV^{-1} = U^{-1}$.
\\
$j \in \mathcal{A}(\mathfrak{F}(\mathcal{D}_3))$ ist durch $j(U) \in \{U, U^{-1}\}$, $j(V) \in \{\pm V, \pm UV, \pm U^{-1}V \}$ definierbar.
\\
$\mathcal{A}(\mathfrak{F}(\mathcal{D}_3))$ enth\"alt also $12$ Elemente,
und $\mathcal{I}(\mathfrak{F}(\mathcal{D}_3)) \cong \mathcal{D}_3 \cong \mathcal{S}_3$ enth\"alt $6$ Elemente.
\\
$\mathfrak{F}(\mathcal{D}_2)^\times = P^{-1}(\mathcal{D}_2)$ wird erzeugt von $U$ und $V$ mit $U^2 = -1$ und $VUV^{-1} = U^{-1}$.
\\
$j \in \mathcal{A}(\mathfrak{F}(\mathcal{D}_2))$  ist durch $j(U)$ und $j(V)$ definierbar. Daf\"ur gibt es $6 \times 4$ M\"oglichkeiten.
\\
$\mathcal{A}(\mathfrak{F}(\mathcal{D}_2))$ enth\"alt also $24$ Elemente,
und $\mathcal{I}(\mathfrak{F}(\mathcal{D}_2)) \cong \mathcal{D}_2 \cong \mathcal{V}_4$ enth\"alt $4$ Elemente.
\\
$\mathfrak{F}(\mathcal{T})^\times = P^{-1}(\mathcal{T}) \supset P^{-1}(\mathcal{D}_2)$
wird erzeugt von $U$, $V$ und $W = (1- U - V - UV)/2$.
\\
$j \in \mathcal{A}(\mathfrak{F}(\mathcal{T}))$  ist durch $j(U)$ und $j(V)$ definierbar. Daf\"ur gibt es $6 \times 4$ M\"oglichkeiten.
\\
$\mathcal{A}(\mathfrak{F}(\mathcal{T}))$ enth\"alt also $24$ Elemente,
und $\mathcal{I}(\mathfrak{F}(\mathcal{T})) \cong \mathcal{T} \cong \mathcal{A}_4$ enth\"alt $12$ Elemente.

\end{enumerate}

\end{proof}

\begin{theorem} \label{Item07-05}

Sei $d \in \mathbb{N}$ quadratfrei, sei $k = \mathbb{Q}(i\sqrt{d})$ mit Diskriminante $D$,
und bezeichne $t$ die Anzahl der verschiedenen Primteiler von $D$.
\\
Sei $M$ eine $k$-Quaternionenalgebra, und sei $\mathfrak{N}$ eine $M$-Maximalordnung.
\begin{enumerate}[(i)]

\item \label{Item07-05.1}
Sei $\mathcal{D}_3 \subset P\Gamma(\mathfrak{N})$ eine  $3$-Diedergruppe.
Dann hat $\mu(\mathcal{D}_3, \mathfrak{N})$ den Wert
\\
$2^{t-1}$, falls $d \equiv 0 \pmod{3}$
\\
$2^t$, falls $d \equiv 1 \pmod{3}$
\\
$2^t$, falls $d \equiv 2 \pmod{3}$, und $d$ einen Primteiler $p \equiv \pm 5 \pmod{12}$ hat
\\
$2^{t+1}$, falls $d \equiv 2 \pmod{3}$, und $p \equiv \pm 1 \pmod{12}$ f\"ur alle Primteiler $p \ne 2$ von $d$

\item \label{Item07-05.2}
Sei $\mathcal{T} \subset P\Gamma(\mathfrak{N})$ eine  Tetraedergruppe.
Dann hat $\mu(\mathcal{T}, \mathfrak{N})$ den Wert
\\
$2^{t-1}$, falls $d \not\equiv 3 \pmod{4}$
\\
$2^t$, falls $d \equiv 3 \pmod{8}$
\\
$2^t$, falls $d \equiv 7 \pmod{8}$, und $d$ einen Primteiler $p \equiv \pm 3 \pmod{8}$ hat
\\
$2^{t+1}$, falls $d \equiv 7 \pmod{8}$, und $p \equiv \pm 1 \pmod{8}$ f\"ur alle Primteiler $p$ von $d$

\item \label{Item07-05.3}
Sei $\mathcal{D}_2 \subset P\Gamma(\mathfrak{N})$ eine maximalendliche $2$-Diedergruppe.
\\
Dann ist $d \not\equiv 3 \pmod{4}$, und $\mu(\mathcal{D}_2, \mathfrak{N})$ hat den Wert $2^{t-1}$.

\end{enumerate}

\end{theorem}

\noindent
$Bemerkung$: Satz \ref{Item07-05} geht schon aus \cite[S\"atze  20.39 und 26.12]{Kr} hervor,
wird dort allerdings unter R\"uckgriff auf $P\Gamma(\mathfrak{N})$-Konjugationsklassen zyklischer Untergruppen bewiesen.

\begin{proof}

Wir berechnen $\mu(\mathcal{G}, \mathfrak{N})$ mit Lemma \ref{Item07-04} mithilfe von Satz \ref{Item07-03}.
Zur Berechnung der Faktoren $C(\mathfrak{F}(\mathcal{G}))$ und $[\mathcal{A}(\mathfrak{N}) : \mathcal{I}(\mathfrak{N})]$
greifen wir auf Satz \ref{Item04-09} und Lemma \ref{Item07-02} zur\"uck.
\begin{enumerate}[(i)]

\item[(\ref{Item07-05.1})]
$\mu(\mathcal{D}_3, \mathfrak{N}) = C(\mathfrak{F}(\mathcal{D}_3))[\mathcal{A}(\mathfrak{N}) : \mathcal{I}(\mathfrak{N})]$.
Es ist $\Lambda(\mathfrak{F}(\mathcal{D}_3)) = 1$ und $F(\mathcal{D}_3)_3 \not\cong M_2(\mathbb{Q}_3)$.
\\
Falls $d \equiv 1 \pmod{3}$, ist $3$ tr\"age in $k$, also $C(\mathfrak{F}(\mathcal{D}_3)) = 2$.
Sonst ist $C(\mathfrak{F}(\mathcal{D}_3)) = 1$.
\\
Falls $d \not\equiv 2 \pmod{3}$, ist $\Sigma_k(F(\mathcal{D}_3)) = 1$.
Daher gilt dann $[\mathcal{A}(\mathfrak{N}) : \mathcal{I}(\mathfrak{N})] = 2^{t-1}$.
\\
Falls $d \equiv 2 \pmod{3}$, ist $\Sigma_k(F(\mathcal{D}_3)) = 3$.
Wenn in diesem Fall $\left(\dfrac{3,-d}{p}\right) = 1$ an allen Stellen $p$ von $\mathbb{Q}$ gilt,
dann ist $[\mathcal{A}(\mathfrak{N}) : \mathcal{I}(\mathfrak{N})] = 2^{t+1}$,
sonst ist $[\mathcal{A}(\mathfrak{N}) : \mathcal{I}(\mathfrak{N})] = 2^t$.
\\
Die Behauptungen folgen jetzt durch Auswertung der Hilbertsymbole.

\item[(\ref{Item07-05.2})]
$\mu(\mathcal{T}, \mathfrak{N}) = C(\mathfrak{F}(\mathcal{T}))[\mathcal{A}(\mathfrak{N}) : \mathcal{I}(\mathfrak{N})]$.
Es ist $\Lambda(\mathfrak{F}(\mathcal{T})) = 1$ und $F(\mathcal{T})_2 \not\cong M_2(\mathbb{Q}_2)$.
\\
Falls $d \equiv 3 \pmod{8}$, ist $2$ tr\"age in $k$, also $C(\mathfrak{F}(\mathcal{T})) = 2$. Sonst ist $C(\mathfrak{F}(\mathcal{T})) = 1$.
\\
Falls $d \not\equiv 7 \pmod{8}$, ist $\Sigma_k(F(\mathcal{T})) = 1$.
Daher gilt dann $[\mathcal{A}(\mathfrak{N}) : \mathcal{I}(\mathfrak{N})] = 2^{t-1}$.
\\
Falls $d \equiv 7 \pmod{8}$, ist $\Sigma_k(F(\mathcal{T})) = 2$.
Wenn in diesem Fall $\left(\dfrac{2,-d}{p}\right) = 1$ an allen Stellen $p$ von $\mathbb{Q}$ gilt,
dann ist $[\mathcal{A}(\mathfrak{N}) : \mathcal{I}(\mathfrak{N})] = 2^{t+1}$,
sonst ist $[\mathcal{A}(\mathfrak{N}) : \mathcal{I}(\mathfrak{N})] = 2^t$.
\\
Die Behauptungen folgen jetzt durch Auswertung der Hilbertsymbole.

\item[(\ref{Item07-05.3})]
$\mu(\mathcal{D}_2, \mathfrak{N}) = C(\mathfrak{F}(\mathcal{D}_2))[\mathcal{A}(\mathfrak{N}) : \mathcal{I}(\mathfrak{N})]/3$.
Es ist $\Lambda(\mathfrak{F}(\mathcal{D}_2)) = 2$ und $F(\mathcal{D}_2)_2 \not\cong M_2(\mathbb{Q}_2)$.
\\
Nach Satz \ref{Item06-06} ist $d \not\equiv 3 \pmod{4}$, also ist $2$ verzweigt in $k$, also gilt $C(\mathfrak{F}(\mathcal{D}_2)) = 3$
\\
und $\Sigma_k(F(\mathcal{D}_2)) = 1$.
Daher gilt $[\mathcal{A}(\mathfrak{N}) : \mathcal{I}(\mathfrak{N})] = 2^{t-1}$.

\end{enumerate}

\end{proof}


\section{Einige Durchschnitte nichtzyklischer endlicher Gruppen} \label{Chapter08}

\begin{lemma} \label{Item08-01}

Sei $d \in \mathbb{N}$ quadratfrei, $d \not\equiv 3 \pmod 4$.
Sei $k = \mathbb{Q}(i\sqrt{d})$ mit der Hauptordnung $\mathfrak{o}$,
und sei $\pi \in \mathfrak{o}_2$ ein Primelement.
Sei $M$ eine $k$-Quaternionenalgebra.
Sei $\mathcal{D}_2 \subset P\Gamma(M)$ eine $2$-Diedergruppe,
und sei $\mathcal{T}$ die Tetraedergruppe mit $\mathcal{D}_2 \subset \mathcal{T} \subset P\Gamma(M)$.
\\
Seien $U, V$ Erzeugende von $P^{-1}(\mathcal{D}_2)$ mit $U^2 = -1$ und $VUV^{-1} = U^{-1}$,
und seien $U, V$ und $W := (1-U-V-UV)$ Erzeugende von $P^{-1}(\mathcal{T})$.
Dann gilt $M \cong M_2(k)$, sowie:
\begin{enumerate}[(i)]

\item \label{Item08-01.1}
Sei $\mathfrak{N}$ die $M$-Maximalordnung  mit $\mathfrak{F}(\mathcal{T}) \subset \mathfrak{N}$.
Dann ist $\{ 1, (1+U)/\pi, (1+V)/\pi, W \}$ eine $\mathfrak{o}_2$-Basis von $\mathfrak{N}_2$.
F\"ur alle $X \in \{U, V, UV \}$ gilt $k[X]_2 \cap \mathfrak{N}_2 = \mathfrak{o}_2[(1+X)/\pi]$.

\item \label{Item08-01.2}
Sei $\mathfrak{N}(U)$ der $\mathfrak{o}$-Modul mit $\mathfrak{o}$-Basis
\\
$\{ 1, (i\sqrt{d} + U)/2, V, (i\sqrt{d}V + UV)/2 \}$, falls $d \equiv 1 \pmod 4$,
\\
$\{ 1, (1 + i\sqrt{d}V + U)/2, V, (i\sqrt{d} + V + UV)/2 \}$, falls $d \equiv 2 \pmod 4$.
\\
Dann sind $\mathfrak{N}(U)$, $\mathfrak{N}(V) := W^{-1}\mathfrak{N}(U)W$ und $\mathfrak{N}(UV) := W^{-1}\mathfrak{N}(V)W$
paarweise
\\
verschiedene $M$-Maximalordnungen.
Ist $\mathfrak{N}'$ eine $M$-Maximalordnung, dann gilt
\\
$\mathfrak{F}(\mathcal{D}_2) =F(\mathcal{D}_2) \cap \mathfrak{N}'$
genau dann, wenn $\mathfrak{N}' = \mathfrak{N}(U)$, $\mathfrak{N}' = \mathfrak{N}(V)$ oder $\mathfrak{N}' = \mathfrak{N}(UV)$.

\item \label{Item08-01.3}
Sei $X \in \{ U, V, UV \}$.
Dann ist $k[X] \cap \mathfrak{N}(X)$ die Hauptordnung von $k[X]$:
\\
$k[X] \cap \mathfrak{N}(X) = \mathfrak{o}[(i\sqrt{d} + X)/2]$, falls $d \equiv 1 \pmod 4$.
\\
$k[X]_2 \cap \mathfrak{N}(X)_2 = \mathfrak{o}_2[(1+X)/\pi]$, falls $d \equiv 2 \pmod 4$.
\\
F\"ur $Y \in \{ U, V, UV \}$ mit $Y \ne X$ gilt dann $k[Y] \cap \mathfrak{N}(X) = \mathfrak{o}[Y]$.

\item \label{Item08-01.4}
$\Lambda(\mathfrak{N} \cap \mathfrak{N}(U)) =
\Lambda(\mathfrak{N} \cap \mathfrak{N}(V)) = \Lambda(\mathfrak{N} \cap \mathfrak{N}(UV)) = 2$.
Ist $\mathfrak{N}'$ eine $M$-Maximalordnung mit $\Lambda(\mathfrak{N} \cap \mathfrak{N}') = 2$,
so gilt umgekehrt $\mathfrak{N}' = \mathfrak{N}(U)$, $\mathfrak{N}' = \mathfrak{N}(V)$ oder $\mathfrak{N}' = \mathfrak{N}(UV)$.
\\
Weiter gilt $\Lambda(\mathfrak{N}(U) \cap \mathfrak{N}(V)) =
\Lambda(\mathfrak{N}(V) \cap \mathfrak{N}(UV)) = \Lambda(\mathfrak{N}(UV) \cap \mathfrak{N}(U)) = 4$.

\end{enumerate}

\end{lemma}

\begin{proof}

Nach \ref{Item06-04}.(\ref{Item06-04.2}) gilt $M \cong M_2(k)$. Nach Satz \ref{Item04-09} gibt es
genau eine $M$-Maximalord\-nung $\mathfrak{N}$ mit $\mathfrak{F}(\mathcal{T}) \subset \mathfrak{N}$
und genau drei $M$-Maxmalordnungen $\mathfrak{N}'$ mit $\mathfrak{F}(\mathcal{D}_2) = F(\mathcal{D}_2) \cap \mathfrak{N}'$.
\\
Die Aussagen in (\ref{Item08-01.1})-(\ref{Item08-01.3}) lassen sich elementar nachpr\"ufen.
Die erste Aussage in (\ref{Item08-01.4}) folgt mit Lemma \ref{Item04-08}.(\ref{Item04-08.2}),
die zweite Aussage l\"asst sich leicht mit \cite[Theorem 6.5.3.3]{Ma} zeigen,
und die dritte Aussage folgt mit $\pi W \in \mathfrak{N}'_2 \smallsetminus \pi\mathfrak{N}'_2$
aus Lemma \ref{Item02-02}.(\ref{Item02-02.3})-(\ref{Item02-02.4}).

\end{proof}

\begin{lemma} \label{Item08-02}

Sei $d \in \mathbb{N}$ quadratfrei, $d \equiv 2 \pmod 4$. Sei $k = \mathbb{Q}(i\sqrt{d})$ mit der Hauptordnung $\mathfrak{o}$.
Sei $U \in SL_2(\mathbb{C})$ mit $U^2 = -1$, und sei $K := k(U)$ mit der Hauptordnung $\mathfrak{O}$.
\\
Dann gilt $[\Gamma(\mathfrak{O}) : \Gamma(\mathfrak{o}[U])] \le 2$,
und insbesondere gilt $\Gamma(\mathfrak{O})^{(2)} \subset \Gamma(\mathfrak{o}[U])$.
\\
Weiter gilt $[\Gamma(\mathfrak{O}) : \Gamma(\mathfrak{o}[U])] = 2$ genau dann,
wenn es $x, y \in \mathbb{Z}$ gibt mit $x^2 - dy^2 = 2$.

\end{lemma}

\begin{proof}

Man sieht leicht ein (siehe Lemma \ref{Item06-02}), dass $U$ die Gruppe $\mathcal{U}$ der Einheitswurzeln von $K$ mit Norm $1$ erzeugt.
Sei $k_+ := \mathbb{Q}(i\sqrt{d}U)$ mit Hauptordnung $\mathfrak{o}_+ = \mathbb{Z}[i\sqrt{d}U] \subset \mathfrak{o}[U]$.
Seien $\epsilon$ eine Grundeinheit von $K$ und $\epsilon_+$ eine Grundeinheit von $k_+$.
\\[4pt]
Falls $d = 2$, wird $\mathfrak{O}^\times$ von $U': = i\sqrt{2}(1-U)/2$ und $\epsilon_+ = 1 + i\sqrt{2}U$ erzeugt, siehe \cite[§ 26 (9)]{Ha}.
Daher wird $\Gamma(\mathfrak{O})$ von $U = U'^2$ und $U'\epsilon_+$ erzeugt.
Wegen $U'\epsilon_+ \not\in \Gamma(\mathfrak{o}[U])$, aber $(U'\epsilon_+)^2 \in \Gamma(\mathfrak{o}[U])$
gilt also $[\Gamma(\mathfrak{O}) : \Gamma(\mathfrak{o}[U])] = 2$,
und mit $x = 2$, $y=1$ ist $x^2-dy^2 = 2$.
\\[4pt]
Falls $d \ne 2$, ist $\mathcal{U}$ die Gruppe aller Einheitswurzeln von $K$.
Nach \cite[§ 20, Satz 14]{Ha} gilt
$[\Gamma(\mathfrak{O}) : \Gamma(\mathfrak{o}[U])] \le [\Gamma(\mathfrak{O}) : \mathcal{U}\Gamma(\mathfrak{o}_+)] \le
[\mathfrak{O}^\times : \mathcal{U}\mathfrak{o}_+^\times] \le 2$.
Es ist $\epsilon = (a + bU)/2$ mit $a, b \in \mathfrak{o}$
und $a = \alpha + \alpha'i\sqrt{d}$, $\epsilon_+ = \gamma + \gamma'i\sqrt{d}U$
mit $\alpha, \alpha', \gamma, \gamma' \in \mathbb{Z}$.
Falls $[\mathfrak{O}^\times : \mathcal{U}\mathfrak{o}_+^\times] =  2$,
nehmen wir $\epsilon_+ = \epsilon^2U \in \Gamma(\mathfrak{o}_+)$ an, siehe \cite[§ 26 (8)]{Ha}.
Es folgt $\gamma 'i\sqrt{d} - \gamma U = \epsilon^2 = (a^2  -b^2 + abU)/4$,
also $2N(\epsilon) + 2\gamma 'i\sqrt{d} = a^2 = (\alpha^2 - d\alpha'^2) + 2\alpha\alpha'i\sqrt{d}$,
und speziell $\alpha^2 - d\alpha'^2 = 2N(\epsilon)$.
\begin{itemize}

\item
Sei $[\Gamma(\mathfrak{O}) : \Gamma(\mathfrak{o}[U])] = 2$.
Dann gilt auch $[\mathfrak{O}^\times : \mathcal{U}\mathfrak{o}_+^\times] = 2$.
W\"are $N(\epsilon) = -1$,
so w\"urde $\Gamma(\mathfrak{O})$ von $U$ und $\epsilon^2 = -\epsilon_+ U$ erzeugt,
also w\"are $\Gamma(\mathfrak{O}) = \mathcal{U}\Gamma(\mathfrak{o}_+)$, Widerspruch.
Daher gilt $N(\epsilon) = 1$, also $\alpha^2 - d\alpha'^2 = 2$.

\item
Sei umgekehrt $N(x+yi\sqrt{d}U) = x^2 - dy^2 = 2$ mit $x, y \in \mathbb{Z}$.
Nach \cite[§ 26 (12$_\text{II}$)]{Ha} gilt dann $[\mathfrak{O}^\times : \mathcal{U}\mathfrak{o}_+^\times] = 2$.
W\"are $N(\epsilon) = -1$, so w\"are $N(\alpha+\alpha' i\sqrt{d}U)) = -2$, also $N(e_+) = -1$, Widerspruch zu \cite[§ 26 (9)]{Ha}.
Also ist $N(\epsilon) = 1$, und offenbar ist $\alpha'$ ungerade.
Damit erkennt man leicht, dass $\epsilon \not\in \Gamma(\mathfrak{o}[U])$.
Also ist $[\Gamma(\mathfrak{O}) : \Gamma(\mathfrak{o}[U])] = 2$.

\end{itemize}

\end{proof}

\begin{theorem} \label{Item08-03}

Sei $d \in \mathbb{N}$ quadratfrei, sei $d \ne 1$, und sei $k = \mathbb{Q}(i\sqrt{d})$.
\\
Seien $M$ eine $k$-Quaternionenalgebra und $\mathfrak{M}$ eine $M$-Maximalordnung.
\\
Sei $\mathcal{T} \subset P\Gamma(\mathfrak{M})$ eine Tetraedergruppe,
und sei $\mathcal{C}_2 \subset \mathcal{T}$ eine Untergruppe der Ordnung $2$.
\begin{enumerate}[(i)]

\item \label{Item08-03.1}
Sei $d \equiv 2 \pmod 4$, und es gebe $x, y \in \mathbb{Z}$ mit $x^2 - dy^2 = 2$.
Dann gibt es eine maximalendliche $2$-Diedergruppe $\mathcal{D}'_2 \subset P\Gamma(\mathfrak{M})$ mit $\mathcal{C}_2  \subset \mathcal{D}'_2$.
Jede Tetraedergruppe $\mathcal{T}'' \subset P\Gamma(\mathfrak{M})$ mit $\mathcal{C}_2  \subset \mathcal{T}''$
ist $P\Gamma(\mathfrak{M})$-konjugiert zu $\mathcal{T}$.
Jede maximalendliche $2$-Diedergruppe $\mathcal{D}''_2 \subset P\Gamma(\mathfrak{M})$ mit $\mathcal{C}_2  \subset \mathcal{D}''_2$
ist $P\Gamma(\mathfrak{M})$-konjugiert zu $\mathcal{D}'_2$.

\item \label{Item08-03.2}
Sei $d \not\equiv 2 \pmod 4$, oder f\"ur alle $x, y \in \mathbb{Z}$ sei $x^2 - dy^2 \ne 2$.
Dann gibt es eine Tetraedergruppe $\mathcal{T}' \subset P\Gamma(\mathfrak{M})$ mit $\mathcal{C}_2  \subset \mathcal{T}'$,
die nicht $P\Gamma(\mathfrak{M})$-konjugiert zu $\mathcal{T}$ ist.
\\
Jede Tetraedergruppe $\mathcal{T}'' \subset P\Gamma(\mathfrak{M})$ mit $\mathcal{C}_2  \subset \mathcal{T}''$
ist $P\Gamma(\mathfrak{M})$-konjugiert zu $\mathcal{T}$ oder zu $\mathcal{T}'$.
Es gibt keine maximalendliche $2$-Diedergruppe $\mathcal{D}'_2 \subset P\Gamma(\mathfrak{M})$ mit $\mathcal{C}_2  \subset \mathcal{D}'_2$.

\end{enumerate}

\end{theorem}

\begin{proof}

Sei $\mathcal{D}_2$ die $2$-Diedergruppe mit $\mathcal{C}_2 \subset \mathcal{D}_2 \subset \mathcal{T}$.
Wir zeigen (\ref{Item08-03.0.1}),
dass es eine \mbox{$2$-Diedergruppe} $\mathcal{D}'_2 \subset P\Gamma(\mathfrak{M})$ mit $\mathcal{C}_2  \subset \mathcal{D}'_2$ gibt,
die nicht $P\Gamma(\mathfrak{M})$-konjugiert zu $\mathcal{D}_2$ ist,
und dass jede $2$-Diedergruppe $\mathcal{D}''_2 \subset P\Gamma(\mathfrak{M})$ mit $\mathcal{C}_2  \subset \mathcal{D}''_2$
dann $P\Gamma(\mathfrak{M})$-konjugiert zu $\mathcal{D}_2$ oder $\mathcal{D}'_2$ ist.
Wir m\"ussen danach nur noch kl\"aren (\ref{Item08-03.0.2}),
ob $\mathcal{D}'_2 \subset P\Gamma(\mathfrak{M})$ eine maximalend\-liche Untergruppe ist,
oder ob es eine Tetraedergruppe $\mathcal{T}' \subset P\Gamma(\mathfrak{M})$ mit $\mathcal{D}'_2 \subset \mathcal{T}'$ gibt.
\begin{enumerate}[(a)]

\item \label{Item08-03.0.1}
$P^{-1}(\mathcal{D}_2)$ wird von $U \in P^{-1}(\mathcal{C}_2)$ und $V$ mit $U^2 = -1$ und $VUV^{-1} = U^{-1}$ erzeugt.
\\
Ist $\mathcal{D}'_2$ eine $2$-Diedergruppe mit $\mathcal{C}_2 \subset \mathcal{D}'_2 \subset P\Gamma(\mathfrak{M})$,
dann wird $P^{-1}(\mathcal{D}'_2)$ von $U$ und $V'$ mit $V'UV'^{-1} = U^{-1}$ erzeugt.
Daraus folgt $X': = V^{-1}V' \in k(U) \cap \Gamma(\mathfrak{M})$.
\\
$V'$ und $X'$ sind durch $\mathcal{D}'_2$ bis auf einen Faktor aus $P^{-1}(\mathcal{C}_2)$ eindeutig bestimmt.
\\
Ist umgekehrt $X' \in k(U) \cap \Gamma(\mathfrak{M})$ und $V' := VX'$,
so ist $V'UV'^{-1} = U^{-1}$, also sind $U$ und $V'$ Erzeugende von $P^{-1}(\mathcal{D}'_2)$ f\"ur eine 2-Diedergruppe $\mathcal{D}'_2$
mit $\mathcal{C}_2 \subset \mathcal{D}'_2 \subset P\Gamma(\mathfrak{M})$.

$\mathcal{D}_2$ und $\mathcal{D}'_2$ sind genau dann $P\Gamma(\mathfrak{M})$-konjugiert,
wenn es ein $U' \in P^{-1}(\mathcal{C}_2)$ und ein $J \in \Gamma(\mathfrak{M})$ gibt mit $JUJ^{-1} = U$ und $JV'J^{-1} = VU'$.
\\
$Zwischenbemerkung$: Wir k\"onnen o.B.d.A. $JUJ^{-1} = U$ annehmen, da die Elemente $X \in P^{-1}(\mathcal{D}_2)$ mit $X^2 = -1$
in $P^{-1}(\mathcal{T})$ jeweils paarweise zueinander konjugiert sind.
\\
Dann ist $J \in k(U) \cap \Gamma(\mathfrak{M})$, und $X' = V^{-1}V' = V^{-1}J^{-1}VU'J = (J^{-1})^*U'J =  J^2U'$.
\\
Ist umgekehrt $X' = J^2U'$ mit $U' \in P^{-1}(\mathcal{C}_2)$ und $J \in k(U) \cap \Gamma(\mathfrak{M})$,
so ist $JUJ^{-1} = U$ und $JV'J^{-1} = JVX'J^{-1} = JVJ^2U'J^{-1} = VJ^*J^2U'J^{-1} = VU'$.
\\
Wegen $d \ne 1$ ist $k(U)$ algebraischer Zahlk\"orper mit zwei komplexen Stellen,
Nach dem Dirichletschen Einheitensatz ist $k(U) \cap \Gamma(\mathfrak{M})$ daher direktes Produkt
einer zyklischen unendlichen Gruppe und der Untergruppe $P^{-1}(\mathcal{C}_2)$ der Einheitswurzeln von $k(U)$ mit Norm $1$.
Deshalb gilt $[(k(U) \cap \Gamma(\mathfrak{M})) : (k(U) \cap \Gamma(\mathfrak{M}))^{(2)}P^{-1}(\mathcal{C}_2)] = 2$.
Sei jetzt $\mathcal{D}'_2$ eine $2$-Diedergruppe mit $\mathcal{C}_2 \subset \mathcal{D}'_2 \subset P\Gamma(\mathfrak{M})$
und $X' \notin (k(U) \cap \Gamma(\mathfrak{M}))^{(2)}P^{-1}(\mathcal{C}_2)$.
Dann ist $\mathcal{D}'_2$ nicht $P\Gamma(\mathfrak{M})$-konjugiert zu $\mathcal{D}_2$.
Ist nun $\mathcal{D}''_2$ eine $2$-Diedergruppe mit $\mathcal{C}_2 \subset \mathcal{D}''_2 \subset P\Gamma(\mathfrak{M})$,
so wird $P^{-1}(\mathcal{D}''_2)$ von $U$ und $V''$ mit $V''UV''^{-1} = U^{-1}$ erzeugt.
Wegen $V'^{-1}V'' = (V^{-1}V')^{-1}V^{-1}V''$
gilt $V^{-1}V'' \in (k(U) \cap \Gamma(\mathfrak{M}))^{(2)}P^{-1}(\mathcal{C}_2)$
oder $V'^{-1}V'' \in (k(U) \cap \Gamma(\mathfrak{M}))^{(2)}P^{-1}(\mathcal{C}_2)$,
und $\mathcal{D}''_2$ ist $P\Gamma(\mathfrak{M})$-konjugiert zu $\mathcal{D}_2$ oder $\mathcal{D}'_2$.

\item \label{Item08-03.0.2}
Falls $d \equiv 3 \pmod 4$,
gibt es nach Satz \ref{Item06-06} keine maximalendliche $2$-Diedergruppe $\mathcal{D}'_2 \subset P\Gamma(\mathfrak{M})$.
Also gilt dann $\mathcal{D}'_2 \subset \mathcal{T}'$ f\"ur eine Tetraedergruppe $\mathcal{T}' \subset P\Gamma(\mathfrak{M})$.
\\
Falls $d \not\equiv 3 \pmod 4$,
folgt mit Lemma \ref{Item08-01}.(\ref{Item08-01.1}), dass $(k(U) \cap \mathfrak{M})_2 = \mathfrak{o}_2[(1+U)/\pi]$ gilt.
\\
Falls $d \equiv 1 \pmod 4$,
gibt es keine maximalendliche $2$-Diedergruppe $\mathcal{D}'_2 \subset P\Gamma(\mathfrak{M})$ mit $\mathcal{C}_2  \subset \mathcal{D}'_2$,
denn nach Lemma \ref{Item08-01}.(\ref{Item08-01.3})
w\"are sonst $k(U) \cap \mathfrak{M} = \mathfrak{o}[(i\sqrt{d} + U)/2]$ oder $k(U) \cap \mathfrak{M} = \mathfrak{o}[U]$.
Also gilt $\mathcal{D}'_2 \subset \mathcal{T}'$ f\"ur eine Tetraedergruppe $\mathcal{T}' \subset P\Gamma(\mathfrak{M})$.
\\
Sei nun $d \equiv 2 \pmod 4$.
$P^{-1}(\mathcal{T})$ wird von $U$, $V$ und $W := (1 - U - V - UV)/2$ erzeugt.
Wenn es eine Tetraedergruppe $\mathcal{T}' \subset P\Gamma(\mathfrak{M})$ mit $\mathcal{D}'_2 \subset \mathcal{T}'$ gibt,
wird $P^{-1}(\mathcal{T}')$ von $U$, $V'$ und $W' := (1 - U - V' - UV')/2$ erzeugt.
Aus $(1+U)V(1-X')/2 = W'-W \in \mathfrak{M}$
folgt dann $N(1-X') \in 2\mathfrak{o}$ oder gleichwertig $S(X') \in 2\mathfrak{o}$,
und daher $X' \in \Gamma(\mathfrak{o}[U])$.
\\
Sei umgekehrt $X' \in \Gamma(\mathfrak{o}[U])$, also $X' = a + bU$ mit $a, b \in \mathfrak{o}$ und $a^2 + b^2 = 1$.
Wegen $[\mathfrak{o}_2 : \pi \mathfrak{o}_2] = 2$ ist $(1-a+b) \in \pi \mathfrak{o}_2$
und $(1 - X')/\pi = (1-a+b)/\pi - b(1+U)/\pi \in \mathfrak{M}_2$,
also $W' := (1-U-V'-UV')/2 = W + (1+U)V(1-X')/2 \in \mathfrak{M}_2$.
Dann sind $U$, $V'$ und $W'$ Erzeugende von $P^{-1}(\mathcal{T}')$
f\"ur eine Tetraedergruppe $\mathcal{T}'$ mit $\mathcal{C}_2 \subset \mathcal{T}' \subset P\Gamma(\mathfrak{M})$.

Genau dann ist $\mathcal{D}'_2 \subset P\Gamma(\mathfrak{M})$ maximalendliche $2$-Diedergruppe,
wenn $X' \not\in \Gamma(\mathfrak{o}[U])$,
wegen $(k(U) \cap \Gamma(\mathfrak{M}))^{(2)}P^{-1}(\mathcal{C}_2) \subset \Gamma(\mathfrak{o}[U])$
also wenn $[(k(U) \cap \Gamma(\mathfrak{M})) : \Gamma(\mathfrak{o}[U])] = 2$,
nach Lemma \ref{Item08-02} also wenn es $x, y \in \mathbb{Z}$ gibt mit $x^2 - dy^2 = 2$.

\end{enumerate}

\end{proof}

\begin{theorem} \label{Item08-04}

Sei $d \in \mathbb{N}$ quadratfrei, sei $d \ne 1$, und sei $k = \mathbb{Q}(i\sqrt{d})$.
\\
Seien $M$ eine $k$-Quaternionenalgebra und $\mathfrak{M}$ eine $M$-Maximalordnung.
\\
Sei $\mathcal{D}_2 \subset P\Gamma(\mathfrak{M})$ eine maximalendliche $2$-Diedergruppe.
\begin{enumerate}[(i)]

\item \label{Item08-04.1}
Sei $d \equiv 2 \pmod 4$, und es gebe $x, y \in \mathbb{Z}$ mit $x^2 - dy^2 = 2$.
\begin{enumerate}[(a)]

\item \label{Item08-04.1.1}
Dann gibt es eine Untergruppe $\mathcal{C}_2 \subset \mathcal{D}_2$ der Ordnung $2$
und eine Tetraedergruppe $\mathcal{T}' \subset P\Gamma(\mathfrak{M})$ mit $\mathcal{C}_2 \subset \mathcal{T}'$.
Die Zuordnung $\mathcal{D}_2 \mapsto \mathcal{C}_2$ ist eindeutig.
\\
Jede Tetraedergruppe $\mathcal{T}'' \subset P\Gamma(\mathfrak{M})$ mit $\mathcal{C}_2  \subset \mathcal{T}''$
ist $P\Gamma(\mathfrak{M})$-konjugiert zu $\mathcal{T}'$.

\item \label{Item08-04.1.2}
Seien $\mathcal{C}'_2, \mathcal{C}''_2 \subset \mathcal{D}_2$ die beiden anderen Untergruppen der Ordnung $2$.
\\
Dann sind $\mathcal{C}'_2$ und $\mathcal{C}''_2$ nicht zu $\mathcal{C}_2$, aber zueinander $P\Gamma(\mathfrak{M})$-konjugiert.

\item \label{Item08-04.1.3}
Jede maximalendliche $2$-Diedergruppe $\mathcal{D}''_2 \subset P\Gamma(\mathfrak{M})$ mit
\\
$\mathcal{C}_2  \subset \mathcal{D}''_2$, $\mathcal{C}'_2 \subset \mathcal{D}''_2$ oder $\mathcal{C}''_2 \subset \mathcal{D}''_2$
ist $P\Gamma(\mathfrak{M})$-konjugiert zu $\mathcal{D}_2$.

\end{enumerate}

\item \label{Item08-04.2}
Sei $d \not\equiv 2 \pmod 4$, oder f\"ur alle $x, y \in \mathbb{Z}$ sei $x^2 - dy^2 \ne 2$.
\begin{enumerate}[(a)]

\item \label{Item08-04.2.1}
Sei $\mathcal{C}_2  \subset \mathcal{D}_2$ eine Untergruppe der Ordnung $2$.
\\
Dann gibt es keine Tetraedergruppe $\mathcal{T}' \subset P\Gamma(\mathfrak{M})$ mit $\mathcal{C}_2  \subset \mathcal{T}'$.

\item \label{Item08-04.2.2}
Seien $\mathcal{C}'_2, \mathcal{C}''_2 \subset \mathcal{D}_2$ die beiden anderen Untergruppen der Ordnung $2$.
\\
Dann sind $\mathcal{C}_2$, $\mathcal{C}'_2$ und $\mathcal{C}''_2$ paarweise nicht zueinander $P\Gamma(\mathfrak{M})$-konjugiert.

\item \label{Item08-04.2.3}
Es gibt eine maximalendliche $2$-Diedergruppe $\mathcal{D}'_2 \subset P\Gamma(\mathfrak{M})$
\\
mit $\mathcal{C}_2  \subset \mathcal{D}'_2$,
die nicht $P\Gamma(\mathfrak{M})$-konjugiert zu $\mathcal{D}_2$ ist.
\\
Jede maximalendliche $2$-Diedergruppe $\mathcal{D}''_2 \subset P\Gamma(\mathfrak{M})$ mit
\\
$\mathcal{C}_2  \subset \mathcal{D}''_2$, $\mathcal{C}'_2  \subset \mathcal{D}''_2$ oder $\mathcal{C}''_2  \subset \mathcal{D}''_2$
ist $P\Gamma(\mathfrak{M})$-konjugiert zu $\mathcal{D}_2$ oder $\mathcal{D}'_2$.

\end{enumerate}

\end{enumerate}

\end{theorem}

\begin{proof}

Nach Satz \ref{Item06-06} ist $d \not\equiv 3 \pmod 4$, also $2$ verzweigt in $k$. Sei $\pi \in \mathfrak{o}_2$ Primelement.
In den Beweisabschnitten (\ref{Item08-04.0.1}) - (\ref{Item08-04.0.3}) zeigen wir jeweils beide Teile (\ref{Item08-04.1}) und (\ref{Item08-04.2}) des Satzes.
\begin{enumerate}[(a)]

\item \label{Item08-04.0.1}
Sei $d \equiv 2 \pmod 4$, und es gebe $x, y \in \mathbb{Z}$ mit $x^2 - dy^2 = 2$.
\\
Daraus folgt $\left(\dfrac{2,-d}{p}\right) = \left(\dfrac{2,d}{p}\right) = 1$ f\"ur alle Stellen $p \ne 2, \infty$ von $\mathbb{Q}$.
Nach Satz \ref{Item06-07} gibt es daher auch eine Tetraedergruppe $\mathcal{T} \subset P\Gamma(\mathfrak{M})$.
Nach Satz \ref{Item08-03}(\ref{Item08-03.1})
gibt es eine maximalendliche $2$-Diedergruppe $\mathcal{D}'_2 \subset P\Gamma(\mathfrak{M})$,
so dass $\mathcal{T} \cap \mathcal{D}'_2$ die Ordnung $2$ hat.
Die $P\Gamma(\mathfrak{M})$-Konjugationsklasse von $\mathcal{D}'_2$
ist durch die $P\Gamma(\mathfrak{M})$-Konjugationsklasse von $\mathcal{T}$ eineindeutig bestimmt.
Nach Satz \ref{Item07-05} gilt $\mu(\mathcal{D}_2, \mathfrak{M}) = \mu(\mathcal{T}, \mathfrak{M})$.
Also gibt es auch eine Tetraedergruppe $\mathcal{T}' \subset P\Gamma(\mathfrak{M})$,
so dass $\mathcal{C}_2 := \mathcal{T}' \cap \mathcal{D}_2$ die Ordnung $2$ hat.
\\
Sei $X \in P^{-1}(\mathcal{C}_2)$ mit $X^2 = -1$.
Nach Lemma \ref{Item08-01}.(\ref{Item08-01.1}) gilt $(k(X) \cap \mathfrak{M})_2 = \mathfrak{o}_2[(1+X)/\pi]$,
Nach Lemma \ref{Item08-01}.(\ref{Item08-01.3}) ist dadurch $X$ bis aufs Vorzeichen, also $\mathcal{C}_2$ eindeutig bestimmt.

Sei umgekehrt $\mathcal{T}' \subset P\Gamma(\mathfrak{M})$ eine Tetraedergruppe,
so dass $\mathcal{T}' \cap \mathcal{D}_2$ die Ordnung $2$ hat.
Nach Satz \ref{Item08-03} ist dann $d \equiv 2 \pmod 4$, und es gibt $x, y \in \mathbb{Z}$ mit $x^2 - dy^2 = 2$.

\item \label{Item08-04.0.2}
$P^{-1}(\mathcal{D}_2)$ wird von $U$ und $V$ mit $U^2 = -1$ und $VUV^{-1} = U^{-1}$ erzeugt.

Sei $d \equiv 2 \pmod 4$, und es gebe $x, y \in \mathbb{Z}$ mit $x^2 - dy^2 = 2$. Dann ist $x$ gerade.
\\
Wir nehmen o.b.d. A. an, dass $U \in P^{-1}(\mathcal{C}'_2)$, $V \in P^{-1}(\mathcal{C}''_2)$ und $UV \in P^{-1}(\mathcal{C}_2)$ gilt.
\\
Gem\"a{\ss} Beweisabschnitt (\ref{Item08-04.0.1}) gilt dann $(k(UV) \cap \mathfrak{M})_2 = \mathfrak{o}_2[(1+UV)/\pi]$.
\\
Sei $\epsilon := (x + yi\sqrt{d}U + yi\sqrt{d}V - xUV)/2 = (1 - UV)x/2 + yU(1-UV)i\sqrt{d}/2 \in \mathfrak{M}$
Man rechnet elementar nach, dass $N(\epsilon) = 1$ und $\epsilon V \epsilon^{-1} = U$ gilt.

Seien umgekehrt zwei Untergruppen von $\mathcal{D}_2$ der Ordnung $2$ zueinander
$P\Gamma(\mathfrak{M})$-konjugiert.
Wir k\"onnen o.B.d.A. $\epsilon V \epsilon^{-1} = U$ mit $\epsilon \in \Gamma(\mathfrak{M})$ annehmen.
Dann gilt $k(V) \cap \mathfrak{M} \cong k(U) \cap \mathfrak{M}$
und mit Lemma \ref{Item08-01} folgt $(1+UV)/\pi \in \mathfrak{M}_2$.
Sei $V' = \epsilon U \epsilon^{-1}$,
und sei $\mathcal{D}'_2 \subset P\Gamma(\mathfrak{M})$ die von $P(U)$ und $P(V')$ erzeugte $2$-Diedergruppe.
Wir k\"urzen $F := F(\mathcal{D}_2)$ und $F' := F(\mathcal{D}'_2) = \epsilon F \epsilon^{-1}$ ab.
Sei $T: = \epsilon \Phi_F(\epsilon)^*$.
Dann gilt $\Phi_F(T) = T^*$,
und f\"ur $A \in F$ gilt $T \Phi_F(\epsilon A \epsilon^{-1}) T^{-1} =
\epsilon \Phi_F(\epsilon)^{-1} \Phi_F(\epsilon) A \epsilon^{-1} = \epsilon A \epsilon^{-1}$.
Daher gilt $F' = F(F,T)$ und $\Phi_{F'}(A) = T \Phi_F(A) T^{-1}$ f\"ur $A \in M$.
Mit Lemma \ref{Item08-01} folgt aus $k(V) \cap \Phi_F(\mathfrak{M}) \cong k(U) \cap \Phi_F(\mathfrak{M})$,
dass $\Phi_F(\mathfrak{M}) = \mathfrak{M}$, und analog $\Phi_{F'}(\mathfrak{M}) = \mathfrak{M}$.
Wegen $U \in F \cap F'$ gilt $TUT^{-1} = U$, also $T \in k(U) \cap \mathfrak{M} = \mathfrak{o}[U]$.
Mit $\Phi_F(T) = T^*$ folgt $T = \gamma + \gamma'i\sqrt{d}U$
mit $\gamma, \gamma' \in \mathbb{Z}$ und $\gamma^2 - d\gamma'^2 = N(T) = 1$, $\gamma$ ungerade, $\gamma'$ gerade.
Mit $T \in k(U)$ folgt weiter $\Phi_{F'}(TV) =  T T^* V T^{-1} = TV$ und $(TV)U(TV)^{-1} = U^{-1}$.
Also erzeugen $P(U)$ und $P(TV)$ eine $2$-Diedergruppe $\mathcal{D}''_2$ mit $P^{-1}(\mathcal{D}''_2) \subset F' \cap \mathfrak{M}$.
Mit Lemma \ref{Item06-03} folgt $P^{-1}(\mathcal{D}''_2) \subset \Gamma(F' \cap \mathfrak{M}) =
\Gamma(\mathfrak{F}(\mathcal{D}'_2)) = P^{-1}(\mathcal{D}'_2)$,
also $\mathcal{D}''_2 = \mathcal{D}'_2$.
Wegen $(1 + UTV)/\pi = (1+ \gamma UV)/\pi - \gamma'i\sqrt{d}V/\pi \in \mathfrak{M}_2$ ist notwendig $TV = \pm V'$.
Aus $\epsilon V \epsilon^{-1} = U$ folgt $\epsilon = (a + bU + bV - aUV)/2$ mit $a, b \in \mathfrak{o}$ und $a^2 + b^2 = 2$,
und $V' = \epsilon U \epsilon^{-1} = ((-a^2 + b^2)/2 - abU)V$,
also $\gamma + \gamma'i\sqrt{d}U = \pm (-a^2 + b^2)/2 - abU)$.
Seien $a = \alpha + \alpha'i\sqrt{d}$ und $b = \beta + \beta'i\sqrt{d}$
mit $\alpha, \alpha', \beta, \beta' \in \mathbb{Z}$.
Damit erhalten wir $-\alpha\alpha' + \beta\beta' = 0$ und $\alpha\beta - \alpha'\beta' d = 0$.
Daraus folgt $\beta'(\beta^2 - d\alpha'^2) = \alpha(\beta^2 - d\alpha'^2) = 0$,
wegen $d \ne 1$ also $\alpha = \beta' = 0$ oder $\alpha' = \beta = 0$.
Dann ist $2 = a^2 + b^2 = \beta^2 - d\alpha'^2$ oder $2 = \alpha^2 - d\beta'^2$,
und weiter folgt $d \not\equiv 1 \pmod 4$, also $d \equiv 2 \pmod 4$.

\item \label{Item08-04.0.3}
Zun\"achst zeigen wir die Behauptung unter der Voraussetzung $\mathcal{C}_2  \subset \mathcal{D}''_2$.
\\
Sei $d \equiv 2 \pmod 4$, und es gebe $x, y \in \mathbb{Z}$ mit $x^2 - dy^2 = 2$.
\\
Dann folgt die Behauptung mit Satz \ref{Item08-04}.(\ref{Item08-04.1})(\ref{Item08-04.0.1}) und Satz \ref{Item08-03}.(\ref{Item08-03.1}).
\\
Sei $d \not\equiv 2 \pmod 4$, oder f\"ur alle $x, y \in \mathbb{Z}$ sei $x^2 - dy^2 \ne 2$.
\\
Der Beweis erfolgt dann wie im Beweisabschnitt (\ref{Item08-03.0.1}) zu Satz \ref{Item08-03},
\\
wenn wir dort die Zwischenbemerkung wie folgt ersetzen:
\\
$Zwischenbemerkung$: Wir k\"onnen o.B.d.A. $JUJ^{-1} = U$ annehmen, da $U$  in $\Gamma(\mathfrak{M})$ zu $U^{-1}$,
nach Satz \ref{Item08-04}.(\ref{Item08-04.2})(\ref{Item08-04.0.2}) aber nicht zu $V^{\pm 1}$ und  $(UV)^{\pm 1}$ konjugiert ist.

Nun zeigen wir die Behauptung unter der Voraussetzung $\mathcal{C}'_2  \subset \mathcal{D}''_2$ aus.
Es gen\"ugt zu zeigen, dass es eine $2$-Diedergruppe $\mathcal{D}'_2$
mit $\mathcal{C}_2  \subset \mathcal{D}'_2 \subset P\Gamma(\mathfrak{M})$ gibt,
die $P\Gamma(\mathfrak{M})$-konjugiert zu $\mathcal{D}''_2$ ist.
(Unter der Voraussetzung $\mathcal{C}''_2  \subset \mathcal{D}''_2$ erfolgt der Beweis analog.)
\\
$P^{-1}(\mathcal{D}_2)$ wird von $U, V$ mit $U^2 = -1$ und $VUV^{-1} = U^{-1}$ erzeugt.
O.B.d.A. nehmen wir $U \in P^{-1}(\mathcal{C}'_2)$ mit $k(U) \cap \mathfrak{M} = \mathfrak{o}[U]$
und $V \in P^{-1}(\mathcal{C}''_2)$  mit $k(V) \cap \mathfrak{M} = \mathfrak{o}[V]$ an.
\\
$P^{-1}(\mathcal{D}''_2)$ wird von $U, V'$ mit $V'UV'^{-1} = U^{-1}$ erzeugt.
Dann gilt $V^{-1}V' \in \mathfrak{o}[U]$.
\\
Sei $V' = V(a + bU)$ mit $a, b \in \mathfrak{o}$, $a^2 + b^2 = 1$.
Wir k\"onnen $V'$ durch $UV' = V(b - aU)$ ersetzen, o.B.d.A. also $(a-1), b \in \pi \mathfrak{o}_2$ annehmen.
W\"are $(a-1) \not\in 2 \mathfrak{o}_2$, so w\"are $a^2 = 1 + \pi^2 + 2\pi + 4a'$ mit $a' \in \mathfrak{o}_2$
und $b^2 = \pi^2 + 4b'$ mit $b' \in \mathfrak{o}_2$, Widerspruch.
Jetzt folgt $(a-1), b \in 2 \mathfrak{o}$.
und $\epsilon := ((a + 1) + bU + (a - 1)V - bUV)/2 \in \Gamma(\mathfrak{M})$
sowie $\epsilon UV' \epsilon^{-1} = UV$.
Sei $\mathcal{D}'_2$ die von $P(\epsilon U \epsilon^{-1})$ und $P(\epsilon V' \epsilon^{-1})$ erzeugte $2$-Diedergruppe.

\end{enumerate}

\end{proof}

\begin{theorem} \label{Item08-05}

Sei $d \in \mathbb{N}$ quadratfrei, sei $d \ne 3$, und sei $k = \mathbb{Q}(i\sqrt{d})$.
\\
Seien $M$ eine $k$-Quaternionenalgebra und $\mathfrak{M}$ eine $M$-Maximalordnung.
Sei $\mathcal{D}_3 \subset P\Gamma(\mathfrak{M})$ eine $3$-Diedergruppe,
und sei $\mathcal{C}_3 \subset \mathcal{D}_3$ deren Untergruppe der Ordnung $3$.
Dann gibt es eine $3$-Diedergruppe $\mathcal{D}'_3 \subset P\Gamma(\mathfrak{M})$ mit $\mathcal{C}_3  \subset \mathcal{D}'_3$,
die nicht $P\Gamma(\mathfrak{M})$-konjugiert zu $\mathcal{D}_3$ ist.
Jede $3$-Diedergruppe $\mathcal{D}''_3 \subset P\Gamma(\mathfrak{M})$ mit $\mathcal{C}_3  \subset \mathcal{D}''_3$
ist $P\Gamma(\mathfrak{M})$-konjugiert zu $\mathcal{D}_3$ oder zu $\mathcal{D}'_3$.

\end{theorem}

\begin{proof}

wie im Beweisabschnitt (\ref{Item08-03.0.1}) zu Satz \ref{Item08-03},
wenn wir dort $\mathcal{D}_2$, $\mathcal{C}_2$, $U^2 = -1, d \ne 1$ durch $\mathcal{D}_3$, $\mathcal{C}_3$, $U^3 = -1, d \ne 3$ ersetzen,
sowie die Zwischenbemerkung wie folgt ersetzen:
\\
$Zwischenbemerkung$: Wir k\"onnen o.B.d.A. $JUJ^{-1} = U$ annehmen,
da die beiden Elemente $X = U^{\pm 1}\in P^{-1}(\mathcal{D}_3) \smallsetminus \mathfrak{o}^\times$
mit $X^3 = -1$ in $P^{-1}(\mathcal{D}_3)$ zueinander konjugiert sind.

\end{proof}

\begin{corollary} \label{Item08-06}

Sei $d \in \mathbb{N}$ quadratfrei, und sei $k = \mathbb{Q}(i\sqrt{d})$.
\\
Seien $M$ eine $k$-Quaternionenalgebra und $\mathfrak{M}$ eine $M$-Maximalordnung.
\\
F\"ur $m = 3$ oder $m = 2$ bezeichne $\lambda^*_m(\mathfrak{M})$ die Anzahl der $P\Gamma(\mathfrak{M})$-Konjugationsklassen von Gruppen der Ordnung $m$,
die in einer $m$-Diedergruppe $\mathcal{D}_m \subset P\Gamma(\mathfrak{M})$ enthalten sind.
\\
Dann ist $2 \lambda^*_3(\mathfrak{M}) = \mu(\mathcal{D}_3, \mathfrak{M})$ f\"ur $d \ne 3$,
und $2 \lambda^*_2(\mathfrak{M}) = \mu(\mathcal{T}, \mathfrak{M}) + 3 \mu(\mathcal{D}_2, \mathfrak{M})$ f\"ur $d \ne 1$.

\end{corollary}

\noindent
$Bemerkung$: $\lambda^*_m(\mathfrak{M})$ entspricht der in \cite[§ 15]{Kr} definierten Bezeichnung $\lambda^*_{2m}(\mathfrak{M})$.
\\
Korollar \ref{Item08-06} geht schon aus \cite[S\"atze  20.39 und 26.12, Tabellen 20.40 und 26.13]{Kr} hervor.
Erstmals explizit formuliert und auf andere Art bewiesen wurde er in \cite [Corollary 23]{Ra}.

\begin{proof}

F\"ur $d \ne 3$ geht die Behauptung $2 \lambda^*_3(\mathfrak{M}) = \mu(\mathcal{D}_3, \mathfrak{M})$ direkt aus Satz \ref{Item08-05} hervor.
\\
Sei nun zun\"achst $d \equiv 2 \pmod 4$, und es gebe $x, y \in \mathbb{Z}$ mit $x^2 - dy^2 = 2$.
\\
Nach Satz \ref{Item08-03}.(\ref{Item08-03.1}) und \ref{Item08-04}.(\ref{Item08-04.1})
schneiden sich dann je eine Tetraedergruppe $\mathcal{T} \subset P\Gamma(\mathfrak{M})$
und eine $2$-Diedergruppe $\mathcal{D}_2 \subset P\Gamma(\mathfrak{M})$ in einer Gruppe $\mathcal{C}_2$ der Ordnung $2$.
Bis auf $P\Gamma(\mathfrak{M})$-Konjugation sind $\mathcal{T}$ und $\mathcal{D}_2$ durch $\mathcal{C}_2$ eindeutig bestimmt,
und bis auf $P\Gamma(\mathfrak{M})$-Konjugation enthalten sie au{\ss}er $\mathcal{C}_2$
keine bzw. genau eine Gruppe der Ordnung $2$.
\\
Sei schlie{\ss}lich $d \not\equiv 2 \pmod 4$ oder $x^2 - dy^2 \ne 2$ f\"ur alle $x, y \in \mathbb{Z}$, und sei $d \ne 1$.
\\
Dann ist der Durchschnitt einer Tetraedergruppe $\mathcal{T} \subset P\Gamma(\mathfrak{M})$
und einer $2$-Diedergruppe $\mathcal{D}_2 \subset P\Gamma(\mathfrak{M})$ trivial,
siehe Satz \ref{Item08-03}.(\ref{Item08-03.2}).
Die Behauptung folgt mit \ref{Item08-03}.(\ref{Item08-03.2}) und \ref{Item08-04}.(\ref{Item08-04.2}).

\end{proof}


\section{Eichler-Ordnungen und Kongruenz-Untergruppen} \label{Chapter09}

\begin{definition} \label{Item09-01}

Sei $d \in \mathbb{N}$ quadratfrei, und sei $k = \mathbb{Q}(i\sqrt{d})$ mit Hauptordnung $\mathfrak{o}$.
\\
Seien $M$ eine $k$-Quaternionenalgebra, $\mathfrak{E}$ eine $M$-Ordnung,
und $\mathfrak{c} \ne \mathfrak{o}$ ein ganzes $\mathfrak{o}$-Ideal.
\begin{enumerate}[(i)]

\item \label{Item09-01.1}

Wir nennen $\mathfrak{E}$ eine $M$-Eichler-Ordnung, wenn es
$M$-Maximalordnungen $\mathfrak{M} \ne \mathfrak{M}'$ gibt mit $\mathfrak{E} = \mathfrak{M} \cap \mathfrak{M}'$.
Dann sei $\mathfrak{n}(\mathfrak{E}) := N(\mathfrak{M}\mathfrak{M}')^{-1}$.

\item \label{Item09-01.2}

Sei $\Gamma_0(\mathfrak{c}) := \left\{ \left. \begin{pmatrix} a & b \\ c & d' \end{pmatrix} \in SL_2(\mathfrak{o}) \right|
c \in \mathfrak{c} \right\}$.
\\
Sei $\Gamma_1(\mathfrak{c}) := \left\{ \left. \begin{pmatrix} a & b \\ c & d' \end{pmatrix} \in \Gamma_0(\mathfrak{c}) \right|
(a-1) \in \mathfrak{c} \right\}$.
\\
Sei $\Gamma(\mathfrak{c}) := \left\{ \left. \begin{pmatrix} a & b \\ c & d' \end{pmatrix} \in \Gamma_1(\mathfrak{c}) \right|
b \in \mathfrak{c} \right\}$.
\\
Wir bezeichnen $P\Gamma_0(\mathfrak{c})$, $P\Gamma_1(\mathfrak{c})$ und $P\Gamma(\mathfrak{c})$ als Bianchi-Kongruenzuntergruppen.

\end{enumerate}

\end{definition}

\noindent
$Bemerkung$: $\mathfrak{n}(\mathfrak{E})$ ist wohldefiniert, mit $N_{abs}(\mathfrak{n}(\mathfrak{E})) = \Lambda(\mathfrak{E})$,
siehe Lemma \ref{Item02-04} und \ref{Item05-05}.(\ref{Item05-05.1}).

\begin{lemma} \label{Item09-02}

Sei $d \in \mathbb{N}$ quadratfrei, und sei $k = \mathbb{Q}(i\sqrt{d})$.
\\
Sei $M$ eine $k$-Quaternionenalgebra, und seien $\mathfrak{M} \ne \mathfrak{M}'$ zwei $M$-Maximalordnungen.
\begin{enumerate}[(i)]

\item \label{Item09-02.1}
Sei $\mathfrak{M}''$ eine $M$-Maximalordnung mit $N(\mathfrak{M}\mathfrak{M}')^{-1} = N(\mathfrak{M}\mathfrak{M}'')^{-1}$.
\\
Dann gibt es ein $A \in \Gamma(\mathfrak{M})$ mit $A\mathfrak{M}'A^{-1} = \mathfrak{M}''$.

\item \label{Item09-02.2}
Seien $\mathfrak{M}$ und $\mathfrak{M}'$ isomorph zueinander.
\\
Dann gibt es ein $A \in M^\times$
mit $A\mathfrak{M}A^{-1} = \mathfrak{M}'$ und $A\mathfrak{M}'A^{-1} = \mathfrak{M}$.

\end{enumerate}

\end{lemma}

\begin{proof}
jeweils mit \cite[Theorem 7.7.5 (Strong Approximation)]{Ma} und $S = \{ \infty \}$ auf Basis der folgenden Aussagen.
Sei $\mathfrak{p}$ eine endliche Stelle von $k$, und sei $\pi \in \mathfrak{o}_\mathfrak{p}$ ein Primelement.

\begin{enumerate}[(i)]

\item[(\ref{Item09-02.1})]
Wir zeigen, dass es ein $A \in \Gamma(\mathfrak{M}_\mathfrak{p})$ gibt
mit $A\mathfrak{M}'_\mathfrak{p}A^{-1} = \mathfrak{M}''_\mathfrak{p}$.
Sei o.B.d.A. $\mathfrak{M}_\mathfrak{p} \ne \mathfrak{M}''_\mathfrak{p}$.
Dann zerf\"allt $M_p$ und es gibt ein $r \in \mathbb{N}$ mit
$N(\mathfrak{M}_\mathfrak{p}\mathfrak{M}''_\mathfrak{p})^{-1} = \pi^r\mathfrak{o}_\mathfrak{p}$
und Matrixeinheiten in $M_\mathfrak{p}$,
bez\"uglich derer $\mathfrak{M}_\mathfrak{p} = M_2(\mathfrak{o}_\mathfrak{p})$ gilt.
F\"ur $m, n \in \mathbb{N}_0$ und $c \in \mathfrak{o}_\mathfrak{p}$
definieren wir $J(m,n,c) := \begin{pmatrix} \pi^m & c \\ 0 & \pi^n \end{pmatrix}$.
Nach \cite[Theorem 6.5.3.3]{Ma} gibt es $m, n \in \mathbb{N}_0$ und $c \in \mathfrak{o}_\mathfrak{p}$
mit $\mathfrak{M}'_\mathfrak{p} = J(m,n,c) \mathfrak{M}_\mathfrak{p}J(m,n,c)^{-1}$
und $m + n = r$, sowie $c \in \mathfrak{o}_\mathfrak{p}^\times$, falls $mn > 0$.
F\"ur $mn > 0$ sei $A(m,n,c) := \begin{pmatrix} c^{-1}(1-\pi^n) & 1 \\ -\pi^n & c \end{pmatrix}$.
Seien $A(r,0,c) := \begin{pmatrix} 0 & -1 \\ 1 & -c \end{pmatrix}$ und $A(0,r,0) := 1$.
Mit $A' := A(m,n,c)$ gilt dann $A'\mathfrak{M}'_\mathfrak{p}A'^{-1} =
\begin{pmatrix} \mathfrak{o}_\mathfrak{p} & \pi^{-r} \mathfrak{o}_\mathfrak{p}
\\ \pi^r \mathfrak{o}_\mathfrak{p} & \mathfrak{o}_\mathfrak{p} \end{pmatrix}$.
Ebenso gibt es $A'' \in \Gamma(\mathfrak{M}_\mathfrak{p})$ mit $A''\mathfrak{M}''_\mathfrak{p}A''^{-1} =
\begin{pmatrix} \mathfrak{o}_\mathfrak{p} & \pi^{-r} \mathfrak{o}_\mathfrak{p}
\\ \pi^{-r} \mathfrak{o}_\mathfrak{p} & \mathfrak{o}_\mathfrak{p} \end{pmatrix}$.
Wir setzen $A := A''^{-1}A'$.

\item[(\ref{Item09-02.2})]
Sei $A \in M^\times$ mit $A\mathfrak{M}A^{-1} = \mathfrak{M}'$.
Wir zeigen, dass es ein $E \in \Gamma(\mathfrak{M}_\mathfrak{p})$ gibt
mit $(AE)\mathfrak{M}_\mathfrak{p}(AE)^{-1} = \mathfrak{M}'_\mathfrak{p}$
und $(AE)\mathfrak{M}'_\mathfrak{p}(AE)^{-1} = \mathfrak{M}_\mathfrak{p}$.
Sei o.B.d.A. $\mathfrak{M}_\mathfrak{p} \ne \mathfrak{M}'_\mathfrak{p}$.
Dann gibt es ein $r \in \mathbb{N}$ und Matrixeinheiten in $M_2(k_\mathfrak{p})$, bez\"uglich derer
$\mathfrak{M}_\mathfrak{p} = \begin{pmatrix} \mathfrak{o}_\mathfrak{p} & \mathfrak{o}_\mathfrak{p} \\
\mathfrak{o}_\mathfrak{p} & \mathfrak{o}_\mathfrak{p} \end{pmatrix}$,
$\mathfrak{M}'_\mathfrak{p} = \begin{pmatrix} \mathfrak{o}_\mathfrak{p} & \pi^r \mathfrak{o}_\mathfrak{p} \\
\pi^{-r}\mathfrak{o}_\mathfrak{p} & \mathfrak{o}_\mathfrak{p} \end{pmatrix}$ gilt.
Mit $A' = \begin{pmatrix} \pi^r & 0 \\ 0 & 1 \end{pmatrix}$ ist $\mathfrak{M}'_\mathfrak{p} = A'\mathfrak{M}_\mathfrak{p}A'^{-1}$
und $A = \pi^s A'B$ mit $s \in \mathbb{Z}$, $B \in \mathfrak{M}_\mathfrak{p}^\times$.
Mit $E := B^{-1}\begin{pmatrix} 0 & N(B) \\ -1 & 0 \end{pmatrix}$ gilt
$AE = \begin{pmatrix} 0 & \pi^{-s}N(A) \\ -\pi^s & 0 \end{pmatrix}$,
$E \in \Gamma(\mathfrak{M}_\mathfrak{p})$ und $S(AE) = 0$,
also $(AE)\mathfrak{M}'_\mathfrak{p}(AE)^{-1} = (AE)^2\mathfrak{M}_\mathfrak{p}(AE)^{-2} = \mathfrak{M}_\mathfrak{p}$.

\end{enumerate}

\end{proof}

\begin{theorem} \label{Item09-03}

Sei $d \in \mathbb{N}$ quadratfrei, und sei $k = \mathbb{Q}(i\sqrt{d})$ mit Hauptordnung $\mathfrak{o}$.
\\
Sei $M$ eine $k$-Quaternionenalgebra,
und seien $\mathfrak{M} \ne \mathfrak{M}'$ zwei $M$-Maximalordnungen.
\\
Dann gilt f\"ur die $M$-Eichler-Ordnung $\mathfrak{E} = \mathfrak{M} \cap \mathfrak{M}'$:
\begin{enumerate}[(i)]

\item \label{Item09-03.1}

$P\Gamma(\mathfrak{E})$ enth\"alt genau dann eine Tetraedergruppe, wenn
\\
$d \equiv 3 \pmod 8$ ist, $\mathfrak{n}(\mathfrak{E}) = 2\mathfrak{o}$ ist,
und $P\Gamma(\mathfrak{M})$ eine Tetraedergruppe enth\"alt.

\item \label{Item09-03.2}

$P\Gamma(\mathfrak{E})$ enth\"alt genau dann eine $3$-Diedergruppe, wenn
\\
$d \equiv 1 \pmod 3$ ist, $\mathfrak{n}(\mathfrak{E}) = 3\mathfrak{o}$ ist,
und $P\Gamma(\mathfrak{M})$ eine $3$-Diedergruppe enth\"alt.

\item \label{Item09-03.3}

$P\Gamma(\mathfrak{E})$ enth\"alt genau dann eine maximalendliche $2$-Diedergruppe, wenn
\\
$d \not\equiv 3 \pmod 4$ ist,
$\mathfrak{n}(\mathfrak{E}) = 2\mathfrak{o}$ oder $\mathfrak{n}(\mathfrak{E})^2 = 2\mathfrak{o}$ ist,
\\
und $P\Gamma(\mathfrak{M})$ oder $P\Gamma(\mathfrak{M}')$ eine maximalendliche $2$-Diedergruppe enth\"alt.

\end{enumerate}

\noindent
Hier gilt $M \cong M_2(k)$, und $\mathfrak{M}, \mathfrak{M}'$ sind bis auf die Reihenfolge eindeutig durch $\mathfrak{E}$ bestimmt.

\end{theorem}

\begin{proof} ${}$
\begin{enumerate}[(i)]

\item[(\ref{Item09-03.1})]

Sei $\mathcal{T} \subset P\Gamma(\mathfrak{E})$ Tetraedergruppe.
Dann gilt $\mathfrak{F}(\mathcal{T}) \subset \mathfrak{M}$ mit $\Delta(\mathfrak{F}(\mathcal{T}))$ = 2.
Nach Satz \ref{Item04-09} gibt es genau dann eine
$M$-Maximalordnung $\mathfrak{M}'' \ne \mathfrak{M}$ mit $\mathfrak{F}(\mathcal{T}) \subset \mathfrak{M}''$,
wenn $2$ in $k$ tr\"age, also $d \equiv 3 \pmod 8$ ist.
Dann ist $\mathfrak{E}' := \mathfrak{M} \cap \mathfrak{M}''$ eine $M$-Eichler-Ordnung mit
$\mathcal{T} \subset P\Gamma(\mathfrak{E}')$ und $\mathfrak{n}(\mathfrak{E}') = 2\mathfrak{o}$.
Mit Lemma \ref{Item04-06}.(\ref{Item04-06.2}) folgt $\mathfrak{M}'' = \mathfrak{M}'$ und $\mathfrak{E}' = \mathfrak{E}$.

Seien umgekehrt $\mathcal{T} \subset P\Gamma(\mathfrak{M})$ eine Tetraedergruppe,
$d \equiv 3 \pmod 8$ und $\mathfrak{n}(\mathfrak{E}) = 2\mathfrak{o}$.
Wie eben gezeigt, gibt es eine $M$-Maximalordnung $\mathfrak{M}'' \ne \mathfrak{M}$
und eine $M$-Eichler-Ordnung $\mathfrak{E}' := \mathfrak{M} \cap \mathfrak{M}''$ mit
$\mathfrak{F}(\mathcal{T}) \subset \mathfrak{E}'$
und $\mathfrak{n}(\mathfrak{E}') = N(\mathfrak{M}\mathfrak{M}'')^{-1} = 2\mathfrak{o}$.
\\
Aus $N(\mathfrak{M}\mathfrak{M}')^{-1} = 2\mathfrak{o}$
folgt mit Lemma \ref{Item09-02}.(\ref{Item09-02.1}), dass es ein $A \in \Gamma(\mathfrak{M})$ gibt mit $A\mathfrak{M}''A^{-1} = \mathfrak{M}'$.
Dann ist $\mathcal{T}' = P(AP^{-1}(\mathcal{T})A^{-1})$ Tetraedergruppe mit $\mathcal{T}' \subset P\Gamma(\mathfrak{E})$.

\item[(\ref{Item09-03.2})]

analog zum Beweis von (\ref{Item09-03.1}), mit $p = 3$ anstelle von $p = 2$.

\item[(\ref{Item09-03.3})]
Sei $\mathcal{D}_2 \subset P\Gamma(\mathfrak{E})$ maximalendliche $2$-Dieder-Untergruppe.
Dann ist $\mathcal{D}_2$ maximal\-endlich in $P\Gamma(\mathfrak{M})$ oder $P\Gamma(\mathfrak{M}')$.
Nach Satz \ref{Item06-06} gilt $d \not\equiv 3 \pmod 4$.
Mit den Bezeichnungen von Lemma \ref{Item08-01} ist
o.B.d.A. entweder $\mathfrak{E} = \mathfrak{N}(U) \cap \mathfrak{N}(V)$ oder $\mathfrak{E} = \mathfrak{N} \cap \mathfrak{N}(U)$,
wegen $\Lambda(\mathfrak{N}(U) \cap \mathfrak{N}(V)) = 4$ und $\Lambda(\mathfrak{N} \cap \mathfrak{N}(U)) = 2$
also $\mathfrak{n}(\mathfrak{E}) = 2\mathfrak{o}$ oder $\mathfrak{n}(\mathfrak{E})^2 = 2\mathfrak{o}$.

Sei umgekehrt o.B.d.A. $\mathcal{D}_2 \subset P\Gamma(\mathfrak{M})$ eine maximalendliche $2$-Dieder-Untergruppe,
$d \not\equiv 3 \pmod 4$ und $\mathfrak{n}(\mathfrak{E}) = 2\mathfrak{o}$ oder $\mathfrak{n}(\mathfrak{E})^2 = 2\mathfrak{o}$.
Mit den Bezeichnungen von Lemma \ref{Item08-01} sei $\mathfrak{M} = \mathfrak{N}(U)$, also
$N(\mathfrak{M}\mathfrak{N}(V)) = N(\mathfrak{M}\mathfrak{M}')$ oder $N(\mathfrak{M}\mathfrak{N}) = N(\mathfrak{M}\mathfrak{M}')$.
Nach Lemma \ref{Item09-02}.(\ref{Item09-02.1}) gibt es ein $A \in \Gamma(\mathfrak{M})$
mit $A\mathfrak{N}(V)A^{-1} = \mathfrak{M}'$ bzw. $A\mathfrak{N}A^{-1} = \mathfrak{M}'$.
Dann ist $\mathcal{D}'_2 = P(AP^{-1}(\mathcal{D}_2)A^{-1}) \subset P\Gamma(\mathfrak{E})$
eine maximalendliche $2$-Diedergruppe.

\end{enumerate}

\noindent
Nach Satz \ref{Item06-04} bzw. \ref{Item06-06} gilt $M \cong M_2(k)$,
und nach Lemma \ref{Item02-04} sind $\mathfrak{M}$ und $\mathfrak{M}'$ bis auf die Reihenfolge eindeutig durch $\mathfrak{E}$ bestimmt,
da sie nur an einer Stelle verschieden sind.

\end{proof}

\begin{theorem} \label{Item09-04}

Sei $d \in \mathbb{N}$ quadratfrei, und sei $k = \mathbb{Q}(i\sqrt{d})$ mit Hauptordnung $\mathfrak{o}$.
Sei $\mathfrak{c} \ne \mathfrak{o}$ ein ganzes $\mathfrak{o}$-Ideal.
Dann ist $\mathfrak{E}_0(\mathfrak{c}) := \begin{pmatrix} \mathfrak{o} & \mathfrak{o} \\ \mathfrak{c} & \mathfrak{o} \end{pmatrix}$
eine $M_2(k)$-Eichler-Ordnung mit $\mathfrak{n}(\mathfrak{E}_0(\mathfrak{c})) = \mathfrak{c}$.
Nach Definition ist $\Gamma_0(\mathfrak{c}) = \Gamma(\mathfrak{E}_0(\mathfrak{c}))$, und damit gilt:
\begin{enumerate}[(i)]

\item \label{Item09-04.1}

$P\Gamma_0(\mathfrak{c})$ enth\"alt genau dann eine Tetraedergruppe,
wenn $d \equiv 3 \pmod{8}$ ist
\\
und $\mathfrak{c} = 2\mathfrak{o}$ ist
und $p \equiv 1$ oder $p \equiv 3 \pmod{8}$ f\"ur alle Primteiler $p$ von $d$ gilt.

\item \label{Item09-04.2}

$P\Gamma_0(\mathfrak{c})$ enth\"alt genau dann eine $3$-Diedergruppe,
\\
wenn $\mathfrak{c} = 3\mathfrak{o}$ ist und $p \equiv 1 \pmod{3}$ f\"ur alle Primteiler $p$ von $d$ gilt.

\item \label{Item09-04.3}

$P\Gamma_0(\mathfrak{c})$ enth\"alt genau dann eine maximalendliche $2$-Diedergruppe,
\\
wenn $d \not\equiv 3 \pmod{4}$ ist und eine der folgenden Bedingungen erf\"ullt ist:
\begin{enumerate}[(a)]

\item \label{Item09-04.3.1}

$\mathfrak{c} = 2\mathfrak{o}$ und $p \equiv 1 \pmod{4}$ f\"ur alle Primteiler $p \not= 2$ von $d$.

\item \label{Item09-04.3.2}

$\mathfrak{c}^2 = 2\mathfrak{o}$ und $n \not\equiv 7 \pmod{8}$ f\"ur alle Teiler $n \in \mathbb{N}$ von $d$.

\end{enumerate}

\end{enumerate}

\end{theorem}

\begin{proof}

Seien $\mathfrak{M} := \begin{pmatrix} \mathfrak{o} & \mathfrak{o} \\ \mathfrak{o} & \mathfrak{o} \end{pmatrix}$,
$\mathfrak{M}' := \begin{pmatrix} \mathfrak{o} & \mathfrak{c}^{-1} \\ \mathfrak{c} & \mathfrak{o} \end{pmatrix}$
und $\mathfrak{E} := \mathfrak{M} \cap \mathfrak{M}' = \mathfrak{E}_0(\mathfrak{c})
= \begin{pmatrix} \mathfrak{o} & \mathfrak{o} \\ \mathfrak{c} & \mathfrak{o} \end{pmatrix}$.
\\
Es gilt $\mathfrak{n}(\mathfrak{E}) = N(\mathfrak{M}\mathfrak{M}')^{-1} = \mathfrak{c}$,
siehe Lemma \ref{Item02-02}.(\ref{Item02-02.1}) und (\ref{Item02-02.2}).
\begin{enumerate}[(i)]

\item[(\ref{Item09-04.1})]

folgt mit Satz \ref{Item09-03}.(\ref{Item09-03.1}) und Satz \ref{Item06-08}.(\ref{Item06-08.2}).

\item[(\ref{Item09-04.2})]

folgt mit Satz \ref{Item09-03}.(\ref{Item09-03.2}) und Satz \ref{Item06-08}.(\ref{Item06-08.1}).

\item[(\ref{Item09-04.3})]

folgt mit Satz \ref{Item09-03}.(\ref{Item09-03.3}) und  Satz \ref{Item06-08}.(\ref{Item06-08.3}).
\begin{enumerate}[(a)]

\item

$\mathfrak{M}$ und $\mathfrak{M}'$ sind isomorph zueinander.

\item

Nach Satz \ref{Item06-07}.(\ref{Item06-07.3}), (\ref{Item06-07.2}) mit Lemma \ref{Item05-05}.(\ref{Item05-05.1}) und Satz \ref{Item06-08}.(\ref{Item06-08.2})
enth\"alt $P\Gamma(\mathfrak{M}')$ maximalendliche $2$-Diedergruppen
genau dann, wenn $P\Gamma(\mathfrak{M})$ Tetraedergruppen enth\"alt,
also $p \equiv 1$ oder $p \equiv 3 \pmod{8}$ f\"ur alle Primteiler $p \not= 2$ von $d$ gilt.
\\
Die Oder-Verkn\"upfung der Bedingungen f\"ur $\mathfrak{M}$ und $\mathfrak{M}'$ liefert die Behauptung.

\end{enumerate}

\end{enumerate}

\end{proof}

\begin{theorem} \label{Item09-05}

Sei $d \in \mathbb{N}$ quadratfrei, und sei $k = \mathbb{Q}(i\sqrt{d})$ mit Hauptordnung $\mathfrak{o}$.
Sei $\mathfrak{c} \ne \mathfrak{o}$ ein ganzes $\mathfrak{o}$-Ideal.
Dann ist $\mathfrak{E}(\mathfrak{c}) := \begin{pmatrix} \mathfrak{o} & \mathfrak{c} \\ \mathfrak{c} & \mathfrak{o} \end{pmatrix}$
eine $M_2(k)$-Eichler-Ordnung mit $\mathfrak{n}(\mathfrak{E}(\mathfrak{c})) = \mathfrak{c}^2$.
\\
$P\Gamma(\mathfrak{E}(\mathfrak{c}))$, $P\Gamma_1(\mathfrak{c})$ und $P\Gamma(\mathfrak{c})$
enthalten weder Tetraeder- noch $3$-Diedergruppen.
\\
Wenn $P\Gamma(\mathfrak{E}(\mathfrak{c}))$ oder $P\Gamma_1(\mathfrak{c})$ oder $P\Gamma(\mathfrak{c})$
eine $2$-Diedergruppe enth\"alt, dann ist
\\
$d \not\equiv 3 \pmod{4}$ und $\mathfrak{c}^2 = 2\mathfrak{o}$,
sowie $\Gamma_1(\mathfrak{c}) = \Gamma_0(\mathfrak{c})$ und $\Gamma(\mathfrak{c}) = \Gamma(\mathfrak{E}(\mathfrak{c}))$.
\\
Weiter gilt dann: $P\Gamma(\mathfrak{c})$ enth\"alt genau dann eine maximalendliche $2$-Diedergruppe,
\\
wenn $p \equiv 1$ oder $p \equiv 3 \pmod{8}$ f\"ur alle Primteiler $p \not= 2$ von $d$ gilt.

\end{theorem}

\begin{proof}

Seien $\mathfrak{M} := \begin{pmatrix} \mathfrak{o} & \mathfrak{c} \\ \mathfrak{c}^{-1} & \mathfrak{o} \end{pmatrix}$
und $\mathfrak{M}' := \begin{pmatrix} \mathfrak{o} & \mathfrak{c}^{-1} \\ \mathfrak{c} & \mathfrak{o} \end{pmatrix}$.
Dann gilt $\mathfrak{E}(\mathfrak{c}) = \mathfrak{M} \cap \mathfrak{M}'$
mit $\mathfrak{n}(\mathfrak{E(\mathfrak{c})}) = N(\mathfrak{M}\mathfrak{M}')^{-1} = \mathfrak{c}^2$,
und $\Gamma(\mathfrak{c}) = \Gamma_1(\mathfrak{c}) \cap \Gamma(\mathfrak{E}(\mathfrak{c}))$.
Weiter gilt $N(\mathfrak{M} M_2(\mathfrak{o}))^{-1} = \mathfrak{c}$.
\\[4pt]
Enth\"alt $P\Gamma(\mathfrak{E}(\mathfrak{c}))$ eine Tetraeder- oder $3$-Diedergruppe,
dann ist $2$ bzw. $3$ nach Satz \ref{Item09-03} tr\"age in k,
also ist $\mathfrak{n}(\mathfrak{E}(\mathfrak{c})) = 2\mathfrak{o}$ bzw. $\mathfrak{n}(\mathfrak{E}(\mathfrak{c})) = 3\mathfrak{o}$
kein Quadrat, Widerspruch.
\\
Enth\"alt $P\Gamma_1(\mathfrak{c})$ eine Tetraeder- oder $3$-Diedergruppe,
dann gibt es $A \in P^{-1}(P\Gamma_1(\mathfrak{c}))$ mit $S(A) = 0$
und $B \in P^{-1}(P\Gamma_1(\mathfrak{c}))$ mit $S(B) = 1$.
F\"ur $X \in P^{-1}(P\Gamma_1(\mathfrak{c}))$ gilt $\pm X \in \Gamma_1(\mathfrak{c})$,
also $(S(X) \pm 2) \in \mathfrak{c}$.
Daher gilt $(S(A) \pm 2) \in \mathfrak{c}$, also $2 \in \mathfrak{c}$,
sowie $(S(B) + 2) = 3 \in \mathfrak{c}$ oder $(S(B) - 2) = 1 \in \mathfrak{c}$.
Damit folgt $3 - 2 = 1 \in \mathfrak{c}$, also $\mathfrak{c} = \mathfrak{o}$, Widerspruch.
\\[4pt]
Enth\"alt $P\Gamma(\mathfrak{E}(\mathfrak{c}))$ eine $2$-Diedergruppe,
dann folgt mit Satz \ref{Item09-03}.(\ref{Item09-03.3}), dass $d \not\equiv 3 \pmod 4$
und $\mathfrak{n}(\mathfrak{E}) = \mathfrak{c}^2 = 2\mathfrak{o}$ gilt.
Enth\"alt $P\Gamma_1(\mathfrak{c})$ eine $2$-Diedergruppe,
so folgt mit $\Gamma_1(\mathfrak{c}) \subset \Gamma_0(\mathfrak{c})$ und Satz \ref{Item09-04}.(\ref{Item09-04.3}),
dass $d \not\equiv 3 \pmod{4}$ und $\mathfrak{c} = 2\mathfrak{o}$ oder $\mathfrak{c}^2 = 2\mathfrak{o}$ gilt.
Wir f\"uhren die Annahme $\mathfrak{c} = 2\mathfrak{o}$ zum Widerspruch.
Sei also $\mathcal{D}_2 \subset P\Gamma_1(2\mathfrak{o})$ eine $2$-Diedergruppe, und sei $2\mathfrak{o} = \mathfrak{p}^2$.
Wegen $-1 \in \Gamma_1(2\mathfrak{o})$ wird $P^{-1}(\mathcal{D}_2)$ erzeugt von
$U = \begin{pmatrix} 1+2a & b \\ 2c & -1-2a \end{pmatrix}$ und
$V = \begin{pmatrix} 1+2a' & b' \\ 2c' & -1-2a' \end{pmatrix}$
mit $UV = -VU$ und $a, b, c, a', b', c' \in \mathfrak{o}$.
Mit $N(U) = 1$ folgt $2a + 2a^2 +bc + 1 = 0$.
Daher gilt $b \not\in \mathfrak{p}$ und $b' \not\in \mathfrak{p}$,
wegen $[\mathfrak{o} : \mathfrak{p}] = 2$ also $(b-b') \in \mathfrak{p}$.
Weiter gilt $2 = N(U-V) = -4(a-a')^2 - 2(b-b')(c-c')$, also $(b-b') \not \in \mathfrak{p}$, Widerspruch.
\\
Also gilt $\mathfrak{c} = \mathfrak{p}$.
F\"ur $X = \begin{pmatrix} a & b \\ c & d' \end{pmatrix} \in \Gamma_0(\mathfrak{c})$ gilt
$c \in \mathfrak{c}$ und $ad' - bc = 1$, also $ad' \not\in \mathfrak{c}$.
Wegen $[\mathfrak{o} : \mathfrak{c}] = 2$ folgt daraus $(a-1) \in \mathfrak{c}$.
Also ist $\Gamma_0(\mathfrak{c}) = \Gamma_1(\mathfrak{c})$,
und genauso $\Gamma(\mathfrak{E}(\mathfrak{c})) = \Gamma(\mathfrak{c})$.
\\[4pt]
Nach Satz \ref{Item06-07} und \ref{Item06-08}.(\ref{Item06-08.2})
enthalten $P\Gamma(\mathfrak{M}), P\Gamma(\mathfrak{M}')$ maximalendliche $2$-Diedergruppen
genau dann, wenn $P\Gamma(M_2(\mathfrak{o}))$ Tetraedergruppen enth\"alt,
also $p \equiv 1$ oder $p \equiv 3 \pmod{8}$ f\"ur alle Primteiler $p \not= 2$ von $d$ gilt.
Die letzte Aussage folgt nun mit Satz \ref{Item09-03}.(\ref{Item09-03.3}).

\end{proof}

\begin{theorem} \label{Item09-06}

Sei $d \in \mathbb{N}$ quadratfrei, und sei $k = \mathbb{Q}(i\sqrt{d})$ mit Hauptordnung $\mathfrak{o}$.
\\
Seien $D$ die Diskriminante von $k$ und $t$ die Anzahl der verschiedenen Primteiler von $D$.
\\
Sei $M$ eine $k$-Quaternionenalgebra, und sei $\mathfrak{E}$ eine $M$-Eichler-Ordnung.
Dann gilt:
\begin{enumerate}[(i)]

\item \label{Item09-06.1}

Sei $\mathcal{T} \subset P\Gamma(\mathfrak{E})$ eine  Tetraedergruppe.
Dann ist $\mu(\mathcal{T}, \mathfrak{E}) = 2^t$.

\item \label{Item09-06.2}

Sei $\mathcal{D}_3 \subset P\Gamma(\mathfrak{E})$ eine  $3$-Diedergruppe.
Dann ist $\mu(\mathcal{D}_3, \mathfrak{E}) = 2^t$.

\item \label{Item09-06.3}

Sei $\mathcal{D}_2 \subset P\Gamma(\mathfrak{E})$ maximalendliche $2$-Diedergruppe.
Dann hat $\mu(\mathcal{D}_2, \mathfrak{E})$ den Wert
\\
$2^t$, falls $\mathfrak{n}(\mathfrak{E}) = 2\mathfrak{o}$ ist oder $p \equiv \pm 1 \pmod{8}$ f\"ur alle Primteiler $p \not= 2$ von $d$ gilt
\\
$2^{t-1}$, falls $\mathfrak{n}(\mathfrak{E})^2 = 2\mathfrak{o}$ ist und $p \equiv \pm 3 \pmod{8}$ f\"ur einen Primteiler $p$ von $d$ gilt

\end{enumerate}

\end{theorem}

\begin{proof}

Sei $\mathcal{G} \subset P\Gamma(\mathfrak{E})$ eine nichttzyklische maximalendliche Untergruppe.
Nach Lemma \ref{Item07-04}.(\ref{Item07-04.1}) ist $\mu(\mathcal{G}, \mathfrak{E}) =
B_1(\mathfrak{F}(\mathcal{G}), \mathfrak{E}) [\mathcal{A}(\mathfrak{F}(\mathcal{G})) : \mathcal{I}(\mathfrak{F}(\mathcal{G}))]^{-1}$.
Da \cite[Theorem 11.6.1]{Ma} auch f\"ur Eichler-Ordnungen gilt, ist $N(\mathfrak{E}^\times) = \mathfrak{o}^\times$ und
$B_1(\mathfrak{F}(\mathcal{G}), \mathfrak{E}) = 2 B(\mathfrak{F}(\mathcal{G}), \mathfrak{E})$, siehe Satz \ref{Item07-03}.
\\[4pt]
Falls $\mathcal{G} = \mathcal{T}$ oder $\mathcal{G} = \mathcal{D}_3$, ist $\mathfrak{E}$
nach Lemma \ref{Item04-06}.(\ref{Item04-06.2}) die einzige $M$-Eichler-Ordnung mit $\mathfrak{F}(\mathcal{T}) \subset \mathfrak{E}$.
Falls $\mathcal{G} = \mathcal{D}_2$, gibt es mit den Bezeichnungen von Lemma \ref{Item08-01} genau drei $M$-Eichler-Ordnungen
$\mathfrak{E}(X) := \mathfrak{N} \cap \mathfrak{N}(X)$
mit $\mathcal{D}_2 \subset P\Gamma(\mathfrak{E}(X))$ und $\Lambda(\mathfrak{E}(X)) = 2$,
und genau drei $M$-Eichler-Ordnungen $\mathfrak{E}(X,Y) := \mathfrak{N}(X) \cap \mathfrak{N}(Y)$
mit $\mathcal{D}_2 \subset P\Gamma(\mathfrak{E}(X,Y))$ und $\Lambda(\mathfrak{E}(X,Y)) = 4$.
Konjugation mit $W$ vertauscht diese Ordnungen jeweils zyklisch.
Analog zum Beweis von Satz \ref{Item07-03} kann man jetzt zeigen, dass
$B(\mathfrak{F}(\mathcal{T}), \mathfrak{E}) = [\mathcal{A}(\mathfrak{E}) : \mathcal{I}(\mathfrak{E})]$,
$B(\mathfrak{F}(\mathcal{D}_3), \mathfrak{E}) = [\mathcal{A}(\mathfrak{E}) : \mathcal{I}(\mathfrak{E})]$ und
$B(\mathfrak{F}(\mathcal{D}_2), \mathfrak{E}) = 3[\mathcal{A}(\mathfrak{E}) : \mathcal{I}(\mathfrak{E})]$ gilt.
\\[4pt]
Seien $\mathfrak{M} \ne \mathfrak{M}'$ die beiden $M$-Maximalordnungen mit $\mathfrak{E} = \mathfrak{M} \cap \mathfrak{M}'$.
Dann induziert die Inklusion $\mathcal{A}(\mathfrak{M}) \cap \mathcal{A}(\mathfrak{M}') \hookrightarrow \mathcal{A}(\mathfrak{M})$
nach Lemma \ref{Item09-02}.(\ref{Item09-02.1}) eine surjektive Abbildung 
$\mathcal{A}(\mathfrak{M}) \cap \mathcal{A}(\mathfrak{M}') \to \mathcal{A}(\mathfrak{M})/\mathcal{I}(\mathfrak{M})$
mit Kern $\mathcal{I}(\mathfrak{M}) \cap \mathcal{A}(\mathfrak{M}')$.
Ist $j \in \mathcal{I}(\mathfrak{M}) \cap \mathcal{A}(\mathfrak{M}')$,
dann gibt es ein $J \in \mathfrak{M}^\times$ mit $j(A) = JAJ^{-1}$ f\"ur $A \in M$
und $J\mathfrak{M}'J^{-1} = \mathfrak{M}'$.
Aus $N(J) \in \mathfrak{o}^\times$ folgt mit \cite[Lemma 6.7.5]{Ma},
dass $J\mathfrak{M}' = \mathfrak{M}'J = \mathfrak{M}'$,
also dass $J \in \mathfrak{M}'^\times$ und $j \in \mathcal{I}(\mathfrak{M}')$.
Mit Lemma \ref{Item07-02} folgt
$[(\mathcal{A}(\mathfrak{M}) \cap \mathcal{A}(\mathfrak{M}')) : (\mathcal{I}(\mathfrak{M}) \cap \mathcal{I}(\mathfrak{M}'))] =
[\mathcal{A}(\mathfrak{M}) : \mathcal{I}(\mathfrak{M})] = 2^{t-1}$.
\\[4pt]
F\"ur $j \in \mathcal{A}(\mathfrak{E})$ ist
entweder $j(\mathfrak{M}) = \mathfrak{M}$ und $j(\mathfrak{M}') = \mathfrak{M}'$
oder $j(\mathfrak{M}) = \mathfrak{M}'$ und $j(\mathfrak{M}') = \mathfrak{M}$.
Ist $\mathfrak{M}$ nicht isomorph zu $\mathfrak{M}'$,
so gilt $j(\mathfrak{M}) = \mathfrak{M}$ und $j(\mathfrak{M}') = \mathfrak{M}'$ f\"ur alle $j \in \mathcal{A}(\mathfrak{E})$.
Wegen $\mathcal{I}(\mathfrak{M} \cap \mathfrak{M}') \cap \mathcal{A}(\mathfrak{M}) \cap \mathcal{A}(\mathfrak{M}') = 
\mathcal{I}(\mathfrak{M}) \cap \mathcal{I}(\mathfrak{M}')$
folgt daraus $[\mathcal{A}(\mathfrak{E}) : \mathcal{I}(\mathfrak{E})] = 2^{t-1}$.
\\
Ist $\mathfrak{M}$ isomorph zu $\mathfrak{M}'$, so gibt es nach Lemma \ref{Item09-02}.(\ref{Item09-02.2})
ein $j \in \mathcal{A}(\mathfrak{E})$, das $\mathfrak{M}$ mit $\mathfrak{M}'$ vertauscht.
Dann ist $[\mathcal{A}(\mathfrak{M} \cap \mathfrak{M}') : (\mathcal{A}(\mathfrak{M}) \cap \mathcal{A}(\mathfrak{M}'))] = 2$
und daher $[\mathcal{A}(\mathfrak{E}) : \mathcal{I}(\mathfrak{E})] = 2^t$.
\\[4pt]
Nach Satz \ref{Item09-03} und Lemma \ref{Item05-05}.(\ref{Item05-05.3})
ist $\mathfrak{M}$ genau dann nicht isomorph zu $\mathfrak{M}'$, wenn
\\
$\Lambda(\mathfrak{M} \cap \mathfrak{M}') = N_{abs}(\mathfrak{n}(\mathfrak{E})) = 2$ ist
und es eine Stelle $p$ von $\mathbb{Q}$ gibt mit $\left( \dfrac{2 ,-d} {p} \right) = -1$.
\\
Die Aussagen des Satzes folgen mit Lemma \ref{Item07-04}.(\ref{Item07-04.2}) und Auswertung der Hilbert-Symbole.

\end{proof}

\begin{theorem} \label{Item09-07}

Sei $d \in \mathbb{N}$ quadratfrei, und sei $k = \mathbb{Q}(i\sqrt{d})$.
\\
Sei $M$ eine $k$-Quaternionenalgebra, und sei $\mathfrak{E}$ eine $M$-Eichler-Ordnung.
\\
Sei $\mathcal{T} \subset P\Gamma(\mathfrak{E})$ eine Tetraedergruppe.
Sei $\mathcal{C}_2 \subset \mathcal{T}$ eine Untergruppe der Ordnung $2$.
Dann gibt es eine Tetraedergruppe $\mathcal{T}' \subset P\Gamma(\mathfrak{E})$ mit $\mathcal{C}_2  \subset \mathcal{T}'$,
die nicht $P\Gamma(\mathfrak{E})$-konjugiert zu $\mathcal{T}$ ist.
Jede Tetraedergruppe $\mathcal{T}'' \subset P\Gamma(\mathfrak{E})$ mit $\mathcal{C}_2  \subset \mathcal{T}''$
ist $P\Gamma(\mathfrak{E})$-konjugiert zu $\mathcal{T}$ oder $\mathcal{T}'$.

\end{theorem}

\begin{proof}

Nach Satz \ref{Item09-03} gibt es in $P\Gamma(\mathfrak{E})$ keine maximalendlichen $2$-Diedergruppen.
\\
Der Beweis erfolgt dann wie der Beweis von Satz \ref{Item08-03}, mit $\mathfrak{E}$ anstelle von $\mathfrak{M}$.

\end{proof}

\begin{theorem} \label{Item09-08}

Sei $d \in \mathbb{N}$ quadratfrei, und sei $k = \mathbb{Q}(i\sqrt{d})$.
\\
Sei $M$ eine $k$-Quaternionenalgebra, und sei $\mathfrak{E}$ eine $M$-Eichler-Ordnung.
\\
Sei $\mathcal{D}_3 \subset P\Gamma(\mathfrak{E})$ eine $3$-Diedergruppe
und $\mathcal{C}_3 \subset \mathcal{D}_3$ die Untergruppe der Ordnung $3$.
Dann gibt es eine $3$-Diedergruppe $\mathcal{D}'_3 \subset P\Gamma(\mathfrak{E})$ mit $\mathcal{C}_3  \subset \mathcal{D}'_3$,
die nicht $P\Gamma(\mathfrak{E})$-konjugiert zu $\mathcal{D}_3$ ist.
Jede $3$-Diedergruppe $\mathcal{D}''_3 \subset P\Gamma(\mathfrak{E})$ mit $\mathcal{C}_3  \subset \mathcal{D}''_3$
ist $P\Gamma(\mathfrak{E})$-konjugiert zu $\mathcal{D}_3$ oder $\mathcal{D}'_3$.

\end{theorem}

\begin{proof}

wie der Beweis von Satz \ref{Item08-05}, mit $\mathfrak{E}$ anstelle von $\mathfrak{M}$

\end{proof}

\begin{theorem} \label{Item09-09}

Sei $d \in \mathbb{N}$ quadratfrei, sei $d \ne 1$, und sei $k = \mathbb{Q}(i\sqrt{d})$ mit Hauptordnung $\mathfrak{o}$.
\\
Sei $M$ eine $k$-Quaternionenalgebra, und sei $\mathfrak{E}$ eine $M$-Eichler-Ordnung.
\\
Sei $\mathcal{D}_2 \subset P\Gamma(\mathfrak{E})$ eine maximalendliche $2$-Diedergruppe.
\begin{enumerate}[(i)]

\item \label{Item09-09.1}

Sei $\mathfrak{n}(\mathfrak{E})^2 = 2\mathfrak{o}$, sei $d \equiv 2 \pmod 4$, und es gebe $x, y \in \mathbb{Z}$ mit $x^2 - dy^2 = 2$.
\\
Dann gilt $p \equiv \pm 1 \pmod 8$ f\"ur alle Primteiler $p \ne 2$ von $d$, und weiter gilt dann:
\begin{enumerate}[(a)]

\item \label{Item09-09.1.1}

Es gibt eine Untergruppe $\mathcal{C}_2 \subset \mathcal{D}_2$ der Ordnung $2$ und eine
$2$-Diedergruppe $\mathcal{D}'_2 \subset P\Gamma(\mathfrak{E})$ mit $\mathcal{C}_2 \subset \mathcal{D}'_2$,
die nicht $P\Gamma(\mathfrak{E})$-konjugiert zu $\mathcal{D}_2$ ist.
\\
Jede $2$-Diedergruppe $\mathcal{D}''_2 \subset P\Gamma(\mathfrak{E})$ mit $\mathcal{C}_2  \subset \mathcal{D}_2$
\\
ist $P\Gamma(\mathfrak{E})$-konjugiert zu $\mathcal{D}_2$ oder zu $\mathcal{D}'_2$.

\item \label{Item09-09.1.2}

Seien $\mathcal{C}'_2, \mathcal{C}''_2 \subset \mathcal{D}_2$ die beiden anderen Untergruppen der Ordnung $2$.
\\
Dann sind $\mathcal{C}'_2$ und $\mathcal{C}''_2$ nicht zu $\mathcal{C}_2$, aber zueinander $P\Gamma(\mathfrak{E})$-konjugiert.
\\
Jede $2$-Diedergruppe $\mathcal{D}''_2 \subset P\Gamma(\mathfrak{E})$
mit $\mathcal{C}'_2 \subset \mathcal{D}''_2$ oder $\mathcal{C}''_2 \subset \mathcal{D}''_2$
\\
ist $P\Gamma(\mathfrak{E})$-konjugiert zu $\mathcal{D}_2$.

\end{enumerate}

\item \label{Item09-09.2}
Sei $\mathfrak{n}(\mathfrak{E}) = 2\mathfrak{o}$,
oder sei $d \not\equiv 2 \pmod 4$, oder f\"ur alle $x, y \in \mathbb{Z}$ sei $x^2 - dy^2 \ne 2$.
\begin{enumerate}[(a)]

\item \label{Item09-09.2.1}

Seien $\mathcal{C}_2, \mathcal{C}'_2, \mathcal{C}''_2 \subset \mathcal{D}_2$ die drei Untergruppen der Ordnung $2$.
\\
Dann sind $\mathcal{C}_2$, $\mathcal{C}'_2$ und $\mathcal{C}''_2$ paarweise nicht zueinander $P\Gamma(\mathfrak{E})$-konjugiert.

\item \label{Item09-09.2.2}

Es gibt eine $2$-Diedergruppe $\mathcal{D}'_2 \subset P\Gamma(\mathfrak{E})$ mit $\mathcal{C}_2  \subset \mathcal{D}'_2$,
\\
die nicht $P\Gamma(\mathfrak{E})$-konjugiert zu $\mathcal{D}_2$ ist.
\\
Jede $2$-Diedergruppe $\mathcal{D}''_2 \subset P\Gamma(\mathfrak{E})$
mit $\mathcal{C}_2  \subset \mathcal{D}''_2$, $\mathcal{C}'_2  \subset \mathcal{D}''_2$ oder $\mathcal{C}''_2  \subset \mathcal{D}''_2$
\\
ist $P\Gamma(\mathfrak{E})$-konjugiert zu $\mathcal{D}_2$ oder zu $\mathcal{D}'_2$.

\end{enumerate}

\end{enumerate}

\end{theorem}

\begin{proof}

Seien $\mathfrak{M},\mathfrak{M}'$ die $M$-Maximalordnungen mit $\mathfrak{E} = \mathfrak{M} \cap\mathfrak{M}'$,
sei $\mathcal{D}_2 \subset P\Gamma(\mathfrak{M})$ maximalendliche Untergruppe,
$\mathcal{T} \subset P\Gamma(M)$ Tetraedergruppe mit $\mathcal{D}_2 \subset \mathcal{T}$,
und $2\mathfrak{o} = \mathfrak{p}^2$.
\begin{enumerate}[(i)]

\item[(\ref{Item09-09.1})]
Sei $p \ne 2$ ein Primteiler von $d$. Dann gilt $x^2 \equiv 2 \pmod p$,
und mit dem zweiten Erg\"anzungssatz zum quadratischen Reziprozit\"atsgesetz folgt $p \equiv \pm 1 \pmod 8$.
\begin{enumerate}[(a)]

\item \label{Item09-09.01.1}

Nach Satz \ref{Item08-04}.(\ref{Item08-04.1})(\ref{Item08-04.0.1}) gibt es
eine Gruppe $\mathcal{C}_2$ der Ordnung $2$, eine $2$-Diedergruppe $\mathcal{D}'_2$ und eine Tetraedergruppe $\mathcal{T}'$
mit $\mathcal{C}_2 \subset \mathcal{D}'_2 \subset \mathcal{T} \subset P\Gamma(\mathfrak{M})$.
Die Zuordnung $\mathcal{D}_2 \mapsto \mathcal{C}_2$ ist eindeutig,
$\mathcal{C}'_2, \mathcal{C}''_2$ sind in $P\Gamma(\mathfrak{E}) \subset P\Gamma(\mathfrak{M})$ nicht konjugiert zu $\mathcal{C}_2$.
Da $P\Gamma(\mathfrak{E})$ keine Tetraedergruppe enth\"alt, erfolgt der Beweis wie der von Satz \ref{Item08-03},
mit $\mathfrak{E}$ anstelle von $\mathfrak{M}$, wenn wir die Zwischenbemerkung wie folgt ersetzen:
\\
$Zwischenbemerkung$: Wir k\"onnen $JUJ^{-1} = U$ annehmen, da $U$  in $\Gamma(\mathfrak{E})$ zu $U^{-1}$,
aber (wie eben nachgewiesen) nicht zu $V^{\pm 1}$ und  $(UV)^{\pm 1}$ konjugiert ist.

\item \label{Item09-09.01.2}

Nach Satz \ref{Item08-04}.(\ref{Item08-04.1})(\ref{Item08-04.0.2}) ist
$\mathcal{C}'_2$ in $P\Gamma(\mathfrak{M})$-konjugiert zu $\mathcal{C}''_2$.
\\
Wegen $\mathcal{D}_2 \subset \mathfrak{M}'$ und $N(\mathfrak{M}\mathfrak{M}')^{-1} = \mathfrak{p}$
folgt mit Lemma \ref{Item08-01}, dass $\mathcal{T} \subset P\Gamma(\mathfrak{M}')$.
Mit den Bezeichnungen des ersten Absatzes von Beweisabschnitt (\ref{Item08-04.0.2}) zu Satz \ref{Item08-04}
gilt daher $(k(UV) \cap \mathfrak{M}')_2 = \mathfrak{o}_2[(1+UV)/\pi]$,
und damit auch $\epsilon \in \mathfrak{M}' \cap \mathfrak{M} = \mathfrak{E}$.
\\
Wie im zweiten Absatz von Beweisabschnitt (\ref{Item08-04.0.3}) zu Satz \ref{Item08-04} folgt,
dass $\mathcal{D}''_2$ in $P\Gamma(\mathfrak{M})$ zu $\mathcal{D}_2$ oder $\mathcal{D}'_2$ konjugiert ist.
Da es nach (\ref{Item09-09.01.1}) keine Tetraedergruppe $\mathcal{T}''$
mit $\mathcal{D}''_2 \subset \mathcal{T}'' \subset P\Gamma(\mathfrak{M})$ gibt,
ist $\mathcal{D}''_2$ in $P\Gamma(\mathfrak{E}) \subset P\Gamma(\mathfrak{M})$ nicht zu $\mathcal{D}'_2$ konjugiert.

\end{enumerate}

\item[(\ref{Item09-09.2})] 
\begin{enumerate}[(a)]

\item \label{Item09-09.02.1}

Falls $d \not\equiv 2 \pmod 4$ oder $x^2 - dy^2 \ne 2$ f\"ur alle $x, y \in \mathbb{Z}$,
sind $\mathcal{C}_2$, $\mathcal{C}'_2$ und $\mathcal{C}''_2$ nach Satz \ref{Item08-04}.(\ref{Item08-04.2})(\ref{Item08-04.0.2}) 
paarweise nicht in $P\Gamma(\mathfrak{M})$, also auch nicht in $P\Gamma(\mathfrak{E})$ konjugiert.
\\
Seien also $\mathfrak{n}(\mathfrak{E}) = 2\mathfrak{o}$, $d \equiv 2 \pmod 4$, $x, y \in \mathbb{Z}$ mit $x^2 - dy^2 = 2$.
Seien $U, V, \epsilon$ wie im ersten Absatz des Beweisabschnitts ({\ref{Item08-04.0.2}}) zu Satz \ref{Item08-04}.
Es gen\"ugt, die Annahme zum Widerspruch zu f\"uhren,
dass es ein $\epsilon' \in \Gamma(\mathfrak{E})$ gibt mit $\epsilon' V \epsilon'^{-1} = U$.
Mit $\epsilon \epsilon'^{-1} U = U \epsilon \epsilon'^{-1}$ folgt
$\epsilon \epsilon'^{-1} \in k(U) \cap \Gamma(\mathfrak{M}) = \Gamma(\mathfrak{o}[U]) \subset \Gamma(\mathfrak{E})$,
also $\epsilon \in \Gamma(\mathfrak{E})$.
Da $x$ gerade ist, gilt $yi\sqrt{d} (U + V)/2 \in \mathfrak{E}$ und $(1 + UV)/\pi \in \mathfrak{E}_2 \subset \mathfrak{M}'_2$,
Widerspruch.

\item \label{Item09-09.02.2}

F\"ur $\mathcal{C}_2  \subset \mathcal{D}''_2$ erfolgt der Beweis
wie im Beweisabschnitt (\ref{Item08-03.0.1}) zu Satz \ref{Item08-03}
mit $\mathfrak{E}$ anstelle von $\mathfrak{M}$, wenn wir dort die Zwischenbemerkung wie folgt ersetzen:
\\
$Zwischenbemerkung$: Wir k\"onnen $JUJ^{-1} = U$ annehmen,,
da $U$  in $\Gamma(\mathfrak{E})$ zu $U^{-1}$,
nach Satz \ref{Item09-09}.(\ref{Item09-09.2})(\ref{Item09-09.02.1}) aber nicht zu $V^{\pm 1}$ und  $(UV)^{\pm 1}$ konjugiert ist.

F\"ur $\mathcal{C}'_2  \subset \mathcal{D}''_2$ oder $\mathcal{C}''_2  \subset \mathcal{D}''_2$
erfolgt der Beweis analog zum zweiten Absatz im Beweisabschnitt (\ref{Item08-04.0.3}) des Beweises von  Satz \ref{Item08-04},
mit $\mathfrak{E}$ anstelle von $\mathfrak{M}$.

\end{enumerate}

\end{enumerate}

\end{proof}


\emph{E-Mail:} kraemer\underline{ }norbert@t-online.de

\end{document}